\documentclass[11pt]{amsart}

\usepackage[T1]{fontenc}
\usepackage[utf8]{inputenc}
\usepackage[english]{babel}
\usepackage[leqno]{amsmath}
\usepackage{amssymb}
\usepackage{tabularx}
\usepackage{a4wide}
\usepackage{url}
\usepackage{tikz-cd}
\usepackage{csquotes}
\usepackage{enumitem}
\usepackage{mathtools}
\usepackage{manfnt}
\usepackage{mathrsfs}
\usepackage{aurical}

\definecolor{cadmiumgreen}{rgb}{0.0, 0.42, 0.24}
\definecolor{darkred}{rgb}{.85,0,0}
\usepackage[
  colorlinks,
  linkcolor=darkred,
  citecolor=cadmiumgreen,
  pdfstartview ={FitV},
]{hyperref}
\hypersetup{
  pdfauthor={Omid Amini, Matthieu Piquerez},
  pdftitle={Homology of tropical fans}}

\usepackage[
  alphabetic,
  msc-links,
  nobysame,
  lite,
]{amsrefs}

\usepackage{scalerel}

\usetikzlibrary{calc}
\tikzset{vcenter/.style={baseline={([yshift=-.8ex]current bounding box.center)}}}
\tikzset{dot/.style={insert path={node {\tikz[baseline=.6pt]\filldraw[black] (0,0) circle (1.2pt);}}}}

\setlist[itemize,1]{itemsep=\smallskipamount}
\setlist[enumerate,1]{itemsep=\smallskipamount, label=\textnormal{(\arabic*)}}

\newtheorem{thm}{Theorem}[section]
\newtheorem{lemma}[thm]{Lemma}
\newtheorem{prop}[thm]{Proposition}
\newtheorem{claim}[thm]{Claim}
\newtheorem{cor}[thm]{Corollary}

\theoremstyle{definition}
\newtheorem{defii}[thm]{Definition}
\newenvironment{defi}
  {\pushQED{\qed}\defii}
  {\popQED\enddefii}
\newtheorem{remm}[thm]{Remark}
\newenvironment{remark}
  {\pushQED{\qed}\remm}
  {\popQED\endremm}
\newtheorem{exx}[thm]{Example}
\newenvironment{example}
  {\pushQED{\qed}\exx}
  {\popQED\endexx}
\newtheorem{question}[thm]{Question}

\numberwithin{equation}{section}

\newcommand{\cf}{cf.\ }
\newcommand{\ie}{i.e.}
\newcommand{\eg}{e.g.}

\newcommand{\resp}{resp.\ }
\newcommand{\myand}{\ \textrm{and}\ }
\renewcommand{\~}{\widetilde}
\renewcommand{\hat}{\widehat}
\newcommand{\Q}{\mathbb{Q}}
\newcommand{\Z}{\mathbb{Z}}
\newcommand{\N}{\mathbb{N}}
\newcommand{\R}{\mathbb{R}}
\newcommand{\C}{\mathbb C}
\newcommand{\vsim}{\rotatebox{90}{$\sim$}}
\newcommand{\longhookrightarrow}{\lhook\joinrel\longrightarrow}

\let\oldforall\forall
\renewcommand{\forall}{\oldforall\:}
\let\oldbigwedge\bigwedge
\renewcommand{\bigwedge}{{\textstyle\oldbigwedge\!}}
\renewcommand{\emptyset}{\varnothing}
\renewcommand{\geq}{\geqslant}
\renewcommand{\leq}{\leqslant}
\renewcommand{\setminus}{\smallsetminus}

\makeatletter
\let\oldsum\sum
\renewcommand{\sum}{\@ifnextchar_\@mysum\oldsum}
\def\@mysum_#1{\oldsum_{\substack{#1}}}
\let\oldbigoplus\bigoplus
\renewcommand{\bigoplus}{\@ifnextchar_\@mybigoplus\oldbigoplus}
\def\@mybigoplus_#1{\oldbigoplus_{\substack{#1}}}
\let\oldprod\prod
\renewcommand{\prod}{\@ifnextchar_\@myprod\oldprod}
\def\@myprod_#1{\oldprod_{\substack{#1}}}
\makeatother

\makeatletter
\let\oldnu\nu
\newlength{\heightnu}
\newlength{\depthnu}
\def\nu#1_#2{{\settoheight{\heightnu}{\hbox{$#2$}}\settodepth{\depthnu}{\hbox{$#2$}}\oldnu\rule[\depthnu-3pt]{0pt}{1pt}#1_{\!#2}}}
\makeatother

\makeatletter
\let\@oldinfty\infty
\newcommand{\@sminfty}{{\scaleto{\@oldinfty}{2.8pt}}} 
\renewcommand{\infty}{{\mathchoice%
  {\displaystyle{\@oldinfty}}%
  {\textstyle{\@oldinfty}}%
  {\scriptstyle{\@sminfty}}%
  {\scriptscriptstyle{\@sminfty}}}
}
\makeatother

\newcommand{\rquot}[2]{#1\big/#2}
\newcommand{\rest}[1]{\raisebox{-1pt}{$\vert$}_{#1}}
\newcommand{\Card}[1]{\mathrm{card}(#1)} 
\newcommand{\simto}{\xrightarrow{\raisebox{-3pt}[0pt][0pt]{\small$\hspace{-1pt}\sim$}}}
\newcommand{\id}{\mathrm{id}} 
\newcommand{\dual}{\star}
\newcommand{\st}{\bigm|} 
\newcommand{\Bigst}{\Bigm|} 
\newcommand{\bul}{\bullet}

\DeclareMathOperator{\rk}{rk} 
\let\hom\relax
\DeclareMathOperator{\hom}{Hom} 
\let\Im\relax
\DeclareMathOperator{\Im}{Im} 
\DeclareMathOperator{\ord}{ord} 
\DeclareMathOperator{\coker}{coker}

\renewcommand{\i}{\mathrm i} 
\renewcommand{\d}{\mathrm d} 
\renewcommand{\k}{\Bbbk} 
\newcommand{\sing}{\mathrm{sing}} 
\newcommand{\hyp}{\mathbb H} 
\newcommand{\BM}{{^{\scaleto{\mathrm{BM}}{3.6pt}}}}
\newcommand{\cubC}{C_{_\square}} 
\newcommand{\hcubC}{C^{^\square}} 
\newcommand{\cubd}{\d_{_\square}} 
\DeclareMathOperator{\PD}{PD} 
\DeclareMathOperator{\Cone}{Cone} 
\DeclareMathOperator{\sign}{sign} 

\newcommand{\T}{\mathbb T} 
\renewcommand{\P}{\mathbb P} 
\newcommand{\CP}{{\C\P}} 
\newcommand{\TP}{{\T\P}} 

\newcommand{\SF}{\mathbf F} 
\newcommand{\RMod}{{\R_+}} 
\newcommand{\eR}{\mathbf R} 
\newcommand{\openC}{{\mathring{C}}} 
\newcommand{\conezero}{{\underline0}}
\newcommand{\e}{\mathfrak e} 
\newcommand{\nvect}{\mathfrak n} 
\newcommand{\chart}{\~} 
\newcommand{\cube}{\square} 
\DeclareMathOperator{\sed}{sed} 

\newcommand{\shiftcomp}[2][0]{{}\mkern#1mu\overline{\mkern-#1mu#2}}
\newcommand{\comp}[1]{\if#1X \shiftcomp[3]{#1}\else\if#1Z \shiftcomp[3]{#1} \else \shiftcomp{#1}\fi\fi} 

\newcommand{\suppaux}[2]{\scalebox{1}[1.4]{$#1\lvert$}#2\scalebox{1}[1.4]{$#1\rvert$}}
  \newcommand{\supp}[1]{\mathpalette\suppaux{#1}}
\newcommand{\dimsaux}[2]{\raisebox{.2ex}{\scalebox{1}[.8]{$#1\lvert$}}#2\raisebox{.2ex}{\scalebox{1}[.8]{$#1\rvert$}}}
  \newcommand{\dims}[1]{\mathpalette\dimsaux{#1}}
\newcommand{\subface}{\prec}
\newcommand{\ssubface}{\mathbin{\mathchoice
  {\subface\!\!\!\cdot}%
  {\subface\!\!\!\cdot}%
  {\subface\!\cdot}%
  {\subface\!\cdot}%
}} 
\newcommand{\supface}{\succ}
\newcommand{\ssupface}{\mathbin{\mathchoice
  {\cdot\!\!\!\supface}%
  {\cdot\!\!\!\supface}%
  {\cdot\!\supface}%
  {\cdot\!\supface}%
}}

\newcommand{\LL}{\mathscr L} 
\newcommand{\eLL}{\hat\LL} 
\newcommand{\I}{\mathscr{I}} 
\newcommand{\1}{{\mathbf 1}} 
\newcommand{\0}{{\mathbf 0}} 

\renewcommand{\div}{\mathrm{div}} 
\newcommand{\Div}{\mathrm{Div}} 
\newcommand{\Prin}{\mathrm{Prin}} 

\newcommand{\MW}{\mathrm{MW}} 
\newcommand{\W}{\mathrm{W}} 
\newcommand{\x}{\scaleto{\mathrm{X}}{5.8pt}}
\newcommand{\chow}{{\scaleto{\mathrm{Chow}}{4pt}}} 

\newcommand{\Sh}{\mathcal Sh} 
\newcommand{\Csh}{\mathscr C} 
\newcommand{\Bsh}{\mathscr B} 
\newcommand{\Ssh}{\mathscr S} 
\newcommand{\tropand}[1]{{\text{\,\Fontauri\slshape T}\scriptstyle{rop}\strut^{\hspace{-1.5ex}\scriptscriptstyle{\mathit{#1}}}}}
\newcommand{\tropf}{\tropand{\hspace{1.5ex}}} 
\newcommand{\simp}{\tropand{simp}} 
\newcommand{\unim}{\tropand{unim}} 
\newcommand{\uqproj}{\tropand{uqproj}} 

\newcommand{\tropmod}[2]{{\mathcal{T\!M}}_{\!#1}(#2)} 
\newcommand{\ctropmod}[2]{{\mathcal{CT\!M}}_{\!#1}(#2)} 
\newcommand{\basetm}[1]{#1_o} 
\newcommand{\symbuptm}{\tikz[scale=.2, baseline=-.1]{\draw(.1,-.05)to[out=50,in=-130](.9,.05) (.5,0)--++(0,.5);}}
\newcommand{\uptm}[1]{{#1}_{\!\symbuptm}} 
\newcommand{\symbcompuptm}{\tikz[scale=.2, baseline=-.1]{\filldraw(.1,-.05)to[out=50,in=-130](.9,.05) (.5,0)--++(0,.5) ++(-.2,0)rectangle++(.4,1pt);}}
\newcommand{\compuptm}[1]{{#1}_{\!\symbcompuptm}} 
\newcommand{\inftm}[1]{#1_\infty} 
\newcommand{\etm}{\e_{\!\symbuptm}} 
\newcommand{\prtm}{\mathfrak p} 

\newcommand{\U}{\mathscr U}

\newcommand{\spec}{\mathrm{Spec}} 
\newcommand{\torus}{T} 
\newcommand{\mult}{\mathfrak m} 
\DeclareMathOperator{\gys}{Gys} 
\newcommand{\an}{{\mathrm{an}}} 
\newcommand{\G}{\mathbb G} 
\newcommand{\Trop}{\mathrm{Trop}} 

\newcommand{\filt}{{\mathscr G}}
\newcommand{\F}{\mathscr F}
\renewcommand{\H}{\mathscr H}

\newcommand{\Ma}{\mathfrak M} 
\newcommand{\Fl}{\mathscr{F}} 
\newcommand{\Cl}{\mathcal C\!\ell} 
\newcommand{\crct}{\mathfrak C} 
\newcommand{\distel}{\ast} 
\newcommand{\bases}{\mathfrak B} 
\newcommand{\ind}{\mathfrak I} 
\newcommand{\rkm}{\mathrm{rk}} 
\newcommand{\contr}[1]{/#1} 
\newcommand{\del}{\setminus} 
\DeclareMathOperator{\cl}{cl} 

\newcommand{\E}{\textnormal{\textsf{E}}}
\newcommand{\Epnop}{\E}
\newcommand{\Ep}[1]{\prescript{}{#1}\Epnop} 
\newcommand{\EpI}[1]{\prescript\downarrow{#1}\Epnop} 
\newcommand{\EpII}[1]{\prescript\rightarrow{#1}\Epnop} 
\DeclareMathOperator{\Tot}{Tot} 

\let\cech\v
\renewcommand{\v}{\mathbf v} 
\newcommand{\Bl}[2]{\mathcal B\ell_{#2}(#1)} 
\newcommand{\cycl}{\mathrm{cl}} 
\newcommand{\ones}{(1, \dots, 1)}
\newcommand{\Ploc}[1]{#1^{\scaleto{\bigstar}{1ex}}} 
\newcommand{\Ext}[2]{\mathrm{Ext}\bigl(#1, #2\bigr)} 

\begin{document}

\allowdisplaybreaks

\title{Homology of tropical fans}

\author{Omid Amini}
\address{CNRS - CMLS, \'Ecole polytechnique, Institut polytechnique de Paris.}
\email{\href{omid.amini@polytechnique.edu}{omid.amini@polytechnique.edu}}

\author{Matthieu Piquerez}
\address{CMLS, \'Ecole polytechnique, Institut polytechnique de Paris.}
\email{\href{matthieu.piquerez@polytechnique.edu}{matthieu.piquerez@polytechnique.edu}}
\date{May 5, 2021}

\begin{abstract} The aim of this paper is to study homological properties of tropical fans and to propose a \emph{notion of smoothness} in tropical geometry, which goes beyond matroids and their Bergman fans and which leads to an enrichment of the category of smooth tropical varieties.

Among the resulting applications, we prove the \emph{Hodge isomorphism theorem} which asserts that the Chow rings of smooth unimodular tropical fans are isomorphic to the tropical cohomology rings of their corresponding canonical compactifications, and prove a slightly weaker statement for any unimodular fan. We furthermore introduce a \emph{notion of shellability} for tropical fans and show that shellable tropical fans are smooth and thus enjoy all the nice homological properties of smooth tropical fans. Several other interesting properties for tropical fans are shown to be shellable. Finally, we obtain a generalization, both in the tropical and in the classical setting, of the pioneering work of Feichtner-Yuzvinsky and De Concini-Procesi on the \emph{cohomology ring of wonderful compactifications} of complements of hyperplane arrangements.

\smallskip
The results in this paper form the basis for our subsequent works on Hodge theory for tropical and non-Archimedean varieties.
\end{abstract}

\maketitle

\setcounter{tocdepth}{1}

\tableofcontents

\section{Introduction} \label{sec:intro}

The work undertaken in this paper is motivated by the following fundamental questions in tropical geometry:

\begin{question}\label{question:main}
\begin{itemize}
\item What it means for a tropical variety to be smooth?
\item What are algebraic varieties which admit a smooth tropicalization?
\end{itemize}
\end{question}

As in differential topology and algebraic geometry, tropical smoothness should be a local notion. Since tropical fans and their supports form the building blocks of more general tropical varieties, the first part of the question will be reduced to the same one for tropical fans. We precise what we mean by a tropical fan.

\smallskip
Let $N \simeq \Z^n$ be a lattice of finite rank and denote by $N_\R$ the vector space generated by $N$. Recall that a fan $\Sigma$ in $N_\R$ is a non-empty collection of strongly convex polyhedral cones which verify the following two properties. First, for any cone $\sigma \in \Sigma$, any face of $\sigma$ belongs to $\Sigma$. Moreover, for two cones $\sigma$ and $\eta$ in $\Sigma$, the intersection $\sigma \cap \eta$ is a common face of both $\sigma$ and $\eta$. In particular, $\Sigma$ always contains the cone $\conezero:=\{0\}$. We denote by $\dims{\sigma}$ the dimension of a face $\sigma$. A fan is called rational if any cone in $\Sigma$ is rational in the sense that it is generated by rays which have non-zero intersection with the lattice $N$. In that case, each cone $\sigma$ in $\Sigma$ defines a sublattice $N_\sigma$ of $N$ of the same rank as the dimension of $\sigma$, and given by the integral points of the vector subspace of $N_\R$ generated by $\sigma$. A fan is pure dimensional if all its maximal cones have the same dimension. The support of a fan $\Sigma$ is denoted by $\supp\Sigma$.

A \emph{tropical fan} is a pure dimensional rational fan $\Sigma$ in $N_\R$ as above which in addition verifies the following foundational principle in tropical geometry called the balancing condition: for any cone $\tau$ of codimension one in $\Sigma$, we have the vanishing in the quotient lattice $\rquot{N}{N_\tau}$ of the following sum
\[ \sum_{\sigma \supset \tau} \e^\tau_\sigma =0 \]
where the sum is over faces $\sigma$ of maximal dimension containing $\tau$ and $\e^\tau_\sigma$ is the vector which generates the quotient $\rquot{(\sigma \cap N)}{(\tau \cap N)} \simeq \Z_{\geq 0}$. This might be regarded as a polyhedral analogue of \emph{orientability} and leads to the definition of a \emph{fundamental class} which plays an important role in the treatment of the duality and other finer geometric properties in polyhedral geometry.

A special class of tropical fans are the \emph{realizable} ones which are those arising from (faithful) tropicalizations of subvarieties of algebraic tori over trivially valued fields. Tropical fans and their supports form the local charts for more general tropical varieties. This includes in particular the realizable ones which are those coming from (faithful) tropicalizations of subvarieties of algebraic tori and more general toric varieties over non-Archimedean fields.

\smallskip
As in classical geometry, the natural properties we can expect for a tropical notion of smoothness are the following:

\begin{enumerate}
\item Tropical smoothness should be a property only of the support.
\item Smoothness for a tropical fan should be itself a local condition: local star fans which appear around each point of a smooth tropical fan should be themselves smooth.
\item Smooth tropical fans should satisfy a polyhedral analogue of the Poincaré duality.
\end{enumerate}

\medskip

The main objective of this paper is to propose a notion of smoothness in tropical geometry with the above desired list of properties. On the way to doing this, we establish several other interesting properties enjoyed by tropical fans in general, and by smooth tropical fans in particular, which we hope could be of independent interest. We also derive some applications of our results to algebraic geometry; more in connection with Hodge theory and in the direction of answering the second part of Question~\ref{question:main} will be elaborated in our forthcoming work.

\smallskip
In the remaining of this introduction, we provide a brief discussion of our results.

\subsection{Tropical smoothness} The definition of tropical smoothness that we propose in this work is homological. In this regard, it is in the spirit of the notion of homology manifold~{\renewcommand{\v}{\cech}\cites{Cech33, Lef33, Beg42, Borel57, Wilder65, BFMW96, BFMW07, Mio00, Wei02}} and the homological characterization of regular local rings~\cites{Serre55, AB56, AB57, Serre97}.

Let $\Sigma$ be a rational fan. We denote by $\SF_p$ the sheaf of $p$-multivectors on $\Sigma$ and by $\SF^p$ its dual. The sheaf $\SF_p$ is simplicial and on each cone $\sigma$ is defined by
\[\SF_p(\sigma):= \sum_{\eta \in \Sigma \\ \eta \supseteq \sigma} \bigwedge^p N_\eta \ \subseteq \bigwedge^p N, \quad \textrm{and} \quad \SF^p(\sigma) := \SF_p(\sigma)^\dual.\]
The dual $\SF^p$ (with real coefficients) can be regarded as the sheaf of tropical holomorphic $p$-forms on $\Sigma$, \cf \cite{JSS19}.

The above coefficient sheaves lead to the definition of cohomology and homology groups $H^{p,q}(\Sigma)$, $H^{p,q}_c(\Sigma)$, $H_{p,q}(\Sigma),$ and $H^\BM_{p,q}(\Sigma)$, with the last one being the tropical analogue of the Borel-Moore homology. These were introduced in~\cite{IKMZ} and further studied in~\cites{JSS19, MZ14, JRS18, GS-sheaf, GS19, AB14, ARS, Aks19, AP-tht, Mika20, AP20-hc}. We recall the relevant definitions in Section~\ref{sec:prel}.

\smallskip
Let $\Sigma$ be a tropical fan of dimension $d$. The Borel-Moore homology group $H_{d,d}^\BM(\Sigma)$ contains a canonical element $\nu_\Sigma$, the analogue of the fundamental cycle in the tropical setting.

Using the cap product $\frown$, we get a natural map
\[ \cdot\frown\nu_\Sigma\colon H^{p,q}(\Sigma) \longrightarrow H^\BM_{d-p,d-q}(\Sigma). \]

We say that a simplicial tropical fan $\Sigma$ verifies the \emph{Poincaré duality} if the above map is an isomorphism for any bidegree $(p,q)$. As we will show later, it is easy to see that this is in fact equivalent to the vanishing of $H_{a, b}^\BM(\Sigma)$ for $b\neq d$, and the \emph{surjectivity} of the natural maps
\[ \SF^p(\conezero) \to H^\BM_{d-p,d}(\Sigma), \qquad \alpha \mapsto \iota_\alpha(\nu_\Sigma),\]
defined by contraction of multivectors. We will see in Section~\ref{sec:smoothness} that this map is always injective for any tropical fan. The question of classification of tropical fans which verify the Poincaré duality is an interesting problem, and some results in this direction has been obtained by Edvard Aksnes in his master thesis~\cite{Aks19}.

\medskip

A tropical fan $\Sigma$ is called \emph{smooth} if for any cone $\sigma \in \Sigma$, the star fan $\Sigma^\sigma$ verifies the Poincaré duality. In this paper, the star fan $\Sigma^\sigma$ refers to the fan in $\rquot {N_\R}{N_{\sigma,\R}}$ induced by the cones $\eta$ in $\Sigma$ which contain $\sigma$ as a face.

\smallskip
The following is a list of properties satisfied in the category of smooth tropical fans.

\begin{enumerate}
\item Smoothness only depends on the support of the tropical fan.
\item Smoothness is a local condition, namely that, a tropical fan is smooth if it is smooth around each of its points.
\item The category of smooth tropical fans is closed under products. In fact, we have the following stronger property: Let $\Sigma_1, \Sigma_2$ be two tropical fans and set $\Sigma = \Sigma_1 \times \Sigma_2$ (which is again a tropical fan). Then $\Sigma$ is smooth if and only if $\Sigma_1 $ and $\Sigma_2$ are both smooth.
\item Any complete rational fan, \ie, a fan with support equal to $N_\R$, is tropically smooth.
\end{enumerate}

The current-in-use notion of smoothness in tropical geometry has been so far based on the notion of \emph{matroids} and their associated \emph{Bergman fans}. To any matroid $\Ma$ over a base set $E$ one associates a polyhedral fan $\Sigma_\Ma$ called the Bergman fan of the matroid in the real vector space $\rquot{\R^E\!}{\R\ones}$. We give the precise definition in Section \ref{subsec:Bergman_fans_shellable}. In the case the matroid is given by an arrangement of hyperplanes, the Bergman fan can be identified with the tropicalization of the complement of the hyperplane arrangement, for the coordinates given by the linear functions which define the hyperplane arrangement. From the calculation of homology groups of Bergman fans, proved in~\cites{Sha13a, JSS19, JRS18}, and the observation that star fans of Bergman fans are themselves Bergman, it follows that

\begin{enumerate}[resume*]
\item Bergman fans of matroids are tropically smooth in the above sense.
\end{enumerate}

This last property shows that the results we will prove in this paper in particular apply to those tropical varieties which are locally modeled by the support of Bergman fans of matroids.

The proof of the Poincaré duality for matroids is based on the use of tropical modifications and the operations of deletion and contraction on matroids, which allow to proceed by induction~\cite{Sha13a}. More generally, we will prove the following result.

\begin{enumerate}[resume*]
\item The category of smooth tropical fans is stable under tropical modifications along smooth divisors.
\end{enumerate}

The notions of divisors and tropical modifications along a divisor will be studied in detail in Sections~\ref{sec:divisors} and~\ref{sec:operations}, and lead to a notion of \emph{tropical shellability} for fans in tropical geometry discussed in Section~\ref{sec:operations-intro} below. We prove that

\begin{enumerate}[resume*]
\item Any shellable tropical fan is smooth.
\end{enumerate}

We will see that Bergman fans are shellable but there exist shellable tropical fans which are not coming from matroids. This shows that the category of smooth tropical fans introduced in this paper is strictly larger than the category of \emph{generalized Bergman fans} which are those fans having the same support as the Bergman fan of a matroid.

\smallskip
From the above properties, we deduce the following global duality theorem.

\begin{thm}[Poincaré duality for smooth tropical varieties] \label{thm:PD-intro}
Let $X$ be a smooth connected tropical variety of dimension $d$. The cap product with the fundamental class $\nu_X\in H_{d,d}^\BM(X)$ leads to the Poincaré duality with integral coefficients
\[\cdot \frown \nu_X\colon H^{p,q}(X) \to H_{d-p,d-q}^\BM(X).\]

In particular, we have $H^{d,d}_c(X, \Q) \simeq \Q$ and the natural pairing
\[H^{p,q}(X, \Q) \times H^{d-p,d-q}_c(X, \Q) \to H^{d,d}_c(X, \Q) \simeq \Q\]
is perfect.
\end{thm}

\subsection{Chow groups of unimodular fans and Hodge isomorphism theorem} Recall that for a unimodular fan $\Sigma$ in $N_\R$ with the set of rays denoted by $\Sigma_1$, the Chow ring $A^\bul(\Sigma)$ is defined as the quotient
\[A^\bul(\Sigma) := \rquot{\Z[\x_\rho \mid \rho\in \Sigma_1]}{\bigl(I + J\bigr)}\]
where
\begin{itemize}
\item $I$ is the ideal generated by the products $\x_{\rho_1}\cdots \x_{\rho_k}$, $k\in \N$, such that $\rho_1, \dots, \rho_k$ do not form a cone in $\Sigma$, and
\item $J$ is the ideal generated by the elements of the form
\[\sum_{\rho\in \Sigma_1} \langle m, \e_\rho\rangle \x_\rho\]
for $m\in M$ where $\e_\rho$ is the generator of $N\cap\rho$.
\end{itemize}

It is well-known that this can be identified with the Chow ring of the toric variety $\P_\Sigma$ associated to $\Sigma$, see~\cites{Dan78, Jur80} for complete unimodular fans and~\cites{BDP90, Bri96, FY04} for the general case.

\smallskip
We prove the following general theorem on the relation between Chow groups of unimodular fans and tropical cohomology groups of their canonical compactifications.

The \emph{canonical compactification} $\comp \Sigma$ of a simplicial fan $\Sigma$ in $\R^N$ is defined as the result of completing any cone of $\Sigma$ to a hypercube, by adding some faces at infinity. So any ray $\varrho$ becomes a segment by adding a point $\infty_\varrho$ at infinity. Any two dimensional cone $\sigma$ with two rays $\varrho_1$ and $\varrho_2$ becomes a parallelogram with vertex set the origin, $\infty_{\varrho_1}, \infty_{\varrho_2}$, and a new point $\infty_\sigma$ associated to $\sigma$. The higher dimensional cones are completed similarly, see Section \ref{sec:prel} for more details. Alternatively, $\comp \Sigma$ is obtained by taking the closure of $\Sigma$ in the partial compactification of $N_\R$ given by the tropical toric variety $\TP_\Sigma$ associated to $\Sigma$, obtained by taking the tropicalization of the toric variety $\P_\Sigma$.

\begin{thm}[Hodge isomorphism for unimodular fans]\label{thm:HI'}
Let $\Sigma$ be a unimodular fan of dimension $d$. Assume that $A^1(\Sigma^\sigma)$ is torsion-free for all $\sigma \in \Sigma$. Then, for any non-negative integer $p\leq d$, there is an isomorphism $A^p(\Sigma) \simto H^{p,p}(\comp\Sigma)$. These isomorphisms together induce an isomorphism of\/ $\Z$-algebras between the Chow ring of\/ $\Sigma$ and the tropical cohomology ring $\bigoplus_p H^{p,p}(\comp\Sigma)$. Moreover, the cohomology groups $H^{p,q}(\comp\Sigma)$ for $p < q$ or $p > q = 0$ are all vanishing. Finally, the result holds with rational coefficients without the torsion-freeness assumption.
\end{thm}

Torsion-freeness of $A^1(\Sigma)$ in the statement of the theorem is equivalent to requiring that the primitive vectors of the rays of $\Sigma$ generate the full lattice induced by $N$ in the vector space defined by $\Sigma$. If this happens for all star fans $\Sigma^\sigma$, we say that the fan $\Sigma$ is \emph{saturated}. Without this condition, the statement cannot be true. In this case, the map is always surjective but we might have a torsion kernel. As stated in the theorem, the result holds for rational coefficients without the torsion-freeness assumption. The same holds with integral coefficients if one works with the version of tropical cohomology given by Gross and Shokrieh~\cite{GS-sheaf}.

\smallskip
Note that in the above theorem, we do not assume that the fan is tropical. In the case $\Sigma$ is a smooth tropical fan, the canonical compactification $\comp\Sigma$ is a smooth tropical variety. By Poincaré duality for $\comp \Sigma$ stated in Theorem~\ref{thm:PD-intro}, and by the vanishing result in Theorem~\ref{thm:HI'}, we get vanishing of the cohomology groups $H^{p,q}(\comp \Sigma)$ for $p>q$. This leads to the following refined theorem in the smooth setting. For each $k\in \Z_{\geq 0}$, define
\[H^k(\comp \Sigma) := \bigoplus_{p+q = k} H^{p,q}(\comp \Sigma).\]

\begin{thm}[Hodge isomorphism for smooth unimodular tropical fans]\label{thm:HI-intro}
Let $\Sigma$ be a saturated unimodular tropical fan in $N_\R$. Suppose in addition that $\Sigma$ is tropically smooth. Then we get an isomorphism of rings $A^\bul(\Sigma) \simto H^{2\bul}(\comp\Sigma)$. Moreover, all the cohomology groups $H^{2\bul+1}(\comp \Sigma)$ are vanishing.
\end{thm}

As we will see later, the Hodge isomorphism property actually leads to an equivalent formulation of the notion of smoothness in tropical geometry. The precise statement is as follows.
\begin{thm}[Alternate characterization of tropical smoothness]\label{thm:smoothness_alternate-intro} A unimodular tropical fan $\Sigma$ is smooth if and only if for any face $\sigma$ of\/ $\Sigma$, the canonically compactified star fan $\comp\Sigma^\sigma$ verifies the Poincaré duality. That is, if any only if, for any $\sigma$, $H^{p,q}(\comp \Sigma^\sigma) = 0$ for $p>q$ and the ring $H^{\bul}(\comp \Sigma^\sigma)$ verifies the Poincaré duality.
\end{thm}
More generally, a tropical fan $\Sigma$ is smooth if any unimodular fan with the same support as $\Sigma$ verifies the content of Theorem~\ref{thm:smoothness_alternate-intro}. This theorem implies that the tropical smoothness studied in this paper reflects a notion of \emph{maximal degeneracy} in the polyhedral setting, consistent with the work of Deligne~\cite{Del-md}. We will elaborate on this in our forthcoming work.

\smallskip
The above results in particular provide on one side a tropical analogue of a generalized version of a theorem of Feichtner and Yuzvinsky~\cite{FY04}, which proves a similar in spirit result for some specific combinatorially defined wonderful compactifications of complements of hyperplane arrangements, and on the other side shed light on the work of Adiprasito-Huh-Katz~\cite{AHK} on Hodge theory for matroids:

\smallskip
\noindent\ -- Concerning the latter, as it was noted in~\cite{AHK}, in the case where the matroid $\Ma$ is non-realizable over any field, the Chow ring $A^\bul(\Sigma_\Ma)$ does not correspond to the Chow ring of any smooth projective variety over any field. So it came somehow as a surprise that in the non-realizable case, the Chow ring $A^\bul(\Sigma_\Ma)$, which is the Chow ring of a non-complete smooth toric variety, verifies all the nice properties enjoyed by the cohomology rings of projective complex manifolds. Our Hodge isomorphism theorem above, combined with the observation that $\comp \Sigma_\Ma$ is a smooth projective tropical variety, explains that although $\Ma$ is non-realizable, the ring $A^\bul(\Sigma_\Ma)$ is still the cohomology ring of a smooth projective variety... in the \emph{tropical world}. The results in our series of work, starting with this and its companions~\cites{AP-tht, AP20-hc}, provide generalizations of the results of~\cite{AHK} to more general classes of tropical varieties.

\medskip

\noindent\ -- Concerning the former, Feichtner and Yuzvinsky proved in~\cite{FY04} that the Chow ring of the realizable matroid $\Ma$ given by an arrangement of complex hyperplanes becomes isomorphic to the cohomology ring of a wonderful compactification of the complement of the hyperplane arrangement. Theorem~\ref{thm:generalized_FY-intro} stated below provides a generalization of this result to any smooth \emph{tropical compactification} (in the sense of Tevelev~\cite{Tev07}, we will recall this later in Section~\ref{sec:comp}) of the complement of the hyperplane arrangement, by doing away with the combinatorial data of building and nested sets, assumed usually in the framework of wonderful compactifications for getting control on cohomology. Since Theorem~\ref{thm:HI'} holds for any saturated unimodular fan and Theorem~\ref{thm:HI-intro} for any smooth saturated tropical fan, which as we mentioned form a richer class compared to generalized Bergman fans, one can expect that Theorem~\ref{thm:generalized_FY-intro} stated below should also hold for more general classes of tropical compactifications.

\smallskip
We mention here that a non-Archimedean version of the vanishing statement in Theorem \ref{thm:HI'}, with rational coefficients, is independently proved in a recent work of Ryota Mikami \cite{Mika20} where building on the work of Liu \cite{Liu20}, and using a non-Archimedean analogue of the Gersten resolution \cite{Mika20b}, he proves that the tropical cohomology groups of the Berkovich analytification $X^\an$ of a smooth variety $X$ over a trivially valued field vanish in bidegrees $(p,q)$ for $p<q$. Mikami also establishes an isomorphism between the Chow groups of $X$ and the tropical cohomology groups in Hodge bidegrees of the analytification $X^\an$. By the work of Jell~\cite{Jel19}, the cohomology groups of the Berkovich analytification $X^\an$ for projective $X$ can be obtained as the limit of its tropicalizations. In this regard, Mikami's results and our theorem stated above are complementing each other.

\subsection{Tropical Deligne resolution and cohomology of wonderful compactifications}

Let $\Sigma$ be a smooth unimodular tropical fan. By definition, for any $\sigma\in \Sigma$, the star fan $\Sigma^\sigma$ is a smooth unimodular tropical fan in $N^\sigma_\R$. We denote by $\comp \Sigma^\sigma$ its canonical compactification and as before set
\[H^k(\comp \Sigma^\sigma) := \bigoplus_{p+q=k} H^{p,q}(\comp \Sigma^\sigma).\]
By our Theorem~\ref{thm:HI-intro}, $H^k(\comp \Sigma^\sigma)$ is non-vanishing only in even degrees in which case it becomes equal to $H^{k/2,k/2}(\comp \Sigma^\sigma)$. In Section~\ref{sec:deligne} we prove the following theorem.
\begin{thm}[Tropical Deligne resolution] Let $\Sigma$ be a smooth unimodular tropical fan. Then for any non-negative integer $p$, we have the following long exact sequence
\[0 \rightarrow \SF^p(\conezero) \rightarrow \bigoplus_{\sigma \in \Sigma \\
\dims{\sigma} =p} H^0(\comp \Sigma^\sigma) \rightarrow \bigoplus_{\sigma \in \Sigma \\
\dims{\sigma} =p-1} H^2(\comp \Sigma^\sigma) \rightarrow \dots \rightarrow \bigoplus_{\sigma \in \Sigma \\
\dims{\sigma} =1} H^{2p-2}(\comp \Sigma^\sigma) \rightarrow H^{2p} (\comp \Sigma) \to 0. \]
\end{thm}
In the above sequence, the notation $\dims\sigma$ stands for the dimension of the face $\sigma$, and the maps between cohomology groups are given by the Gysin maps in tropical geometry, as we will expand later in the paper.

The above theorem could be regarded as a \emph{cohomological version of the inclusion-exclusion principle} in the sense that the cohomology groups are described with the help of coefficient sheaves and the coefficient sheaves themselves can be recovered from the cohomology groups.

The origin of the name given to the theorem comes from the Deligne spectral sequence in Hodge theory~\cite{Del-hodge2} which puts a mixed Hodge structure on the cohomology of smooth quasi-projective complex varieties. In the case of a complement of an arrangement of complex hyperplanes, by the results of Orlik-Solomon~\cite{OT13} and Zharkov~\cite{Zha13}, the cohomology in degree $p$ coincides with $\SF^p(\conezero)$ of the Bergman fan $\Sigma_\Ma$, for the matroid $\Ma$ associated to the hyperplane arrangement, and a wonderful compactification of the hyperplanes complement produces a long exact sequence as above with the cohomology groups of the strata in the compactification, see~\cite{IKMZ}. This is a consequence of a theorem of Shapiro~\cite{Shap93} which states that the complement of an arrangement of complex hyperplanes has a pure Hodge structure of Hodge-Tate type, concentrated in bidegree $(p,p)$ for the degree $p$ cohomology. The above theorem generalizes this to any smooth tropical fan.

Note that the terms appearing in the statement of our theorem above are with integral coefficients, so even in the case of a complex hyperplane arrangement, we get a refinement over the proof obtained by using Deligne spectral sequence.

\smallskip
The essence of the above stated theorem is thus the assertion that the coefficient group $\SF^p(\conezero)$ of any smooth tropical fan $\Sigma$ should be regarded as the $p$-th cohomology group of a variety $Y$ with a Hodge structure of Hodge-Tate type concentrated in bidegree $(p,p)$. This variety is precisely the support of the tropical fan $Y=\supp \Sigma$ which has cohomology $H^p(Y) = H^{p,0}(\Sigma) = \SF^p(\conezero)$ for any $p$. We will elaborate further on the connection to the Hodge theory of tropical varieties in our upcoming work.

\smallskip
We note here that a similar in spirit sequence dealing with the Stanley-Reisner rings of simplicial complexes and their quotients by (generic) linear forms play a central role in the recent works~\cites{Adi18, AY20} on partition complexes, and have found several interesting applications.

\smallskip
We provide in Section~\ref{sec:comp} an application of the above result to the cohomology of compactifications of complements of hyperplane arrangements. (A second application is discussed in the next section below.) Namely, consider a collection of hyperplanes $H_0, H_1, \dots, H_m$ in $\CP^r$ given by linear forms $\ell_0, \ell_1, \dots, \ell_m$. Assume furthermore that the vector space generated by $\ell_0,\ell_1, \dots, \ell_m$ has maximum rank, \ie, the intersection $H_0\cap H_1\cap \dots \cap H_m$ is empty. Let $X$ be the complement $\CP^r \setminus \bigcup_{j=0}^m H_j$ and consider the embedding $X \hookrightarrow \torus$ for the torus $\torus = \spec(\C[M])$ with
\[M = \ker(\deg\colon \Z^{m+1}\to \Z), \quad \deg(x_0, \dots, x_m) = x_0+ \dots + x_m.\]
The embedding is given by coordinates $[\ell_0: \dots :\ell_m]$. Denote by $N$ the dual lattice to $M$.

By the theorem of Ardila and Klivans~\cite{AK06}, the tropicalization of $X$ with respect to the above embedding coincides with the support of the Bergman fan $\Sigma_\Ma$ for the matroid $\Ma$ associated to the arrangement. We prove the following theorem.

\begin{thm}\label{thm:generalized_FY-intro}
Let $\Sigma$ be a unimodular fan with support $\Trop(X) = \supp{\Sigma_\Ma}$ and denote by $\comp X$ the corresponding compactification obtained by taking the closure of\/ $X$ in the toric complex variety $\CP_\Sigma$. The cohomology ring of\/ $\comp X$ is concentrated in even degrees and it is of Hodge-Tate type, \ie, the Hodge-decomposition in degree $2p$ is concentrated in bidegree $(p,p)$. Moreover, the restriction map
\[A^\bul(\Sigma) \simeq A^\bul(\CP_\Sigma) \longrightarrow H^{2\bul}(\comp X)\]
is an isomorphism of rings.
\end{thm}
Combined with Theorem~\ref{thm:HI-intro}, this leads to an isomorphism
\[H^{2p}(\comp \Sigma) \simeq H^{2p}(\comp X).\]

As we mentioned previously, the theorem provides a \emph{tropical} generalization of the results of Feichtner-Yuzvinsky \cite{FY04} and De Concini-Procesi \cite{DP95}, by doing away with the combinatorial data of building and nested sets, assumed usually in the theory of wonderful compactifications. In this regard, it would be interesting to extend the above theorem to the setting of arrangements of more general subvarieties, for examples those considered in~\cites{Li09, FM94, Kee93, DD15, CD17, CDDMP20, DG18, DG19, DS18}.

\subsection{Tropical shellability} \label{sec:operations-intro}

In order to show the category of smooth tropical fans is large, we introduce a tropical notion of shellability. This will be based on two types of operations on tropical fans: \emph{star subdivision} and \emph{tropical modification}. Star subdivision corresponds to the fundamental notion of blow-up in algebraic geometry. And the importance of tropical modifications in tropical algebraic geometry was pointed out by Mikhalkin~\cites{Mik06, Mik07}, and resides in the possibility of producing richer tropicalizations out of the existing ones by introducing \emph{new coordinates}. A tropical fan will be called \emph{shellable} if, broadly speaking, it can be obtained from a collection of \emph{basic tropical fans} by only using these two operations.

We give the idea of the definition here and refer to Section~\ref{sec:operations} for more precisions. Let $\Csh$ be a class of tropical fans. Let $\Bsh$ be a subset of $\Csh$ that we call the \emph{base set}. The set of \emph{tropically $\Csh$-shellable fans over $\Bsh$} denoted by $\Sh_\Csh(\Bsh)$ is defined as the smallest class of tropical fans all in $\Csh$ which verifies the following properties.

\begin{itemize}
\item (Any element in the base set is tropically $\Csh$-shellable) $\Bsh \subseteq \Sh_\Csh(\Bsh)$.
\item (Closeness under products) If $\Sigma, \Sigma' \in \Sh_\Csh(\Bsh)$ and the product $\Sigma \times \Sigma'$ belongs to $\Csh$, then $\Sigma \times \Sigma' \in \Sh_\Csh(\Bsh)$.
\item (Closeness under tropical modifications along a tropically shellable divisor) If $\Sigma \in \Sh_\Csh(\Bsh)$, then a tropical modification $\~\Sigma$ of $\Sigma$ along a divisor $\Delta$ which is in $\Sh_\Csh(\Bsh)$ belongs to $\Sh_\Csh(\Bsh)$ provided that it remains in $\Csh$.
\item (Closeness under blow-ups and blow-downs with shellable center) If $\Sigma \in \Csh$, for any cone $\sigma \in \Sigma$ whose star fan $\Sigma^\sigma$ belongs to $\Sh_\Csh(\Bsh)$ and such that the star subdivision of $\Sigma$ at $\sigma$ belongs to $\Csh$, we have $\Sigma \in \Sh_\Csh(\Bsh)$ if and only if the star subdivision belongs to $\Sh_\Csh(\Bsh)$.
\end{itemize}

The main examples of classes $\Csh$ of tropical fans which are of interest to us are \emph{all}, \resp \emph{simplicial}, \resp \emph{unimodular}, \resp \emph{quasi-projective}, \resp \emph{principal}, \resp \emph{div-faithful}, and \resp \emph{locally irreducible} tropical fans. If $\Csh$ is the class of all tropical fans, we just write $\Sh(\Bsh)$. The last three classes and their properties will be introduced and studied later in the paper.

An important example of the base set is the set $\Bsh_0=\{\{\conezero\}, \Lambda\}$ where $\Lambda$ is the complete fan in $\R$ with three cones $\conezero, \R_{\geq 0}$, and $\R_{\leq 0}$. In fact, already in this case, $\Sh(\Bsh_0)$ contains many interesting fans. For example, complete and generalized Bergman fans are all within this class, but $\Sh(\Bsh_0)$ is strictly larger.

\begin{defi}
\begin{itemize}
\item If $\Csh$ is a class of tropical fans, a tropical fan $\Sigma$ is called \emph{shellable in $\Csh$} if $\Sigma \in \Sh_\Csh(\Bsh_0)$. If $\Csh = \tropf$, we simply say that $\Sigma$ is shellable.

\item If $P$ is a predicate on tropical fans and $\Csh$ is a class of tropical fans, then $P$ is called \emph{shellable in $\Csh$} if the subclass $\Ssh$ of fans of $\Csh$ verifying $P$ verifies $\Sh_\Csh(\Ssh) = \Ssh$. \qedhere
\end{itemize}
\end{defi}

As we will show later in Sections~\ref{sec:operations}, \ref{sec:chow_ring}, and \ref{sec:homology_tropical_modification}, several nice properties of tropical fans are shellable:
\begin{enumerate}
\item Normality, local irreducibility, and div-faithfulness are all shellable.

\item Principality is shellable in the class of locally irreducible fans.

\item Poincaré duality for the Chow ring is shellable in the class of div-faithful unimodular tropical fans.
\end{enumerate}

\noindent For the algebraic geometric terminology \emph{normal, irreducible, principal, div-faithful} for tropical fans, we refer to Sections~\ref{sec:smoothness} and~\ref{sec:divisors}.

\smallskip
Most notably, we prove the following result.

\begin{thm} \label{thm:smooth_shellable_intro} Tropical smoothness is shellable.
\end{thm}
The proof of this theorem is based on results we prove on cohomology groups of tropical modifications, which we hope could be of independent interest, and which we obtain as a second application of the tropical Deligne resolution theorem.

\smallskip
We will prove in~\cite{AP-tht} the following theorem which provides a generalization of the main result of Adiprasito, Huh, and Katz~\cite{AHK} and its refinement for generalized Bergman fans, proved independently in~\cite{ADH} by Ardila, Denham, and Huh, and in an earlier version of our work~\cite{AP-tht}. (Further related results on geometry of matroids and their Bergman fans can be found in~\cites{BHMPW, BHMPW20b, BES, BEST, Eur20, DR21, Ardila18, Baker18, Cham18}.)

\begin{thm} The Hard Lefschetz property for quasi-projective unimodular tropical fans is shellable. The same statement holds for Hodge-Riemann relations.
\end{thm}

This local result furthermore leads to the development of Hodge theory for smooth projective tropical varieties, the account of which can be found in our work~\cite{AP-tht}.

\subsection{Reformulation of the tropical smoothness} After finishing the writing of this paper, we learned from Edvard Aksnes that in his forthcoming paper~\cite{Aks21} he proves the following theorem.

\begin{thm}
Let $\Sigma$ be a tropical fan of dimension $n\geq 2$. Suppose that for any face $\sigma$ of\/ $\Sigma$, the Borel-Moore homology groups $H_{a,b}^\BM(\Sigma^\sigma, \Q)$ are vanishing for all pairs $a,b$ with $b \neq \dims \sigma$.

Then $\Sigma$ satisfies the Poincaré duality if each of the star-fans $\Sigma^{\sigma}$ for $\sigma \neq \conezero$ satisfies the Poincaré duality.
\end{thm}

From the above theorem, proceeding by induction, and working with rational coefficients (instead of integral), one can deduce the following reformulation of the smoothness.
\begin{thm} A tropical fan $\Sigma$ is smooth if and only if it is smooth in codimension one and the vanishing property $H_{a,b}^\BM(\Sigma^\sigma, \Q)=0$ holds for any face $\sigma$ of\/ $\Sigma$ and any pair of integers $a,b$ with $b\neq \dims \sigma$.
\end{thm}

\subsection{Organization of the paper} In Section~\ref{sec:prel} we provide necessary background in polyhedral geometry, recall the definition of canonical compactifications and give the definitions of tropical homology and cohomology groups.

Tropical fans are introduced in Section~\ref{sec:smoothness}. That section introduces the notions of normality, (local) irreducibility, and smoothness for tropical fans and discuss basic properties enjoyed by smooth tropical fans.

Divisors on tropical fans and their characteristics are studied in Section~\ref{sec:divisors}. In that section we introduce the class of principal and div-faithful tropical fans and show that smooth tropical fans are both principal and div-faithful.

In Section~\ref{sec:operations}, we recall the definition of tropical modification along a divisor on tropical fans and introduce the above mentioned notion of shellability for tropical fans. Basic properties of shellability are studied in this section.

In Section~\ref{sec:chow_ring} we define the cycle class map from the Chow groups of a tropical fan to the tropical homology groups of the canonical compactification, and show the compatibility of tropical modifications with the cycle class map. Moreover, we show the invariance of the Chow ring with respect to tropical modifications in the class of div-faithful tropical fans. We furthermore prove the shellability of principality and div-faithfulness, and use these results to show the property for the Chow ring of a tropical fan to satisfy the Poincaré duality is shellable in the class of unimodular div-faithful tropical fans.

In Section~\ref{sec:hodge_isomorphism}, we prove the Hodge isomorphism theorem for saturated unimodular fans, and derive the equivalent formulation of smoothness stated in Theorem~\ref{thm:smoothness_alternate-intro}.

Tropical Deligne resolution for smooth tropical fans is proved in Section~\ref{sec:deligne}.

The generalization of the theorem of Feichtner and Yuzvinsky on cohomology rings of tropical compactifications of complements of hyperplane arrangement is presented in Section~\ref{sec:comp}.

Finally, in Section~\ref{sec:homology_tropical_modification} we prove Theorem~\ref{thm:smooth_shellable_intro} on shellability of tropical smoothness. The technical part is to show that smoothness is preserved by tropical modifications along smooth divisors, which we do by using the tropical Deligne resolution.

\subsection{Basic notations}\label{sec:intro-basic-notations} The set of natural numbers is denoted by $\N = \{1, 2, 3, \dots\}$. For any natural number $n$, we denote by $[n]$ the set $\{1,\dots, n\}$. For a set $E$, a subset $A \subseteq E$ and an element $a\in E$, if $a\in A$, we denote by $A-a$ the set $A \setminus \{a\}$; if $a\notin A$, we denote by $A+a$ the set $A \cup\{a\}$.

\smallskip
The set of non-negative real numbers is denoted by $\R_+$.

\smallskip
We use basic results and constructions in multilinear algebra. In particular, for a lattice $N$, we view its dual $M=N^\dual$ as (linear) forms on $N$ and refer to elements of the exterior algebras $\bigwedge^\bul N$ and $\bigwedge^\bul M$ as multivectors and multiforms, respectively. If $\ell \in M$ is a linear map on $N$ and if $\v \in \bigwedge^p N$ is a multivector, we denote by $\iota_\ell(\v) \in \bigwedge^{p-1} N$ the classical contraction of $\v$ by $\ell$. For $v_1, \dots, v_p \in N$ and $\v = v_1 \wedge \dots \wedge v_p$, this is given by
\[ \iota_\ell(v_1 \wedge \dots \wedge v_p) = \sum_{i=1}^p (-1)^{i-1} \ell(v_i) v_1 \wedge \dots \wedge \hat{v_i} \wedge \dots \wedge v_p \in \bigwedge^{p-1} N, \]
where the notation $\hat{v_i}$ indicates, as usual, that the factor $v_i$ is removed from the wedge product. We extend this definition to $k$-forms, for any $k\in \N$, by setting, recursively, $\iota_{\ell \wedge \alpha} := \iota_\alpha \circ \iota_\ell$ for any $\alpha \in \bigwedge^k M$ and $\ell\in M$.

\smallskip
Given a poset $(P, \subface)$ and a functor $\phi$ from $P$ to a category $\mathcal C$, if $\phi$ is covariant (\resp contravariant), then for a pair of elements $\tau \subface \sigma$ in $P$, we denote by $\phi_{\tau \subface \sigma}$ (\resp $\phi_{\sigma \supface \tau}$), the corresponding map $\phi(\tau) \to \phi(\sigma)$ (\resp $\phi(\sigma) \to \phi(\tau)$) in $\mathcal C$, the idea being that in the subscript of the map $\phi_{\bul}$ representing the arrow in $\mathcal C$, the first item refers to the source and the second to the target. This convention will be in particular applied to the poset of faces in a polyhedral complex.

\smallskip
In this article, unless otherwise stated, we will work with homology and cohomology with integral coefficients.

\section{Preliminaries} \label{sec:prel}

The aim of this section is to introduce basic notations and definitions which will be used all through the paper.

\smallskip
Throughout, $N$ will be a free $\Z$-module of finite rank and $M=N^\dual = \hom(N, \Z)$ will be the dual of $N$. We denote by $N_\R$ and $M_\R$ the corresponding real vector spaces, so we have $M_\R = N_\R^\dual$. For a rational polyhedral cone $\sigma$ in $N_\R$, we use the notation $N_{\sigma, \R}$ to denote the real vector subspace of $N_\R$ generated by elements of $\sigma$ and set $N^\sigma_\R := \rquot{N_\R}{N_{\sigma, \R}}$. Since $\sigma$ is rational, we get natural lattices $N^\sigma$ and $N_\sigma$ in $N^\sigma_\R$ and $N_{\sigma, \R}$, respectively, which are both of full rank.

For the ease of reading, we adopt the following convention. We use $\sigma$ (or any other face of $\Sigma$) as a superscript where referring to the quotient of some space by $N_{\sigma, \R}$ or to the elements related to this quotient. In contrast, we use $\sigma$ as a subscript for subspaces of $N_{\sigma,\R}$ or for elements associated to these subspaces.

\smallskip
We denote by $\eR := \R \cup \{\infty\}$ the extended real line with the topology induced by that of $\R$ and a basis of open neighborhoods of infinity given by intervals $(a, \infty]$ for $a\in \R$. Extending the addition of $\R$ to $\eR$ in a natural way, by setting $\infty + a = \infty$ for all $a \in \eR$, gives $\eR$ the structure of a topological monoid called the monoid of \emph{tropical numbers}. We denote by $\eR_+ := \R_+ \cup\{\infty\}$ the submonoid of non-negative tropical numbers with the induced topology. Both monoids admit a natural scalar multiplication by non-negative real numbers (setting $0\cdot\infty=0$). Moreover, the multiplication by any factor is continuous. As such, $\eR$ and $\eR_+$ can be seen as modules over the semiring $\R_+$. Recall that modules over semirings are quite similar to classical modules over rings except that, instead of being abelian groups, they are commutative monoids. Another important collection of examples of topological modules over $\R_+$ are the cones. We can naturally define the tensor product of two modules over $\R_+$.

\subsection{Fans} \label{subsec:fans}

Let $\Sigma$ be a fan of dimension $d$ in $N_\R$. The \emph{dimension} of a cone $\sigma$ in $\Sigma$ is denoted by $\dims\sigma$. The set of $k$-dimensional cones of $\Sigma$ is denoted by $\Sigma_k$, and elements of $\Sigma_1$ are called \emph{rays}. We denote by $\conezero$ the cone $\{0\}$. Any $k$-dimensional cone $\sigma$ in $\Sigma$ is determined by its set of rays in $\Sigma_1$. The \emph{support} of $\Sigma$ denoted $\supp \Sigma$ is the closed subset of $N_\R$ obtained by taking the union of the cones in $\Sigma$. A \emph{facet} of $\Sigma$ is a cone which is maximal for the inclusion. $\Sigma$ is \emph{pure dimensional} if all its facets have the same dimension. The \emph{$k$-skeleton of\/ $\Sigma$} is by definition the subfan of $\Sigma$ consisting of all the cones of dimension at most $k$, and we denote it by $\Sigma_{(k)}$.

\smallskip
Let $\Sigma$ be a fan in $N_\R$. The set of \emph{linear functions on $\Sigma$} is defined as the restriction to $\supp \Sigma$ of linear functions on $N_\R$; such a linear function is defined by an element of $M_\R$. In the case $\Sigma$ is rational, a linear function on $\Sigma$ is called \emph{integral} if it is defined by an element of $M$.

Let $f\colon \supp{\Sigma} \to \R$ be a continuous function. We say that $f$ is \emph{conewise linear on $\Sigma$} if on each face $\sigma$ of $\Sigma$, the restriction $f\rest\sigma$ of $f$ to $\sigma$ is linear. In such a case, we simply write $f\colon \Sigma \to \R$, and denote by $f_\sigma$ the linear form on $N_{\sigma,\R}$ which coincides with $f\rest\sigma$ on $\sigma$. If the linear forms $f_\sigma$ are all integral, then we say $f$ is \emph{conewise integral linear}.

\smallskip
A conewise linear function $f \colon \Sigma \to \R$ is called \emph{convex}, \resp \emph{strictly convex}, if for each face $\sigma$ of $\Sigma$, there exists a linear function $\ell$ on $\Sigma$ such that $f-\ell$ vanishes on $\sigma$ and is non-negative, \resp strictly positive, on $U \setminus \sigma$ for an open neighborhood $U$ of the relative interior of $\sigma$ in $\supp \Sigma$.

A fan $\Sigma$ is called \emph{quasi-projective} if it admits a strictly convex conewise linear function. A \emph{projective} fan is a fan which is both projective and complete. In the case $\Sigma$ is rational, it is quasi-projective, \resp projective, if and only if the toric variety $\P_\Sigma$ is quasi-projective, \resp projective.

\smallskip
A rational fan $\Sigma$ is called \emph{saturated at $\sigma$} for $\sigma\in\Sigma$ if the set of integral linear functions on $\Sigma^\sigma$ coincides with the set of linear functions which are conewise integral. This amounts in asking the lattice generated by $\Sigma^\sigma\cap N^\sigma$ to be saturated in $N^\sigma$. The fan $\Sigma$ is called \emph{saturated} if it is saturated at all its faces.

We recall that in this paper the star fan $\Sigma^\sigma$ refers to the fan in $\rquot {N_\R}{N_{\sigma,\R}}$ induced by the cones $\eta$ in $\Sigma$ which contain $\sigma$ as a face. This is consistent with the terminology used in~\cite{AHK} and differs from the one in~\cites{Kar04, BBFK02} where this is called \emph{transversal fan}.

\smallskip
For any fan $\Sigma$ of pure dimension $d$, we define the \emph{connectivity-through-codimension-one graph} of $\Sigma$ as follows. This is the graph $G=(\Sigma_d, E)$ whose vertex set is equal to the set of facets $\Sigma_{d}$ and has an edge connecting any pair of facets $\sigma$ and $\eta$ which share a codimension one face in $\Sigma$. We say $\Sigma$ is \emph{connected through codimension one} if the connectivity-through-codimension-one graph of $\Sigma$ is connected.

\subsection*{Convention} In this paper we work with fans modulo isomorphisms. For a rational fan $\Sigma$, we denote by $N_\Sigma$ the restriction of the ambient lattice to the subspace spanned by $\Sigma$. Two rational fans $\Sigma$ and $\Sigma'$ are called isomorphic if there exists an integral linear isomorphism $\phi\colon N_\Sigma\otimes\R \simto N_{\Sigma'}\otimes\R$ inducing an isomorphism between $N_\Sigma$ and $N_{\Sigma'}$ such that for each face $\sigma \in \Sigma$, $\phi(\sigma)$ is a face of $\Sigma'$ and for each face $\sigma' \in \Sigma'$, $\phi^{-1}(\sigma')$ is a face of $\Sigma$.

Any rational fan is isomorphic to a fan in the space $\R^n$ endowed with the lattice $\Z^n$ for a sufficiently big $n$. Hence we can talk about the set of isomorphism classes of rational fans. In practice, by an abuse of the language, we will make no difference between an isomorphism class of rational fans and one of its representative.

\subsection{Canonical compactification} In this section, we describe canonical compactifications of fans, and describe their combinatorics. In the case the fan is rational, the compactification is the extended tropicalization of the corresponding toric variety. A more detailed presentation of the constructions is given in \cites{AP-tht, OR11}.

\smallskip
Let $\Sigma$ be a fan in $N_\R$. We do not need to suppose for now that $\Sigma$ is rational.

\smallskip
For any cone $\sigma$, denote by $\sigma^\vee$ the \emph{dual cone} defined by
\[\sigma^\vee := \Bigl\{m \in M_\R \:\bigm|\: \langle m, a \rangle \geq 0 \:\textrm{ for all } a \in \sigma\Bigr\}. \]

\smallskip
The \emph{canonical compactification} $\comp\sigma$ of $\sigma$ (called sometimes the \emph{extended cone} of $\sigma$) is defined by
\[\comp\sigma := \sigma \otimes_{\R_+} \eR_+. \]

Alternatively, $\comp\sigma$ is given by $\hom_{\RMod}(\sigma^\vee, \eR_+)$, \ie, by the set of morphisms $\sigma^\vee \to \eR_+$ in the category of $\R_+$-modules. In both definitions, we can naturally identify $\sigma$ with the corresponding subset of $\comp\sigma$.

The topology on $\comp\sigma$ is the natural one, \ie, the finest one such that the projections
\[ \begin{array}{rcl}
  (\sigma \times \eR_+)^k &\longrightarrow& \comp\sigma \\
  (x_i, a_i)_{1\leq i\leq k} & \longmapsto & \sum_i x_i \otimes a_i,
\end{array} \]
for $k\in\N$, are continuous. This makes $\comp\sigma$ a compact topological space whose induced topology on $\sigma$ coincides with the usual one.

There is a special point $\infty_\sigma$ in $\comp\sigma$ defined by $x\otimes\infty$ for any $x$ in the relative interior of $\sigma$ (equivalently by the map $\sigma^\vee\to\eR$ which is zero on $\sigma^\perp$ and $\infty$ elsewhere). The definition does not depend on the chosen $x$. Note that for the cone $\conezero$, we have $\infty_{\conezero} = 0$.

For an inclusion of cones $\tau \subseteq \sigma$, we get an inclusion map $\comp \tau \subseteq \comp \sigma$. This identifies $\comp \tau$ as the topological closure of $\tau$ in $\comp \sigma$.

The \emph{canonical compactification} $\comp\Sigma$ is defined as the union of extended cones $\comp\sigma$, $\sigma\in \Sigma$, where for an inclusion of cones $\tau \subseteq \sigma$ in $\Sigma$, we identify $\comp\tau$ with the corresponding subspace of $\comp\sigma$. The topology of $\comp\Sigma$ is the induced quotient topology. Note that each extended cone $\comp\sigma$ naturally embeds as a subspace of $\comp\Sigma$.

\subsection{Partial compactification \texorpdfstring{$\TP_\Sigma$}{} of \texorpdfstring{$N_\R$}{the ambient space}} \label{sec:part-comp} The canonical compactification $\comp\Sigma$ of a fan $\Sigma$ naturally lives in a partial compactification $\TP_\Sigma$ of $N_\R$ defined by the fan $\Sigma$. Moreover, this partial compactification allows to enrich the conical stratification of $\comp\Sigma$ that we will define in Section~\ref{sec:conical_stratification} into an extended polyhedral structure which will be described in Section \ref{subsec:extended_polyhedral_struct}.

\smallskip
We define $\TP_\Sigma$ as follows. For any cone $\sigma$ in $\Sigma$, we consider the space $\chart\sigma$ defined as the pushout $N_\R +_\sigma \comp\sigma$ in the category of $\RMod$-modules, for the inclusions $N_\R\hookleftarrow \sigma \hookrightarrow \comp\sigma$, that we endow with the finest topology for which the sum $N_\R \times \comp\sigma \to \chart\sigma$ becomes continuous. Alternatively, one can define $\chart\sigma$ as $\hom_{\RMod}(\sigma^\vee,\eR)$. Notice that $\hom_{\RMod}(\sigma^\vee, \R)\simeq N_\R$. We set $N^\sigma_{\infty,\R}:= N_\R + \infty_\sigma \subseteq \~\sigma$. Note that with these notations we have $N^\conezero_{\infty,\R}=N_\R$.

The space $\chart\sigma$ is naturally stratified into a disjoint union of subspaces $N^\tau_{\infty,\R}$ each isomorphic to $N^\tau_\R$ for $\tau$ running over faces of $\sigma$. Moreover, the inclusions $\~\tau\subseteq\~\sigma$ for pairs of elements $\tau \subseteq \sigma$ in $\Sigma$ allow to glue these spaces and to define the space $\TP_\Sigma$.

The partial compactification $\TP_\Sigma$ of $N_\R$ is naturally stratified as the disjoint union of $N^\sigma_{\infty,\R} \simeq N^\sigma_\R$ for $\sigma \in \Sigma$.

\smallskip
We have a natural inclusion of $\comp\sigma$ into $\chart\sigma$. We thus get an embedding $\comp\Sigma \subseteq \TP_\Sigma$ which identifies $\comp \Sigma$ as the closure of $\Sigma$ in $\TP_\Sigma$.

\subsection{Conical stratification of \texorpdfstring{$\comp\Sigma$}{the canonical compactification}} \label{sec:conical_stratification} The canonical compactification admits a natural stratification into cones that we will enrich later into an extended polyhedral structure.

\smallskip
Consider a cone $\sigma \in \Sigma$ and a face $\tau$ of $\sigma$. Let $C^\tau_\sigma$ be the subset of $\comp\sigma$ defined by
\[C^\tau_\sigma := \{\infty_\tau + x \mid x\in \sigma\}.\]
Under the natural isomorphism $N^\tau_{\infty,\R}\simeq N^\tau_\R$, one can see that $C^\tau_\sigma$ becomes isomorphic to the projection of the cone $\sigma$ into the linear space $N^\tau_\R \simeq N^\tau_{\infty,\R}$; the isomorphism is given by adding $\infty_\tau$. We denote by $\openC^\tau_\sigma$ the relative interior of $C^\tau_\sigma$. The following proposition gives a precise description of how these different sets are positioned together in the canonical compactification.

\begin{prop}\label{prop:conical_stratification} Let $\Sigma$ be a fan in $N_\R$.
\begin{itemize}
\item The canonical compactification $\comp\Sigma$ is a disjoint union of (open) cones $\openC^\tau_\sigma$ for pairs $(\tau, \sigma)$ of elements of $\Sigma$ with $\tau \subseteq \sigma$. The linear span of the cone $C^\tau_\sigma$ is the real vector space $N^\tau_{\infty, \sigma,\R}$, \ie, the projection of $N_{\sigma, \R}$ into $N^\tau_{\infty, \R}$.
\item For any pair $(\tau, \sigma)$, the closure $\comp C^\tau_\sigma$ of $C^\tau_\sigma$ in $\comp\Sigma$ is the union of all the (open) cones $\openC^{\tau'}_{\sigma'}$ with $\tau \subseteq \tau' \subseteq \sigma' \subseteq \sigma$.
\end{itemize}
\end{prop}

Note that in order to lighten the notations, we write $\comp C^\tau_\sigma$ for the compactification of $C^\tau_\sigma$ instead of the more correct form $\comp{C^\tau_\sigma}$. We use this kind of simplifications all through the paper.

\begin{proof} The proof is a consequence of the tropical orbit-stratum correspondence theorem in the tropical toric variety $\TP_\Sigma$ and the observation we made previously that $\comp \Sigma$ is the closure of $\Sigma$ in $\TP_\Sigma$. We omit the details.
\end{proof}

The cones $\openC^\tau_\sigma$ form the \emph{open strata} of what we call the \emph{conical stratification} of $\comp\Sigma$. We refer to the topological closures $\comp C^\tau_\sigma$ as the \emph{closed strata} of $\comp\Sigma$.

\smallskip
In the case $\Sigma$ is simplicial, closed stratum of dimension $k$ are isomorphic to the hypercube $\eR_+^k$. In this case, we prefer to use the notation $\cube^\tau_\sigma$ instead of $\comp C^\tau_\sigma$.

\subsection{Stratification of \texorpdfstring{$\comp\Sigma$}{the canonical compactification} into fans} There is a second stratification of $\comp\Sigma$ into fans that we describe now.

\begin{defi}[Fans at infinity]
Let $\tau$ be a cone in $\Sigma$. The \emph{fan at infinity based at $\infty_\tau$} denoted by $\Sigma^\tau_\infty$ is the fan in $N^\tau_{\infty,\R}$ which consists of all the cones $C^\tau_\sigma$ for $\sigma$ in $\Sigma$ with $\sigma \supseteq \tau$. Note that we have $\Sigma^\conezero_\infty = \Sigma$, and more generally, $\Sigma^\tau_\infty$ is isomorphic to the star fan $\Sigma^\tau$ of $\tau$.
\end{defi}

From the above descriptions, we get the following proposition.

\begin{prop} The collection of fans at infinity $\Sigma_\infty^\tau$, $\tau \in \Sigma$, provides a partition of\/ $\comp\Sigma$ into locally closed subspaces. The closure of\/ $\Sigma_\infty^\tau$ in $\comp\Sigma$ is naturally isomorphic to the canonical compactification $\comp\Sigma^\tau$\! of the star fan $\Sigma^\tau$.
\end{prop}
\begin{proof} The proof follows directly from Proposition~\ref{prop:conical_stratification}.
\end{proof}

\subsection{Combinatorics of the conical and fan stratifications of \texorpdfstring{$\comp\Sigma$}{the canonical compactification}} For a fan $\Sigma$, we denote by $\LL_\Sigma$ the face poset of $\Sigma$ in which the partial order $\subface$ is given by the inclusion of faces: we write $\tau \subface \sigma$ if $\tau \subseteq \sigma$. Set $\0 := \conezero$. The \emph{extended poset $\eLL_\Sigma$} is defined as $\eLL_\Sigma := \LL_\Sigma \sqcup\{\1\}$ obtained by adding an element $\1$ and extending the partial order to $\eLL_\Sigma$ by declaring $\sigma \subface \1$ for all $\sigma \in \LL_\Sigma$.

The \emph{join} and \emph{meet} operations $\vee$ and $\wedge$ on $\LL_\Sigma$ are defined as follows. For two cones $\sigma$ and $\tau$ of $\Sigma$, we set $\sigma \wedge \tau := \sigma \cap \tau$. To define the operation $\vee$, note that the set of cones in $\Sigma$ which contain both $\sigma$ and $\tau$ is either empty or has a minimal element $\eta \in \Sigma$. In the former case, we set $\sigma \vee \tau := \1$, and in the latter case, $\sigma \vee \tau := \eta$. The two operations are extended to the augmented poset $\eLL_\Sigma$ by $\sigma \wedge \1 = \sigma$ and $\sigma \vee \1 = \1$ for any cone $\sigma$ of $\Sigma$.

\subsection*{Notations} The above discussion leads to the following notations. Let $\tau$ and $\sigma$ be a pair of faces in $\Sigma$. We say $\sigma$ \emph{covers} $\tau$ and write $\tau \ssubface \sigma$ if $\tau\subface\sigma$ and $\dims{\tau} = \dims\sigma-1$. A family of faces $\sigma_1,\dots,\sigma_k$ are called \emph{comparable} if $\sigma_1 \vee \cdots \vee \sigma_k \neq \1$. Moreover, we use the notation $\sigma\sim\sigma'$ for two faces $\sigma$ and $\sigma'$ if $\sigma\wedge\sigma' = \0$ and $\sigma\vee\sigma' \neq \1$.

\medskip

The following proposition is straightforward.
\begin{prop} The poset $(\eLL_\Sigma, \subface)$ with the operations $\wedge$ and $\vee$, and minimal and maximal elements $\0$ and $\1$, respectively, is a lattice.
\end{prop}
The term lattice in the proposition is understood in the sense of order theory~\cite{MMT}.

\smallskip
For a poset $(\LL, \leq)$, and two elements $a,b$ in $\LL$, the (\emph{closed}) \emph{interval} $[a, b]$ is the set of all elements $c$ with $a \leq c \leq b$. By an abuse of the notation, we denote the empty interval by $\0$. Note in particular that the interval $[\0,\0] =\{\0\}$ is different from $\0$. We denote by $\I(\LL)$ the \emph{poset of intervals} of $\LL$ ordered by inclusion.

\smallskip
For a lattice $\eLL$, the poset of intervals $\I(\eLL)$ becomes a lattice for the join and meet operations $\vee$ and $\wedge$ on the intervals defined by
\begin{itemize}
\item $[a_1, b_1] \wedge [a_2, b_2] := [a_1\vee a_2, b_1 \wedge b_2] $, and
\item $[a_1, b_1] \vee [a_2, b_2] := [a_1 \wedge a_2, b_1\vee b_2]$.
\end{itemize}

\smallskip
We have the following characterization of the combinatorics of open and closed strata in the above two stratifications of the canonical compactification of a fan.

\begin{thm} Let $\Sigma$ be a fan in $N_\R$ and consider its canonical compactification $\comp\Sigma$.
\begin{enumerate}
\item There is a bijection between the open strata in the conical stratification of $\comp\Sigma$ and the non-empty elements of the interval lattice $\I(\LL_\Sigma)$. Under the above bijection, the poset of closed strata of $\comp\Sigma$ (for the partial order given by inclusion) becomes isomorphic to the poset $\I(\LL_\Sigma) -\0$ of non-empty intervals in the lattice $\LL_\Sigma$.
\item Consider now the fan stratification of $\comp\Sigma$. The induced poset by closed strata $\comp\Sigma^\tau$ in this stratification is isomorphic to the opposite poset of $\LL_\Sigma$.
\end{enumerate}
\end{thm}
\begin{proof} Both the claims can be proved by straightforward verification.
\end{proof}

\begin{question} Characterize the interval lattices associate to the face lattices of polyhedra. In particular, is it true that the interval lattice of the face lattice of a polyhedron is itself the face lattice of a polyhedron?
\end{question}

\subsection{Extended polyhedral structure} \label{subsec:extended_polyhedral_struct} The closed strata $\comp C^\tau_\sigma$, for $\tau\subface\sigma$ a pair of faces of $\Sigma$, endow $\comp\Sigma$ with an \emph{extended polyhedral structure}. We give a brief description of these structures here and refer to~\cites{JSS19, MZ14, IKMZ, ACP, AP-tht} for more details.

\smallskip
Extended polyhedral structures in \cites{JSS19, MZ14} are called polyhedral spaces with a face structure. In our context, the \emph{faces} of $\comp\Sigma$ are the closed strata in the conical stratification we defined above. These faces verify the same axioms as those of polyhedral complexes with the difference that each face is now an \emph{extended polyhedron}. An extended polyhedron is the closure of a polyhedron in the partial compactification of the ambient vector space induced by this polyhedron. For instance, using the notations of Section~\ref{sec:part-comp}, $\comp C^\tau_\sigma$ is the closure of the polyhedron $C^\tau_\sigma\subseteq N^\tau_{\infty,\R}$ in the associated partial compactification $\~C^\tau_\sigma$ of $N^\tau_{\infty,\R}$, where the cone $C^\tau_\sigma$ is seen as a face of $\Sigma^\tau_\infty \simeq \Sigma^\tau$.

In the rest of this article, all canonical compactifications $\comp\Sigma$ are considered as extended polyhedral structures. We extend to $\comp\Sigma$ the notations introduced for simplicial complexes. In particular, $\delta \in \comp\Sigma$ means that $\delta$ is a face of $\comp\Sigma$, $\supp{\comp\Sigma}$ denotes the support of $\comp\Sigma$, and $\comp\Sigma_k$ denotes the set of faces of dimension $k$ in $\comp\Sigma$.

\smallskip
In the setting considered in this paper, the \emph{sedentarity} of a face $\delta = \comp C^\tau_\sigma$ denoted by $\sed(\delta)$ is by definition the element $\tau$ of $\Sigma$, and refers to the point $\infty_\tau$ at which the stratum $C^\tau_\sigma$ is based. By an abuse of the terminology, remembering only the dimension of $\sed(\delta)$, we might say sometimes that $\delta$ has sedentarity $k$ where $k$ is the dimension of $\tau$.

\subsection{Canonical elements, unit normal vectors and orientations} \label{subsec:orientation}

Let now $\Sigma$ be a rational fan of pure dimension $d$. Let $\sigma$ be a cone of $\Sigma$ and let $\tau$ be a face of codimension one in $\sigma$. Then, $N_{\tau,\R}$ cuts $N_{\sigma,\R}$ into two closed half-spaces only one of which contains $\sigma$. Denote this half-space by $H_\sigma$. By a \emph{unit normal vector to $\tau$ in $\sigma$} we mean any vector $v$ of $N_\sigma \cap H_\sigma$ such that $N_\tau + \Z v = N_\sigma$. We usually denote such an element by $\nvect_{\sigma/\tau}$ and note that it induces a well-defined vector in $N^\tau_\sigma = \rquot{N_\sigma}{N_\tau}$ that we denote by $\e^\tau_\sigma$. We naturally extend the definition to similar pair of faces in $\comp\Sigma$ having the same sedentarity. In the case $\sigma$ is a ray (and $\tau$ is a point of the same sedentarity as $\sigma$), we also use the notation $\e_\sigma$ instead of $\nvect_{\sigma/\tau}$.

On each face $\sigma$ of $\Sigma$, we choose a generator of $\bigwedge^{\dims{\sigma}}N_\sigma$ denoted by $\nu_\sigma$. We call this element the \emph{canonical multivector of $\sigma$}. The dual notion in $\bigwedge^{\dims{\sigma}}N_\sigma^\dual$ is denoted $\omega_\sigma$ and is called the \emph{canonical form on $\sigma$}. In the rest of the article, we assume that $\nu_\conezero = 1$ and that $\nu_\rho = \e_\rho$ for any ray $\rho$.

From the above discussion, we obtain in particular \emph{the canonical element} of $\Sigma$ that we denote by $\nu_\Sigma$ as
\[ \nu_\Sigma := (\nu_\eta)_{\eta \in \Sigma_d} \in \bigoplus_{\eta \in \Sigma_d} \bigwedge^d N_\eta. \]
As we will see later, in the case of tropical fans, this defines an element of the Borel-Moore homology $H^\BM_{d,d}(\Sigma)$ (defined in Section \ref{subsec:homology} below) that we call the \emph{fundamental class of $\Sigma$}.

We extend the above definitions to $\comp\Sigma$ by choosing an element $\nu^\tau_\sigma \in \bigwedge^{\dims{\sigma}-\dims{\tau}}N^\tau_\sigma$ for each polyhedron $\comp C^\tau_\sigma$ (as we see below, for simplicial $\Sigma$, there is a natural choice).

These choices define an \emph{orientation} on $\comp\Sigma$. In particular, to any pair of closed faces $\gamma \ssubface \delta$ of $\comp\Sigma$, we associate a \emph{sign} defined as follows. If both faces have the same sedentarity, then $\sign(\gamma,\delta)$ is the sign of $\omega_\delta\bigl(\nvect_{\delta/\gamma} \wedge \nu_\gamma\bigr)$. Otherwise, there exists a pair of cones $\tau' \ssubface \tau$ and a cone $\sigma$ in $\Sigma$ such that $\gamma = \comp C^{\tau}_\sigma$ and $\delta = \comp C^{\tau'}_\sigma$. In this case, $\sign(\gamma,\delta)$ is the sign of $-\omega_\delta(\e^{\tau'}_\tau \wedge \pi_{\tau'\ssubface\tau}^*(\nu_\gamma))$ where $\pi_{\tau'\ssubface\tau}\colon N_\infty^{\tau'} \to N_\infty^\tau$ is the natural projection and $\e^{\tau'}_\tau$ refers, as above, to the primitive vector of the ray $C^{\tau'}_\tau$ in $C^{\tau'}_\sigma$.

\smallskip
In the case $\Sigma$ is simplicial, there is a natural choice for extending the orientation of $\Sigma$ to $\comp\Sigma$. If $\tau \subface \sigma$ is a pair of cones, there exists a unique minimal $\tau'$ such that $\sigma = \tau \vee \tau'$. Then we set $\nu^\tau_\sigma$ to be the image of $\nu_{\tau'}$ via the projection $N_\R \to N_{\infty, \R}^\tau\simeq N^\tau_\R$.

\subsection{Tropical homology and cohomology groups} \label{subsec:homology} Let $\Sigma$ be a rational fan. The extended polyhedral structure on $\comp\Sigma$ leads to the definition of tropical homology groups and cohomology groups introduced in~\cite{IKMZ} and further studied in~\cites{JSS19, MZ14, JRS18, GS-sheaf, GS19, AB14, ARS, Aks19, AP-tht, Mika20, AP20-hc}.

We recall the definition of the multi-tangent and multi-cotangent (integral) spaces $\SF_p$ and $\SF^p$ for $\comp\Sigma$ for the face structure given by the closed strata $\comp C^\tau_\sigma$ in the conical stratification of $\comp\Sigma$, for pairs of faces $\tau \subface \sigma$ in $\Sigma$, and give a combinatorial complex which calculates the tropical homology and cohomology groups of $\comp \Sigma$. The definitions are naturally adapted to faces in $\Sigma$.

For any face $\delta = \comp C^\tau_\sigma$, we set $N_\delta = N^\tau_{\infty, \sigma} \simeq \rquot{N_{\sigma,\R}}{N_{\tau, \R}}$. For any non-negative integer $p$, the \emph{$p$-th multi-tangent} and the \emph{$p$-th multicotangent space of\/ $\comp\Sigma$ at $\delta$} denoted by $\SF_p(\delta)$ and $\SF^p(\delta)$, respectively, are given by
\[\SF_p(\delta):= \hspace{-.5cm} \sum_{\eta \supface \delta \\ \sed(\eta) = \sed(\delta) } \hspace{-.5cm} \bigwedge^p N_\eta \ \subseteq \bigwedge^p N^\tau, \quad \textrm{and} \quad \SF^p(\delta) := \SF_p(\delta)^\dual.\]

For an inclusion of faces $\gamma \subface \delta$ in $\comp\Sigma$, we get maps $\i_{\delta \supface \gamma}\colon \SF_p(\delta) \to \SF_p(\gamma)$ and $\i^*_{\gamma \subface \delta}\colon \SF^p(\gamma) \to \SF^p(\delta)$ defined as follows. If $\gamma$ and $\delta$ have the same sedentarity, the map $\i_{\delta \supface \gamma}$ is just an inclusion. If $\gamma=\comp C_\sigma^{\tau'}$ and $\delta=\comp C_\sigma^\tau$ with $\tau \subface \tau' \subface \sigma$, then the map $\i_{\delta\supface \gamma}$ is induced by the projection $N^\tau_{\infty,\R} \to N^{\tau'}_{\infty,\R}$. In the general case, $\i_{\delta\supface \gamma}$ is given by the composition of the projection and the inclusion; the map $\i^*_{\gamma \subface \delta}$ is the dual of $\i_{\delta\supface \gamma}$.

\medskip

Let $X$ be a fan or its compactification. For a pair of non-negative integers $p,q$, define
\[ C_{p,q}(X) := \bigoplus_{\delta \in X \\ \dims\delta=q \\ \delta\text{ compact}} \SF_p(\delta) \]
and consider the corresponding complexes
\[C_{p,\bullet}\colon \quad \dots\longrightarrow C_{p, q}(X) \xrightarrow{\ \partial_q\ } C_{p,q-1}(X) \xrightarrow{\ \partial_{q-1}\ } C_{p,q-2} (X)\longrightarrow\cdots\]
where the differential is given by the sum of maps $\sign(\gamma,\delta)\cdot\i_{\delta\supface \gamma}$ with the signs corresponding to a chosen orientation on $X$ as explained in Section \ref{subsec:orientation}.

\smallskip
The \emph{tropical homology} of $X$ is defined by
\[ H_{p,q}(X) := H_q(C_{p,\bullet}). \]

Similarly, we have a cochain complex
\[C^{p,\bullet}\colon \quad \dots\longrightarrow C^{p, q-2}(X) \xrightarrow{\ \d^{q-2}\ } C^{p,q-1}(X) \xrightarrow{\ \d^{q-1}\ } C^{p,q}(X) \longrightarrow \cdots\]
where
\[C^{p,q}(X) := C_{p,q}(X)^\dual \simeq \! \bigoplus_{\delta \in X \\ \dims\delta = q \\ \delta \text{ compact}} \!\! \SF^p(\delta)\]
and the \emph{tropical cohomology} of $X$ is defined by
\[H^{p,q}(X) := H^q(C^{p,\bullet}).\]

We can also define the \emph{compact-dual versions} of tropical homology and cohomology by allowing non-compact faces. These are called Borel-Moore homology and cohomology with compact support, and are defined as follows. We define
\[ C^\BM_{p,q}(X) := \bigoplus_{\delta \in X \\ \dims\delta=q} \SF_p(\delta) \quad \text{and} \quad C_c^{p,q}(X) := \bigoplus_{\delta \in X \\ \dims\delta = q } \SF^p(\comp\delta). \]
We get the corresponding (co)chain complexes $C^\BM_{p,\bul}$ and $C_c^{p,\bul}$, and the \emph{Borel-Moore tropical homology} and the \emph{tropical cohomology with compact support} are respectively
\[ H^\BM_{p,q}(X) := H_q(C^\BM_{p,\bul}(X)) \quad \text{and} \quad H_c^{p,q}(X) := H^q(C_c^{p,\bul}(X)). \]
If $X$ is compact, then both notions of homology and both notions of cohomology coincide.

\smallskip
Here, we defined the cellular versions of tropical homology and cohomology. As in the classical setting, there exist other ways of computing the same groups: for instance using either of singular, cubical, or sheaf cohomologies. We note in particular that the homology and cohomology only depends on the support. Moreover, when it is more convenient, we sometimes use one of these alternative versions.

\medskip

Homology and cohomology in this paper refer to the tropical ones, so we usually omit the mention of the word tropical.

\section{Smoothness in tropical geometry}\label{sec:smoothness}

In this section we study a natural notion of smoothness in tropical geometry. As in differential topology, the smoothness is a local notion, and we consider therefore local charts for tropical varieties which are \emph{supports of tropical fans}.

\subsection{Balancing condition and tropical fans} We start by giving the definition of tropical fans and present some of their basic properties. More results related to algebraic and complex geometry can be found in~\cites{Mik06, Mik07, AR10, GKM09, KM09, Katz12, Bab14, BH17, Gro18} and~\cites{MS15, MR09, BIMS}.

\smallskip
A \emph{tropical fan} $\Sigma$ is a rational fan of pure dimension $d$, for a natural number $d\in \Z_{\geq0}$, which in addition verifies the so-called \emph{balancing condition}: namely, for any cone $\tau$ of dimension $d-1$ in $\Sigma$, we require
\[ \sum_{\sigma \ssupface \tau} \nvect_{\sigma/\tau} \in N_\tau. \]
Equivalently, this means, for any $\tau$ of dimension $d-1$, we have
\begin{align}\label{eq:codimension_one_balancing}
\sum_{\sigma \ssupface \tau} \e^\tau_{\sigma} =0 \quad \textrm{in $N^\tau_\sigma$}.
\end{align}

If the above relation between the vectors $\nvect_{\sigma/\tau}$ is the only one up to multiplication by a scalar, that is, if for any $\tau \in \Sigma_{d-1}$ and real numbers $a_\sigma$ for $\sigma \ssupface \tau$, the relation
\[ \sum_{\sigma \ssupface \tau} a_\sigma\nvect_{\sigma/\tau} \in N_\tau \]
implies that all the coefficients $a_\sigma$ are equal, then we say that $\Sigma$ is \emph{tropically normal}, or simply \emph{normal} if it is understood that $\Sigma$ is a tropical fan, or that $\Sigma$ is \emph{smooth in codimension one}. This is equivalent to requiring the one-dimensional star fan $\Sigma^{\tau}$ to be a tropical line (\ie, its rays be generated by independent unit vectors $e_1, \dots, e_k$ and a last vector $e_0=-e_1-\dots-e_k$), for any codimension one face $\tau$ of $\Sigma$.

\smallskip
Under the \emph{normality condition} for a tropical fan $\Sigma$, we also qualify the canonical compactification $\comp\Sigma$ as being \emph{normal}.

\begin{remark}[Balancing condition and tropicalization] The balancing condition is extended naturally to weighted pure dimensional fans (or polyhedral complexes) where each facet is endowed with an (integral) number called \emph{weight} of the facet. This is the point of view taken in~\cites{GKM09, AR10}, for example, for defining tropical fans and studying their intersection theory. Weighted tropical fans arise naturally in connection with tropicalizations of subvarieties of algebraic tori. In this regard, tropical fans we consider in this paper, with weights all equal to one, are the abstract generalization of those which arise from \emph{faithful tropicalizations}: roughly speaking, if the tropicalization of $X \subseteq T$, for an algebraic tori $T = \spec(\k[M])$, over a trivially valued field $\k$, is a tropical fan $\Sigma \subseteq N_\R$ with all weights equal to one, then $\Sigma$ naturally embeds in the Berkovich analytification $X^{\an}$ of $X$ and forms a skeleton in $X^{\an}$, see~\cites{GRW17, GRW16} for more precise discussion.

Some of the results of this article extend naturally to the more general setting of tropical fans in the presence of a nontrivial weight function. Since our primary purpose in this paper is to discuss the smoothness properties in tropical geometry, in which case all the weights should be equal to one, in order to simplify the presentation, we restrict ourselves to the \emph{reduced} case meaning that we require all the weights to be equal to one on the ambient tropical fan. This being said, we note that even in this case, the weighted balancing condition appears naturally when we will talk later about divisors and Minkowski weights.
\end{remark}

As a follow-up to the above remark we make the following definition.

\begin{defi}[Realizable tropical fan]
A tropical fan $\Sigma$ which arises as the tropicalization of a subvariety $X$ of an algebraic torus $T = \spec(\k[M])$ over a trivially valued field $\k$ is called \emph{realizable over $\k$}. A tropical fan is called \emph{realizable} if it can be realized over some field $\k$.
\end{defi}

\begin{remark}[Balancing condition and smoothness] In continuation of the previous remark, we note that the balancing condition is a crucial ingredient for extracting a tropical notion of smoothness. In fact, it can be seen as a tropical analogue of orientability which allows to define a fundamental class for the fan. In the foundational works~\cites{CLD, GK17, BH17, JSS19}, this is used to establish a tropical analogue of the Stoke's formula, ultimately leading to the Poincaré duality in~\cite{JSS19} under the extra condition for the tropical varieties in question to be \emph{matroidal}. It plays a central role as well in almost all the results we prove in this paper.
\end{remark}

\subsubsection{Localness of the balancing condition}
A property $P$ on rational fans is called \emph{local} or \emph{stellar-stable} if for any fan $\Sigma$ verifying $P$, all the star fans $\Sigma^\sigma$ for $\sigma \in \Sigma$ also verify $P$.

A property $P$ is called a \emph{property of the support}, or we say \emph{$P$ only depends on the support}, if a rational fan $\Sigma$ inside $N_\R$ verifies $P$ if and only if any rational fan $\Sigma'$ considered with a possibly different lattice $N'$ such that $\supp\Sigma = \supp{\Sigma'}$ and $\supp\Sigma \cap N = \supp{\Sigma'}\cap N'$ verifies $P$.

We have the following proposition which summarizes basic properties of tropical fans.

\begin{prop} \label{prop:trop_normal_support_local}
\begin{itemize}
\item \textup{(Support property)} Both the balancing condition and the tropical normality are properties of the support.

\item \textup{(Localness of the balancing condition)} To be tropical, \resp tropically normal, is a local property.
\end{itemize}
\end{prop}
\begin{proof} The first part can be obtained by direct verification. The second part is obtained by identifying star fans of cones in a star fan $\Sigma^\eta$, for $\eta \in\Sigma$, as appropriate star fans in the original fan $\Sigma$.
\end{proof}

\subsubsection{Stability under products} We have the following basic result.

\begin{prop} \label{prop:trop_normal_product} The category of tropical fans and their supports is closed under products. The same holds for the category of normal tropical fans.
\end{prop}
\begin{proof} A codimension one cone in the product $\Sigma \times \Sigma'$ of two fans $\Sigma$ and $\Sigma'$ is of the form $\sigma \times \tau'$ or $\tau \times \sigma'$ for facets $\sigma$ and $\sigma'$ of $\Sigma$ and $\Sigma'$, and codimension one faces $\tau$ and $\tau'$ of $\Sigma$ and $\Sigma'$, respectively. It follows that the star fan of a codimension one face in the product fan can be identified with the star fan of a codimension one face in either $\Sigma$ or $\Sigma'$, from which the result follows.
\end{proof}

\subsubsection{The fundamental class of a tropical fan}

Recall from Section~\ref{sec:prel} that to a rational fan $\Sigma$ of pure dimension $d$, we associate its canonical element $\nu_{\Sigma}$
\[ \nu_\Sigma = (\nu_\sigma)_{\sigma \in \Sigma_d} \in \bigoplus_{\sigma \in \Sigma_d} \bigwedge^d N_\sigma. \]
Here, for each facet $\sigma$ of $\Sigma$, $\nu_\sigma$ is the canonical multivector of $\sigma$ which is a generator of $\bigwedge^{\dims{\sigma}}N_\sigma = \SF_{\dims\sigma}(\sigma)$. We have the following well-known property, see for example~\cites{BH17, CLD}.

\begin{prop}\label{prop:trop_canonical_cycle} The fan $\Sigma$ is tropical if and only if $\nu_\Sigma$ is a cycle of $C_{d,d}^\BM(\Sigma)$.
\end{prop}
In this case, by an abuse of the notation, we denote by $\nu_\Sigma$ the corresponding element in $H^\BM_{d,d}(\Sigma)$ and call it the \emph{fundamental class of\/ $\Sigma$}.

\subsubsection{Tropical irreducibility} We now present a tropical notion of \emph{irreducibility} which as we show below will be directly related to (but slightly stronger than) the normality of tropical fans introduced in the previous section. This is a property which we will need to assume at some occasions in the paper.

\smallskip
\begin{defi}[Tropical irreducibility] A tropical fan $\Sigma$ is called \emph{irreducible at $\eta$}, for a cone $\eta\in\Sigma$, if the homology group $H^\BM_{d-\dims\eta, d-\dims\eta}(\Sigma^\eta)$ is generated by the fundamental cycle $\nu_{\Sigma^\eta}$. A tropical fan irreducible at $\conezero$ is simply called \emph{irreducible}.

The fan $\Sigma$ is called \emph{locally irreducible} if it is irreducible at any cone $\eta\in\Sigma$.
\end{defi}

\begin{prop}[Localness of local irreducibility]
Being locally irreducible is a local property of tropical fans.
\end{prop}
\begin{proof}
This is tautological from the definition.
\end{proof}

The following theorem gives a link between normality and local irreducibility.

\begin{thm}[Characterization of locally irreducible tropical fans] \label{thm:characterization_irreducible} A tropical fan $\Sigma$ is locally irreducible if and only if it is normal and each star fan $\Sigma^\eta$, $\eta \in \Sigma$, is connected through codimension one.
\end{thm}

\begin{proof} Suppose $\Sigma$ is locally irreducible. Consider a cone $\tau$ of codimension one in $\Sigma$. The star fan $\Sigma^\tau$ is an irreducible tropical fan of dimension one. Let $\rho_1, \dots, \rho_k$ be the rays of $\Sigma^\tau$. Then the balancing condition for $\Sigma^\tau$ implies that the vectors $\e_{\rho_1}, \dots, \e_{\rho_k}$ sum up to zero, and the local irreducibility implies that the canonical cycle $\nu_{\Sigma^\tau}$ is a generator of $H_{1,1}^\BM(\Sigma^\tau)$. For any relation in $N^\tau$ of the form
\[\sum_{j=1}^k a_j \e_{\rho_j} =0\]
for scalars $a_j$, we get an element $\alpha \in H_{1,1}^\BM(\Sigma^\tau)$. Since $\Sigma^\tau$ is irreducible, $\alpha$ is a multiple of $\nu_{\Sigma^\tau}$ which implies that the scalars $a_j$ are all equal. This proves that $\Sigma$ is normal.

Consider now a cone $\eta\in \Sigma$. It remains to prove that $\Sigma^\eta$ is connected through codimension one. Suppose this is not the case. Then we can find a partition of the facets of $\Sigma^\eta_d$ into a disjoint union $S_1 \sqcup \dots \sqcup S_l$, for some $l\geq 2$, so that $S_j$ form the connected components of the connectivity-through-codimension-one graph of the facets of $\Sigma^\eta$. It follows that each element $\alpha_j:=(\nu_{\sigma})_{\sigma \in S_j}$ form a cycle in $H_{d,d}^\BM(\Sigma^\eta)$ which contradicts the irreducibility of $\Sigma^\eta$.

\smallskip
We now prove the reverse implication. Suppose $\Sigma$ is normal and that for each $\eta\in \Sigma$, $\Sigma^\eta$ is connected through codimension one. Let $\eta$ be a face of $\Sigma$ and let $\alpha = \sum_{\sigma\in \Sigma^\eta_{d-\dims\eta}} a_\sigma\nu_{\sigma}$ be an element in $H_{d-\dims \eta,d-\dims \eta}^\BM(\Sigma^\eta)$. Applying the normality condition to a codimension one face $\tau$ of $\Sigma^\eta$, we infer that the scalars $\alpha_\sigma$ are all equal for $\sigma \ssupface \tau$. Using now the connectivity of $\Sigma^\eta$ through codimension one, we conclude that the scalars $a_\sigma$ are all equal. This shows $\alpha$ is a multiple of $\nu_{\Sigma^\eta}$ and the irreducibility at $\eta$ follows. Since this holds for any $\eta$, we conclude the local irreducibility of $\Sigma$.
\end{proof}

Theorem \ref{thm:characterization_irreducible} justifies the following definition via Corollary \ref{cor:irreducible_components_irreducible}.

\begin{defi}[Irreducible components of a normal tropical fan] Let $\Sigma$ be a normal tropical fan of dimension $d$. Consider the connectivity-through-codimension-one graph $G = (\Sigma_{d}, E)$ of $\Sigma$. Any connected component of $G$ with vertex set $V \subseteq \Sigma_{d}$ defines a subfan of $\Sigma$ which is tropical, normal and connected through codimension one. We refer to these as the \emph{irreducible components} of $\Sigma$.
\end{defi}

\begin{cor} \label{cor:irreducible_components_irreducible}
Let $\Sigma$ be a normal tropical fan of dimension $d$. The irreducible components of\/ $\Sigma$ are irreducible. Moreover, they induces a partition of the facets $\Sigma_d$.
\end{cor}

\begin{proof}
This follows from the proof of Theorem \ref{thm:characterization_irreducible} and from the definition.
\end{proof}

Note that the statement in the corollary claims the irreducibility only at $\conezero$ (this is weaker than local irreducibility by Example \ref{ex:PD_chow_not_smooth}).

\begin{remark}
Another natural definition for irreducible components of a tropical fan $\Sigma$ is the following. A subfan $\Delta$ of $\Sigma$ is an irreducible component of $\Sigma$ if it is the support of a nonzero element in $H^\BM_{d,d}(\Sigma)$, and if it is minimal among the subfans with this property. In the case $\Sigma$ is normal, this definition coincides with the one given above. However, for general tropical fans, irreducible components might not induce a partition of $\Sigma_d$ (\cf Example \ref{ex:irreducible_components_not_well_defined}).
\end{remark}

\begin{prop} \label{prop:irreducibility_product}
The category of locally irreducible tropical fans is closed under products.
\end{prop}
\begin{proof}
This follows from Theorem~\ref{thm:characterization_irreducible} and the same property for the category of tropical normal fans and fans which are connected through codimension one.
\end{proof}

\subsection{Poincaré duality and cap product} \label{subsec:poincare_duality} Recall that the coefficient ring of our homology and cohomology groups is $\Z$ unless otherwise stated.

Let $\Sigma$ be a simplicial tropical fan. Then $H_{d,d}^\BM(\Sigma)$ contains a canonical element $\nu_\Sigma$. Using the cap product $\frown$
\[ \frown\colon H^{p_1,q_1}(\Sigma) \times H^\BM_{p_2,q_2}(\Sigma) \longrightarrow H^\BM_{p_2-p_1, q_2-q_1}(\Sigma), \]
that we describe below, we get a natural map
\[ \cdot\frown\nu_\Sigma\colon H^{p,q}(\Sigma) \longrightarrow H^\BM_{d-p,d-q}(\Sigma). \]

\begin{defi}[Integral Poincaré duality] We say that a simplicial tropical fan $\Sigma$ verifies the \emph{integral Poincaré duality} if the above map is an isomorphism for any bidegree $(p,q)$.
\end{defi}

Since the homology and cohomology groups of $\Sigma$ only depend on the support $\supp\Sigma$, we get the following.

\begin{prop}
The integral Poincaré duality is a property of the support.
\end{prop}

The above notion of duality can be defined in the same way for the compactification of $\Sigma$. In this case, $H^{d,d}(\comp\Sigma)$ still contains a canonical element $\nu_{\comp\Sigma}$. Since $\Sigma$ is simplicial, $\comp\Sigma$ is a cubical complex, and so the cap product can be defined. For unimodular tropical fans, we get the following theorem. (See Theorem~\ref{thm:smooth_PD} for a more general statement.)
\begin{thm}\label{thm:PD_canonical_compactification}
Let $\Sigma$ be a unimodular tropical fan. Assume that for any face $\sigma$ of\/ $\Sigma$ (including~$\conezero$), $\Sigma^\sigma$ verifies the integral Poincaré duality. Then $\comp\Sigma$ verifies the integral Poincaré duality.
\end{thm}

\begin{proof}
There is a natural cover of $\Sigma$ by open sets $\U_\sigma$ where $\U_\sigma$ is the union of all relative interior of faces containing $\infty_\sigma$. Then $\U_\sigma$ is isomorphic to a nice open set of $\Sigma^\sigma \times \eR^{\dims{\sigma}}$ and the arguments of \cite{JRS18}*{Section 5} can be applied to prove the Poincaré duality for the $\U_\sigma$ and to glue them together. An alternate proof can be obtained from the results we prove later in Section~\ref{sec:hodge_isomorphism}, see Theorem~\ref{thm:smoothness_alternate} and below.
\end{proof}

\subsubsection{Description of the cohomology of fans and of the cap product} Let $\Sigma$ be a tropical fan. The cohomology of $\Sigma$ is very simple since the only compact face of $\Sigma$ is $\conezero$. Thus,
\[ H^{p,q}(\Sigma) \simeq C^{p,q}(\Sigma) \simeq \begin{cases}
  \SF^p(\conezero) & \quad \text{if $q = 0$,} \\
  0              & \quad \text{otherwise.}
\end{cases} \]

Keeping this in mind, we now describe the cap product on $\Sigma$. Recall from Section~\ref{sec:intro-basic-notations} that for a $k$-form $\alpha \in \bigwedge^k M$, we denote by $\iota_\alpha$ the contraction of multivectors by $\alpha$. Contraction naturally extends to chains by linearity and gives a map
\[\iota_\alpha \colon C_{p,q}^\BM(\Sigma) \to C_{p-k,q}^\BM(\Sigma).\]
Identifying now $C^{k,0}(\Sigma)$ with $\SF^k(\conezero) \subseteq \bigwedge^k M$, we can view this as the cap product.

\subsubsection{Reformulation of the integral Poincaré duality} \label{subsubsec:reformulation_PD} The discussion from the previous section leads to a map
\[ \SF^p(\conezero) \to H^\BM_{d-p,d}(\Sigma), \qquad \alpha \mapsto \iota_\alpha(\nu_\Sigma). \]
For any $p$, this is an injective map. The surjectivity of this map in any degree as well as the vanishing of $H_{p,q}^\BM(\Sigma)$ for $q\neq d$ is the essence of the integral Poincaré duality for $\Sigma$.

\subsubsection{Stable invariance of the Poincaré duality} A property $P$ of fans which only depends on the support is called \emph{stably invariant} if for a fan $\Sigma$ and a positive integer $k$, fans supported on $\supp\Sigma$ verify $P$ if and only if those supported on $\supp{\R^k\times \Sigma}$ verify $P$.

\begin{prop}\label{prop:stable_invariance_pd} Poincaré duality is stably invariant.
\end{prop}

We actually prove the following stronger property.
\begin{prop}\label{prop:pd_product} Let $\Sigma$ and $\Sigma'$ be two tropical fans and consider the product $\Sigma\times \Sigma'$ which is a tropical fan. Then $\Sigma \times \Sigma'$ verifies the Poincaré duality if and only if\/ $\Sigma$ and $\Sigma'$ verify the Poincaré duality.
\end{prop}

\begin{proof}[Proof of Proposition~\ref{prop:stable_invariance_pd}] Any fan $\Sigma'$ supported in $\R^k$ verifies the Poincaré duality. Applying the above proposition to $\Sigma$, $\Sigma'$ and the product $\Sigma\times \Sigma'$ leads to the result.
\end{proof}

\begin{proof}[Proof of Proposition~\ref{prop:pd_product}] Let $d$ and $d'$ be the dimensions of $\Sigma$ and $\Sigma'$, respectively.

We first prove that if both $\Sigma$ and $\Sigma'$ verify the Poincaré duality, then so does $\Sigma \times \Sigma'$.

By Künneth decomposition~\cite{GS-sheaf}, since the cohomology groups $H^{\bul,\bul}(\Sigma)$ and $H^{\bul,\bul}(\Sigma')$ are torsion-free, and since by Poincaré duality the same holds for the Borel-Moore homology groups, we get decompositions
\[H^{p,q}(\Sigma \times \Sigma') \simeq \bigoplus_{0\leq r\leq p\\ 0\leq s \leq q} H^{r,s}(\Sigma) \otimes H^{p-r, q-s}(\Sigma'), \qquad \textrm{and} \]
\[H^\BM_{p,q}(\Sigma \times \Sigma') \simeq \bigoplus_{0\leq r\leq p\\ 0\leq s \leq q} H^\BM_{r,s}(\Sigma) \otimes H^\BM_{p-r, q-s}(\Sigma').\]
The Poincaré duality pairing (the cap product) is compatible with the Künneth decomposition. This leads to the assertion that $\Sigma \times \Sigma'$ verifies the Poincaré duality.

\smallskip
For the other direction, assume $\Sigma \times \Sigma'$ verifies the Poincaré duality. We show this holds as well for $\Sigma$, and by symmetry, we get the result for $\Sigma'$. The cohomology $H^{\bul,\bul}(\Sigma\times\Sigma')$ is torsion-free, and the Poincaré duality implies that $H^\BM_{\bul,\bul}(\Sigma\times\Sigma')$ is torsion-free. Applying Künneth formula, we get that the cohomology groups and Borel-Moore homology groups of $\Sigma$ and of $\Sigma'$ are torsion-free and that the previous decompositions hold.

Once again, the cap product respects these decompositions. Hence, for any integers $p, q$, the isomorphism
\[ H^{p,q}(\Sigma\times\Sigma') \simto H^\BM_{d+d'-p,d+d'-q}(\Sigma\times\Sigma') \]
induces an isomorphism on each part: for each pair $r,s$ with $0\leq r \leq p$ and $0 \leq s \leq q$,
\[ H^{r,s}(\Sigma) \otimes H^{p-r,q-s}(\Sigma') \simto H^\BM_{d-r,d-s}(\Sigma) \otimes H^\BM_{d'-p+r,d'-q+s}(\Sigma'). \]
In particular, for $r = p$ and $s = q$, we get the Poincaré duality for $\Sigma$ in bidegree $(p,q)$.
\end{proof}

\subsection{Tropical smoothness}

\begin{defi}[Tropical smoothness]
A tropical fan $\Sigma$ is called \emph{smooth} if for any cone $\sigma \in \Sigma$, the star fan $\Sigma^\sigma$ verifies the Poincaré duality.
\end{defi}

\begin{thm}\label{thm:smoothness-support} The property of being smooth depends only on the support of the tropical fan.
\end{thm}
\begin{proof}
Let $\Sigma$ and $\Sigma'$ be two tropical fans with the same support $X$. Let $x$ be a point of $X$ and denote by $\sigma\in\Sigma$ and $\tau\in \Sigma'$ the cones of $\Sigma$ and $\Sigma'$ which contain $x$ in their relative interiors. It follows that we have the isomorphism $\R^{\dims\sigma}\times \Sigma^\sigma \simeq \R^{\dims{\tau}} \times \Sigma'^{\tau}$. Applying Proposition~\ref{prop:stable_invariance_pd}, we infer that $\Sigma^\sigma$ verifies the Poincaré duality if and only if $\Sigma'^{\tau}$ does. Since this happens for all points of $X$, the result follows.
\end{proof}

\begin{prop} \label{prop:smoothness_product} Let $\Sigma$ and $\Sigma'$ be two tropical fans. Then $\Sigma\times\Sigma'$ is smooth if and only if $\Sigma$ and $\Sigma'$ are smooth. In particular, the category of smooth tropical fans is stable under products.
\end{prop}
\begin{proof}
These follow from Proposition \ref{prop:pd_product}, and the fact that for any face $\sigma\times\sigma'$ of $\Sigma\times\Sigma'$, $\bigl(\Sigma\times\Sigma'\bigr)^{\sigma\times\sigma'} \simeq \Sigma^\sigma \times \Sigma'^{\sigma'}$.
\end{proof}

\begin{defi} \begin{itemize}
\item Let $X$ be the support of a tropical fan. We say $X$ is smooth if one, and so any, tropical fan $\Sigma$ with support $X$ is smooth.

\item A tropical variety $X$ is called smooth if each point of $X$ has an open neighborhood which is isomorphic to an open neighborhood of a point in the canonical compactification $\comp\Sigma$ of a smooth simplicial tropical fan $\Sigma$. \qedhere
\end{itemize}
\end{defi}

A tropical variety in the above definition means a connected topological space with a compatible atlas of charts each identified with an open subset of a canonically compactified unimodular tropical fan such that the maps which give transition between the charts are affine integral linear. We omit the formal definition here and refer to~\cite{AP-tht} and \cite{JSS19} for more details.

\subsection{Integral Poincaré duality for smooth tropical varieties} As in the local setting, a tropical variety $X$ of dimension $d$ comes with its canonical class $\nu_X \in H_{d,d}^{^{\sing,}\BM}(X)$. Here, $H_{\bul,\bul}^{^{\sing,}\BM}$ refers to singular Borel-Moore tropical homology as defined in \cite{JRS18}. (As previously mentioned, in the case $X = \Sigma$ or $\comp \Sigma$,  this coincides with the previously introduced Borel-Moore tropical homology.)

We say $X$ verifies the Poincaré duality with integral coefficients if the cap product
\[\cdot \frown \nu_X\colon H_{_\sing}^{p,q}(X) \to H_{d-p, d-q}^{^{\sing,}\BM}(X)\]
is an isomorphism.

The following is the general form of Theorem~\ref{thm:PD_canonical_compactification}.
\begin{thm}[Poincaré duality for smooth tropical varieties] \label{thm:smooth_PD} A smooth tropical variety verifies the Poincaré duality.
\end{thm}
\begin{proof}
Combined with the results which preceded, concerning the calculation of the cohomology and the integral Poincaré duality for tropical fans, the proof is similar to the one given in~\cite{JRS18} for matroidal tropical varieties. We refer as well to~\cite{GS-sheaf} for a sheaf-theoretic approach and to~\cite{JSS19} for a proof of the duality for cohomology with rational coefficients. Both these proofs are presented for tropical varieties which are locally matroidal, but combined with the results we proved in this section, they can be applied to the general setting considered in this paper.
\end{proof}

\begin{remark} We will later give Example~\ref{subsec:cube} which is a normal unimodular tropical fan $\Sigma'$ with $\rk(H^{2,1}(\comp\Sigma')) = 2$ but $\rk(H^{1,2}(\comp\Sigma')) = 0$. This shows the theorem does not hold in general only assuming unimodularity, and the smoothness assumption is needed in general to ensure the Poincaré duality for the canonical compactification $\comp\Sigma$ of a tropical fan $\Sigma$.
\end{remark}

\section{Tropical divisors}\label{sec:divisors}

In this section, we consider a tropical fan $\Sigma$ of pure dimension $d$ in $N_\R \simeq \R^n$ for some $n$ and study divisors associated to conewise linear functions.

\subsection{Minkowski weights on rational fans} Let $p, d$ be two non-negative integers with $p\leq d$. Let $\Sigma$ be a rational fan of dimension $d$ in $N_\R$.

Assume for each cone $\sigma$ of $\Sigma$ of dimension $d-p$ we are given a weight which is an integer denoted by $w(\sigma)$. Let $C := (\Sigma_{(d-p)}, w)$ be the corresponding weighted fan with the weight function $w$ on the facets of $\Sigma_{(d-p)}$. The weight function $w$ is called a \emph{Minkowski weight} of dimension $d-p$ on $\Sigma$ if the following balancing condition is verified:
\[\forall\:\tau\in \Sigma_{d-p-1},\qquad \sum_{\sigma\ssupface\tau} w(\sigma) \nvect_{\sigma/\tau} = 0 \in N^{\tau}.\]
We denote by $\MW_{d-p}(\Sigma)$ the set of all Minkowski weights of dimension $d-p$ on $\Sigma$. Addition of weights cell by cell turns $\MW_{d-p}(\Sigma)$ into a group.

Note that in the case where $\Sigma$ is a tropical fan, we get a \emph{canonical element} denoted by $[\Sigma]$ in $\MW_d(\Sigma)$. This is given by the weight function 1 on facets of $\Sigma$.

\subsection{Divisors on tropical fans} Let now $\Sigma$ be a tropical fan. A \emph{divisor} of $\Sigma$ is the data of a pair $(\Delta, w)$ consisting of a (possibly empty) subfan $\Delta \subset \Sigma$ of dimension $d-1$ and a weight function $w\colon \Delta_{d-1} \to \Z\setminus \{0\}$ such that the weighted fan $\Delta$ is balanced. This means for any cone $\eta \in \Delta$ of dimension $d-2$, the balancing condition
\[\sum_{\tau \in \Delta_{d-1} \\ \tau \ssupface \eta} w(\tau) \nvect_{\tau/\eta} \in N_\eta\]
holds. Equivalently, a divisor in $\Sigma$ is any element of $\MW_{d-1}(\Sigma)$. In this case, the fan $\Delta$ will be the support of $w$, \ie, the subfan of $\Sigma$ defined by all $\tau \in \Sigma_{d-1}$ with $w(\tau)\neq 0$.

\smallskip
In what follows, we denote by $\Div(\Sigma)$ the group of divisors on $\Sigma$ and note that we have $\Div(\Sigma) = \MW_{d-1}(\Sigma)$.

\subsection{Principal divisor associated to a conewise integral linear function}

Let $f\colon \Sigma \to \R$ be a conewise integral linear function on $\Sigma$ as defined in Section \ref{subsec:fans}, \ie, $f$ is continuous and it is integral linear on each cone $\sigma$ of $\Sigma$. Recall that for each face $\eta$ of $\Sigma$, we denote by $f_\eta$ the linear form induced by $f$ on $N_{\eta,\R}$.

Let $\tau$ be a face of codimension one in $\Sigma$. The \emph{order of vanishing of $f$ along $\tau$} denoted by $\ord_\tau(f)$ is defined as
\[ \ord_\tau(f) := -\sum_{\sigma \ssupface \tau} f_\sigma(\nvect_{\sigma/\tau}) + f_\tau\Bigl(\sum_{\sigma \ssupface \tau} \nvect_{\sigma/\tau}\Bigr)\]
with the sums running over over all cones $\sigma \in \Sigma_d$ containing $\tau$.

\begin{prop} Notations as above, the order of vanishing $\ord_\tau(f)$ is well-defined, that is, it is independent of the choice of normal vectors $\nvect_{\sigma/\tau} \in N_\sigma$.
\end{prop}
\begin{proof} For each pair $\sigma \ssupface \tau$, two different choices $\nvect_{\sigma/\tau}$ and $\nvect'_{\sigma/\tau}$ of normal vectors differ by a vector in $N_\tau$. It follows that
\begin{align*}
-\sum_{\sigma \ssupface \tau} f_\sigma(\nvect_{\sigma/\tau}) &+ f_\tau\Bigl(\sum_{\sigma \ssupface \tau} \nvect_{\sigma/\tau}\Bigr) + \sum_{\sigma \ssupface \tau} f_\sigma(\nvect'_{\sigma/\tau}) - f_\tau\Bigl(\sum_{\sigma \ssupface \tau} \nvect'_{\sigma/\tau}\Bigr)\\
&= - \sum_{\sigma \ssupface \tau} f_\sigma(\nvect_{\sigma/\tau} - \nvect'_{\sigma/\tau}) + f_\tau\Bigl(\sum_{\sigma \ssupface \tau} \nvect_{\sigma/\tau} - \sum_{\sigma \ssupface \tau} \nvect'_{\sigma/\tau}\Bigr)\\
&= - \sum_{\sigma \ssupface \tau} f_\tau(\nvect_{\sigma/\tau} - \nvect'_{\sigma/\tau}) + f_\tau\Bigl(\sum_{\sigma \ssupface \tau} \nvect_{\sigma/\tau} - \sum_{\sigma \ssupface \tau} \nvect'_{\sigma/\tau}\Bigr) =0. \qedhere
\end{align*}
\end{proof}

The order of vanishing function gives a weight function $\ord(f)\colon \Sigma_{d-1} \to \Z$. We associate to $f$ the data of the pair $(\Delta, w)$ where $\Delta$ is the fan defined by the support of $\ord(f)$, \ie, by those cones $\tau$ of dimension $d-1$ in $\Sigma$ for which we have $\ord_\tau(f) \neq 0$, and the weight function $w \colon \Delta_{d-1} \to \Z\setminus\{0\}$ is given by $w(\tau) = \ord_\tau(f)$. We have the following well-known result, see for example \cite{AR10}.
\begin{prop}\label{prop:balancing_divisor} Notations as above, the pair $(\Delta, w)$ is a divisor.
\end{prop}

\begin{proof} For the sake of completeness, we provide a proof. We need to prove the balancing condition around each cone $\eta$ in $\Delta$ of dimension $d-2$. It will be enough to work in $\Sigma$ and extend the weight function $w$ by $0$ on all cones which are not included in $\Delta$.

Replacing $f$ with $f-f_\eta$, and passing to the quotient by $N_{\eta,\R}$, we can assume that $\Sigma$ is a tropical fan of dimension two and $\eta = \conezero$. In this case, we have to show that
\[\sum_{\tau \in \Sigma_1} \ord_{\tau}(f) \e_\tau =0.\]

We choose primitive vectors $\nvect_{\sigma/\tau}$ in $N$ for any pair of non-zero cones $\sigma\ssupface \tau$ in $\Sigma$ such that we have the balancing condition
\begin{equation}\label{eq:bal}
\forall \tau\in \Sigma_1, \qquad \sum_{\sigma\ssupface \tau} \nvect_{\sigma / \tau} =0.
\end{equation}
Consider a pair $\sigma \ssupface \tau$ and denote by $\tau'$ the other ray of $\sigma$. We can write
\[ b_{\tau, \tau'}\, \nvect_{\sigma/\tau} = \e_{\tau'} + a_{\tau, \tau'}\,\e_\tau\]
for a unique pair of integers $a_{\tau, \tau'}$ and $b_{\tau, \tau'}$. Note that $b_{\tau, \tau'} \neq 0$. Equation~\eqref{eq:bal} now reads
\[\sum_{\tau' \sim \tau} \frac 1{b_{\tau, \tau'}}\e_{\tau'} = -\Bigl(\sum_{\tau'\sim \tau} \frac{a_{\tau,\tau'}}{b_{\tau, \tau'}}\Bigr) \e_\tau, \]
where, we recall, the notation $\tau'\sim \tau$ means $\tau'$ and $\tau$ are the two rays of a two dimensional cone in $\Sigma$. To conclude note that
\begin{align*}
\sum_{\tau \in \Sigma_1} \ord_\tau(f) \e_\tau =\hspace{-1cm}&\hspace{1cm} -\sum_{\tau\in \Sigma_1}\,\sum_{\sigma \ssupface \tau} f_\sigma(\nvect_{\sigma/\tau}) \e_\tau + \sum_{\tau \in \Sigma_1} f_\tau\Bigl(\sum_{\sigma\ssupface \tau}\nvect_{\sigma/\tau}\Bigr) \e_\tau \\
  &= -\sum_{\tau\in \Sigma_1}\,\sum_{\sigma \ssupface \tau} f_\sigma(\nvect_{\sigma/\tau}) \e_\tau \\
  & = - \sum_{\tau\in \Sigma_1}\,\sum_{\tau' \sim \tau} f_\sigma\Bigl( \frac{\e_{\tau'}}{b_{\tau, \tau'}} + \frac{a_{\tau,\tau'} \e_\tau}{b_{\tau, \tau'}}\Bigr) \e_\tau & \hspace{-.5cm} \textrm{(with $\sigma = \tau \vee \tau'$)}\\
  &= - \sum_{\tau\in \Sigma_1}\,\sum_{\tau' \sim \tau} \frac {f_{\tau'}(\e_{\tau'})}{b_{\tau, \tau'}} \e_\tau - \sum_{\tau\in \Sigma_1}\,\sum_{\tau' \sim \tau} \frac{f_{\tau}(\e_\tau)a_{\tau,\tau'}}{b_{\tau, \tau'}} \e_\tau & \hspace{-.5cm} \textrm{(since $f_\sigma\rest{\tau} = f_\tau$ for $\sigma\ssupface \tau$)}\\
  &= - \sum_{\tau\in \Sigma_1}\,\sum_{\tau' \sim \tau} \frac{f_{\tau'}(\e_{\tau'})}{b_{\tau, \tau'}} \e_\tau + \sum_{\tau\in \Sigma_1}\,f_{\tau}(\e_\tau) \Bigl(-\sum_{\tau' \sim \tau} \frac{a_{\tau,\tau'}}{b_{\tau, \tau'}}\Bigr) \e_\tau\\
  &= - \sum_{\tau\in \Sigma_1}\,\sum_{\tau' \sim \tau} \frac{f_{\tau'}(\e_{\tau'})}{b_{\tau, \tau'}} \e_\tau + \sum_{\tau\in \Sigma_1}\,f_{\tau}(\e_\tau) \sum_{\tau'\sim \tau} \frac 1{b_{\tau, \tau'}}\e_{\tau'} \\
  &= - \sum_{\tau\in \Sigma_1}\,\sum_{\tau' \sim \tau} \frac{f_{\tau'}(\e_{\tau'})}{b_{\tau, \tau'}} \e_\tau + \sum_{\tau\in \Sigma_1}\sum_{\tau'\sim \tau}\,\frac{f_{\tau}(\e_\tau)}{b_{\tau, \tau'}} \e_{\tau'}=0.
\end{align*}
To see the last equality, note that $b_{\tau, \tau'} = b_{\tau', \tau}$, both terms being equal to the covolume of the sublattice of $N$ generated by $\e_\tau$ and $\e_{\tau'}$. Using this and the symmetry, we infer that the difference of the last sums is zero.
\end{proof}

\begin{defi}[Principal divisors] Let $\Sigma$ be a tropical fan. For any conewise integral linear function on $\Sigma$, we denote by $\div(f)$ the divisor associated to $f$. Such divisors are called \emph{principal}. Principal divisors form a subgroup of $\Div(\Sigma)$ that we denote by $\Prin(\Sigma)$.
\end{defi}

\begin{defi}[Effective and reduced divisors]
A divisor $(\Delta, w)$ on $\Sigma$ is called \emph{effective} if all the coefficients $w(\tau)$ for $\tau$ a $(d-1)$-dimensional cone in $\Delta$ are positive. An effective divisor is called \emph{reduced} if all the weights are equal to one, \ie, $w(\tau)=1$ for any facet $\tau$ of $\Delta$. In such a case, we simply omit the mention of the weight function.
\end{defi}

\subsection{Principal and div-faithful tropical fans}
We now define a class of tropical fans on which divisors behave nicely.

\begin{defi} Let $\Sigma$ be a tropical fan and let $\eta$ be a cone in $\Sigma$.
\begin{itemize}
\item We say that $\Sigma$ is \emph{principal at $\eta$} if any divisor on $\Sigma^\eta$ is principal. We call the tropical fan $\Sigma$ \emph{principal} if $\Sigma$ is principal at any cone $\eta \in \Sigma$.
\item We say that $\Sigma$ is \emph{divisorially faithful at $\eta$} or simply \emph{div-faithful at $\eta$} if for any conewise integral linear function $f$ on $\Sigma^\eta$, if $\div(f)$ is trivial in $\Div(\Sigma^\eta)$, then $f$ is a linear function on $\Sigma^\eta$. We call the tropical fan $\Sigma$ \emph{div-faithful} if $\Sigma$ is div-faithful at any cone $\eta\in \Sigma$. \qedhere
\end{itemize}
\end{defi}

The importance of div-faithfulness in our work relies on the fact that tropical modifications behave very nicely on div-faithful tropical fans, \cf Section \ref{sec:star_fans_tropical_modifications}.

\begin{prop}\label{prop:div-faithful_local} The properties of being principal and div-faithful are both local.
\end{prop}
\begin{proof} The statement is tautological.
\end{proof}

\begin{remark} Examples \ref{ex:non-principal_fan} and \ref{ex:PD_chow_not_smooth} show that being principal, \resp div-faithful, at $\conezero$ does not imply the fan is principal, \resp div-faithful.
\end{remark}

\subsection{Characterization of principal and div-faithful saturated unimodular tropical fans} \label{subsec:characterization_principal_div-faithful}

We now provide a characterization of principality and divisorial faithfulness in the case the tropical fan is \emph{saturated and unimodular}. This will be given via the map \[\cycl\colon A^1(\Sigma) \to \Div(\Sigma) \simeq \MW_{d-1}(\Sigma)\]
that we describe now. Here, as in the introduction, $A^\bul(\Sigma)$ refers to the Chow ring of $\Sigma$.

We start by recalling the definition of the above map. First, for any $k$, we get a natural isomorphism
\[ \MW_{k}(\Sigma) \simeq A^k(\Sigma)^\dual.\]
This is~\cite{AHK}*{Proposition 5.6}, and is a consequence of Localization Lemma~\ref{lem:kernel_Z^k}, see Remark~\ref{rem:chow_to_mw}.

For $k=d$, this gives the \emph{degree map}
\[\deg\colon A^d(\Sigma) \to \Z, \qquad \alpha \to \int_{\Sigma} \alpha\]
where the notation $\int_\Sigma$ means evaluation at the canonical element $[\Sigma]$ of $\MW_d(\Sigma)$ given by the weight function $1$ on $\Sigma_d$, see Remark~\ref{rem:chow_to_mw}.

The composition of the product map in the Chow ring with the degree
\[ \begin{tikzcd}[column sep=small, row sep=0pt]
A^k(\Sigma) \times A^{d-k}(\Sigma) \rar& A^d(\Sigma) \rar& \Z,\qquad\\
(\alpha, \beta) \ar[rr, mapsto]&& \deg(\alpha \cdot \beta),
\end{tikzcd} \]
gives the map
\[\cycl^k\colon A^k(\Sigma) \to A^{d-k}(\Sigma)^\dual \simeq \MW_{d-k}(\Sigma).\]
For $k=1$, this gives the map
\[\cycl\colon A^1(\Sigma) \to \Div(\Sigma).\]

\smallskip
Notice also that, since $\Sigma$ is unimodular, $Z^1(\Sigma)$ coincides with the set of conewise integral linear forms on $\Sigma$. Moreover, since $\Sigma$ is saturated, $A^1(\Sigma)$ is isomorphic to conewise integral linear forms modulo globally linear ones.

\begin{prop} \label{prop:map_cl} The map $\cycl$ is described as follows. Consider an element $\alpha=\sum_{\rho\in \Sigma_1} a_\rho \x_\rho$ in $A^1(\Sigma)$ with coefficients $a_\rho \in \Z$. Let $f$ be the conewise linear function on $\Sigma$ which takes value $a_\rho$ at $\e_\rho$, for any ray $\rho$. Then we have $\cycl(\alpha) = -\div(f)$.
\end{prop}

\begin{proof} Notations as above, we need to show that for each $\tau \in \Sigma_{d-1}$, we have the equality
\[\deg(\alpha \cdot \x_\tau) = -\ord_\tau(f). \]
Let $\ell$ be an integral linear function on $N$ which is equal to $f$ on $\tau$ and let $\x_\ell = \sum_{\rho\in \Sigma_1} \ell(\e_\rho)\x_\rho$. Since $\x_\ell$ is zero in the Chow ring, replacing $\alpha$ by $\alpha -\x_\ell$, we can assume that $a_\rho=0$ for any ray $\rho$ in $\tau$. The proposition now follows by observing that
\[\deg(\alpha \cdot \x_\tau) = \sum_{\rho\in \Sigma_1 \\ \rho\sim \tau} a_\rho \deg(\x_\rho \x_\tau) = \sum_{\rho\in \Sigma_1 \\ \rho\sim \tau} a_\rho = -\ord_\tau(f-\ell) = -\ord_\tau(f). \qedhere \]
\end{proof}

\begin{thm}[Characterization of principal and div-faithful fans]\label{thm:char_div-faithful} Let $\Sigma$ be a unimodular tropical fan of dimension $d$ and let $\eta$ be a cone in $\Sigma$. Assume $\Sigma$ is saturated at $\eta$. Consider the map $\cycl_\eta\colon A^1(\Sigma^\eta) \to \Div(\Sigma^\eta)$. We have
\begin{itemize}
\item $\Sigma$ is principal at $\eta$ if and only if $\cycl_\eta$ is surjective.
\item $\Sigma$ is div-faithful at $\eta$ if and only if $\cycl_\eta$ is injective.
\end{itemize}
\end{thm}
\begin{proof} Both parts follow directly from Proposition~\ref{prop:map_cl} applied to $\Sigma^\eta$.
\end{proof}

\begin{remark} \label{rem:independent_exterior_lattice}
Notice that principality and div-faithfulness are independent of the lattice outside the support of the fan. By this, we mean that a tropical fan $\Sigma$ in $N_\R$ is principal, \resp div-faithful, if and only if for any other lattice $N'$ which verifies $\supp\Sigma \cap N = \supp\Sigma \cap N'$, $\Sigma$ is principal, \resp div-faithful, when considered equipped with this new lattice.

For such properties, replacing the lattice $N$ by $\SF_1(\Sigma)$ if necessary, there is no harm in assuming that the fan is saturated at $\conezero$ (though we cannot assume global saturation by Example \ref{ex:non_saturated_smooth_fan}).
\end{remark}

\begin{thm} \label{thm:smooth_principal_div-faithful}
A smooth tropical fan is both principal and div-faithful.
\end{thm}
\begin{proof}
Let $\Sigma$ be a smooth fan. We first assume that $\Sigma$ is saturated and unimodular. By the Hodge Isomorphism Theorem~\ref{thm:Hodge_isomorphism} for smooth tropical fans, for any cone $\eta \in \Sigma$, $A^\bul(\Sigma^\eta)$ verifies the Poincaré duality. Hence, we get an isomorphism $A^k(\Sigma^\eta) \simeq A^{d-k}(\Sigma^\eta)^\dual \simeq \MW_{d-k}(\Sigma^\eta)$. This implies that the map $\cycl_\eta\colon A^1(\Sigma^\eta) \to \Div(\Sigma^\eta)$ is bijective. We infer from the previous theorem that $\Sigma$ is both principal and div-faithful at $\eta$. This being the case for any $\eta \in \Sigma$, the theorem follows.

We now argue that the saturation condition is not needed by Remark \ref{rem:independent_exterior_lattice}. What is more, we do not need to require unimodularity, nor even simpliciality: Theorem \ref{thm:shellability_factorization} combined with Theorem \ref{thm:examples_shellable} \ref{thm:examples_shellable:div-faithful} and \ref{thm:examples_shellable:principal} imply that div-faithfulness and principality are properties of the support provided the support is locally irreducible. This finishes the proof of our theorem.
\end{proof}

\section{Tropical shellability} \label{sec:operations}

As in the previous section, $\Sigma$ will be a tropical fan of pure dimension $d$ in $N_\R \simeq \R^n$ for some natural number $n$. We moreover fix an orientation $(\nu_\sigma)_{\sigma \in \Sigma}$, and consider the canonical compactification $\comp\Sigma$ of $\Sigma$ with a compatible orientation.

\subsection{Tropical modifications} \label{subsec:tropical_modification} We start by recalling the definition of tropical modifications and introducing some variants of it in context related to canonical compactifications. A survey of results and references related to the concept can be found in~\cite{Kal15}.

Let $f\colon \supp\Sigma \to \R$ be a conewise integral linear function on $\Sigma$. Assume that the principal divisor $\div(f)$ is reduced. Denote by $\Delta$ the corresponding subfan of $\Sigma$. We allow the case the divisor $\div(f)$ is trivial, in which case, $\Delta$ will be empty. In this section, we will define three types of tropical modifications defined by $\Delta$, or more precisely by $f$, that we call \emph{open, closed}, and \emph{extended}, as described below.

\subsubsection{Open tropical modification}
We start by defining the open tropical modification of $\Sigma$ with respect to $f$. This will be a fan in $\~N_\R \simeq \R^{n+1}=\R^{n} \times \R$, for the lattice $\~N := N \times \Z$, that we will denote by $\tropmod{f}{\Sigma}$.

Consider the graph of $f$ which is the map $\Gamma_f$ defined as
\[ \begin{array}{rccc}
\Gamma=\Gamma_f\colon & \supp\Sigma & \longrightarrow & \~N_\R=N_\R \times \R, \\
                      & x           & \longmapsto     & (x, f(x)).
\end{array} \]
For each cone $\sigma$ of $\Sigma$, we consider the cone $\basetm\sigma$ in $\~N_\R$ which is the image of $\sigma$ by $\Gamma_f$, \ie, $\basetm\sigma:= \Gamma_f(\sigma) \subset \~N_\R$. Moreover, to each face $\delta$ of $\Delta$, we associate the face $\uptm\delta := \basetm\delta + \R_+ \etm$, where $\etm = (0, 1)\in N_\R \times \R$. Here $0$ in $\etm=(0,1)$ refers to the origin in $N_\R$.

The \emph{open tropical modification of\/ $\Sigma$ along $\Delta$ with respect to $f$}, or simply the tropical modification of $\Sigma$ along $\Delta$ if the other terms are understood from the context, is the fan $\tropmod{f}{\Sigma}$ in $\~N_\R\simeq \R^{n+1}$ defined as
\[\tropmod{f}{\Sigma} := \bigl\{\,\basetm\sigma \mid \sigma\in\Sigma\,\bigr\} \cup \bigl\{\,\uptm\delta \mid \delta\in\Delta\,\bigr\}.
\]

\begin{prop} \label{prop:balancing_trop_modif}
The tropical modification $\tropmod{f}{\Sigma}$ is a tropical fan. Moreover, we have a natural projection map
\[\prtm \colon \supp{\tropmod{f}{\Sigma}} \to \supp\Sigma\]
which is conewise integral linear.
\end{prop}

Before going through the proof, we make some remarks and introduce some notations. First, we observe that the fan $\tropmod{f}\Sigma$ is rational with respect to the lattice $\~N$, and that for each face $\sigma \in \Sigma$, the lattice $\~N_{\basetm\sigma}$ can be identified with the image $\Gamma_f(N_\sigma)$. This shows that for an inclusion of cones $\tau \ssubface \sigma$, we can pick
\[\nvect_{\basetm\sigma/\basetm\tau} = (\nvect_{\sigma/\tau},f_\sigma(\nvect_{\sigma/\tau})) \in \~N = N \times \Z,\]
where, as before, $f_\sigma \in N_\sigma^\dual$ denotes the linear form which coincides with $f$ on $\sigma$. Second, we observe that for each $\delta\in\Delta$, the lattice $\~N_{\uptm\delta}$ can be identified with $\~N_{\basetm\delta}\times\Z \simeq N_\delta \times \Z$, and that we can choose $\nvect_{\uptm\delta/\basetm\delta} = \etm$. Finally, for an inclusion of cones $\tau \ssupface \delta$ in $\Delta$, we can set $\nvect_{\uptm\tau/\uptm\delta} = (\nvect_{\tau/\delta},0) \in \~N$.

\begin{proof}[Proof of Proposition~\ref{prop:balancing_trop_modif}] We need to prove the balancing condition around any codimension one face of $\tropmod{f}\Sigma$. These are of two kinds, namely faces of the form $\basetm\tau$ for $\tau \in \Sigma_{d-1}$ and those of the form $\uptm\delta$ for $\delta \in \Delta_{d-2}$.

\smallskip
First let $\tau \in \Sigma_{d-1}$, and consider the codimension one face $\basetm\tau$ of $\tropmod{f}\Sigma$. Two cases happen:

\begin{itemize}[leftmargin=1em]
\item Either $\tau \in \Delta$, in which case the $d$-dimensional faces around $\basetm\tau$ are of the form $\basetm\sigma$ for $\sigma \ssupface \tau$ in $\Sigma$ as well as the face $\uptm\tau$. The balancing condition in this case amounts to showing that the vector
\begin{align*}
\nvect_{\uptm\tau/\tau} + \sum_{\sigma\ssupface \tau\\ \sigma\in \Sigma} \nvect_{\basetm\sigma/\basetm\tau}
  = \etm + \sum_{\sigma\ssupface \tau\\ \sigma\in \Sigma} (\nvect_{\sigma/\tau}, f_\sigma(\nvect_{\sigma/\tau}))
  = (0,1)+\Bigl(\,\sum_{\sigma\ssupface \tau\\ \sigma\in \Sigma} \nvect_{\sigma/\tau}, \sum_{\sigma\ssupface \tau\\ \sigma\in \Sigma} f_\sigma(\nvect_{\sigma/\tau}) \Bigr)\end{align*}
belongs to $\~N_{\basetm\tau}$. Since $\ord_\tau(f) =1$, the term on the right hand side of the above equality becomes equal to
\[ \Bigl(\,\sum_{\sigma\ssupface \tau\\ \sigma\in \Sigma} \nvect_{\sigma/\tau}, f_\tau\bigl(\sum_{\sigma\ssupface \tau\\ \sigma\in \Sigma} \nvect_{\sigma/\tau}\bigr) \Bigr),\]
which by balancing condition at $\tau$ in $\Sigma$ belongs to $\~N_{\basetm\tau}= \Gamma_f(N_\tau)$.

\item Or we have $\tau \notin \Delta$, \ie, $\ord_\tau(f)=0$. In this case, the facets around $\basetm\tau$ are of the form $\basetm\sigma$ for $\sigma \ssupface \tau$ in $\Sigma$, and we get
\begin{align*}
\sum_{\sigma\ssupface \tau\\ \sigma\in \Sigma} \nvect_{\basetm\sigma/\basetm\tau}
  = \sum_{\sigma\ssupface \tau\\ \sigma\in \Sigma} (\nvect_{\sigma/\tau}, f_\sigma(\nvect_{\sigma/\tau}))
  =\Bigl(\sum_{\sigma\ssupface \tau \\ \sigma\in \Sigma} \nvect_{\sigma/\tau}, \sum_{\sigma\ssupface \tau \\ \sigma\in \Sigma} f_\sigma(\nvect_{\sigma/\tau}) \Bigr)
  = \Bigl(\sum_{\sigma\ssupface \tau\\ \sigma\in \Sigma} \nvect_{\sigma/\tau}, f_\tau\bigl(\sum_{\sigma\ssupface \tau\\ \sigma\in \Sigma} \nvect_{\sigma/\tau}\bigr) \Bigr),
\end{align*}
which again belongs to $\~N_{\basetm\tau}= \Gamma_f(N_\tau)$.
\end{itemize}

It remains to check the balancing condition around a codimension one face of the form $\uptm\delta$ in $\tropmod{f}\Sigma$ with $\delta \in \Delta_{d-2}$. Facets around $\uptm\delta$ are all the cones $\uptm\tau$ for $\tau \in \Delta_{d-1}$ and $\tau \ssupface \delta$. Using the balancing condition in $\Delta =\div(f)$ around $\delta$, Proposition~\ref{prop:balancing_divisor}, we see that the sum
\[\sum_{\tau\ssupface \delta\\
\tau\in \Delta} \nvect_{\uptm\tau/\uptm\delta} = \Bigl(\sum_{\tau\ssupface \delta \\ \tau \in \Delta} \nvect_{\tau/\delta}, 0\Bigr)\]
belongs to $N_\delta \times \Z \subseteq \~N_{\uptm\delta}$, and the assertion follows.

\smallskip
The second statement is straightforward.
\end{proof}

We naturally endow $\tropmod{f}{\Sigma}$ with the orientation induced by the one of $\Sigma$ on faces of the form $\basetm\sigma$, for $\sigma\in\Sigma$, and extend it so that $\sign(\basetm\delta, \uptm\delta) = +1$ for any $\delta\in\Delta$.

\subsubsection{Closed tropical modification} The closed tropical modification lives in $N_\R \times \eR$, and is defined as the closure of the open tropical modification.

We still denote by $\basetm\sigma$ the image in $\~N_\R = N_\R \times \R$ of a cone $\sigma$ of $\Sigma$. For each $\delta \in \Delta$, we consider the cone $\uptm\delta = \basetm\delta + \R_+\etm $ defined in the previous section, and denote by $\compuptm\delta$ the closure of $\uptm\delta$ in $\R^n \times \eR$. Moreover, for any $\delta \in \Delta$, we denote by $\inftm\delta$ the extended polyhedron $\basetm\delta + \infty\cdot\etm$, which is also equal to $\compuptm\delta \setminus (N_\R \times \R)$.

\smallskip
The \emph{closed tropical modification of\/ $\Sigma$ with respect to $f$} is the extended polyhedral complex in $N_\R \times \eR$ which consists of faces $\{\basetm\sigma \mid \sigma\in\Sigma\} \cup \{\compuptm\delta \mid \delta\in\Delta\} \cup \{\inftm\delta \mid \delta\in\Delta \}$. We denote it by $\ctropmod{f}{\Sigma}$. We extend the orientation induced by the one on $\tropmod{f}{\Sigma}$, adding the condition $\sign(\inftm\delta, \uptm\delta) = -1$ for any $\delta\in\Delta$.

\subsubsection{Extended tropical modification} We now define the extended tropical modification. Let $\comp\Sigma$ be the canonical compactification of $\Sigma$, and let $\comp\Delta$ be the compactification of $\Delta$ inside $\comp\Sigma$. Then the \emph{closed extended tropical modification of\/ $\comp\Sigma$ with respect to $f$} denoted by $\ctropmod{f}{\comp\Sigma}$ is defined as the canonical compactification of the fan $\tropmod{f}{\Sigma}$. The set of faces of $\ctropmod{f}{\comp\Sigma}$ is exactly
\[ \{\basetm\sigma \mid \sigma\in\comp\Sigma\} \cup \{\compuptm\delta \mid \delta\in\comp\Delta\} \cup \{\inftm\delta \mid \delta\in\comp\Delta \}, \]
where $\basetm\sigma, \compuptm\delta$ and $\inftm\delta$ are defined as above with the map $\Gamma_f$ extended by continuity. The projection $\prtm$ extends to this context.

\smallskip
We also define the \emph{open extended tropical modification of $\comp\Sigma$ with respect to $f$} denoted by $\tropmod{f}{\comp\Sigma}$ as follows. Let $\~\Sigma$ be the tropical modification $\tropmod{f}{\Sigma}$. Then, $\tropmod{f}{\comp\Sigma}$ is defined as the restriction of $\ctropmod{f}{\comp\Sigma}$ to the space $\TP_{\~\Sigma} \setminus \TP_{\~\Sigma_\infty^{\etm}}$. The set of faces of $\tropmod{f}{\comp\Sigma}$ is exactly
\[ \{\basetm\sigma \mid \sigma\in\comp\Sigma\} \cup \{\uptm\delta \mid \delta\in\comp\Delta\}. \]

\subsubsection{Degenerate tropical modifications} In the case where $\Delta$ is trivial, the tropical modification of $\Sigma$ with respect to $f$ still has a meaning. In this case, the tropical modification is called \emph{degenerate}. The closed and the open tropical modifications coincide: they become both equal to the graph of $f$. The faces of the tropical modification are in one-to-one correspondance with faces of $\Sigma$. However, unless $f$ is a linear form on $\Sigma$, the tropical modification is not isomorphic to $\Sigma$. We refer to Example \ref{ex:nontrivial_degenerate_tropical_modification} in Section~\ref{sec:examples} which explains this phenomenon. On the contrary, if $f$ is integral linear, then we just obtain the image of $\Sigma$ by a linear map, and in this case, the fans $\Sigma$ and $\tropmod{f}{\Sigma}$ are isomorphic. In particular, if $\Sigma$ is div-faithful at $\conezero$, then the vanishing of $\div(f)$ implies that $f$ is linear, and in this case, the tropical modification becomes isomorphic to $\Sigma$.

\smallskip
Unless otherwise stated, in this article we allow tropical modifications to be degenerate.

\subsubsection{The case of tropical fans which are both principal and div-faithful at $\conezero$}

Assume $\Sigma$ is a tropical fan which is principal and div-faithful at the cone $\conezero \in \Sigma$. In this case any divisor $\Delta$ is the divisor of a conewise integral linear function $f$ on $\Sigma$, and in addition, if $g$ is another integral conewise linear function on $\Sigma$ such that $\div(g) = \Delta$, then $f-g$ is linear. Therefore, when $\Delta$ is reduced, the two tropical modifications $\tropmod{f}{\Sigma}$ and $\tropmod{g}{\Sigma}$ are isomorphic via the affine map which sends the point $x \in \R^{n+1}$ to the point $x + \ell(\prtm(x))\etm$ where $\ell \in M$ is a linear form restricting to $f-g$ on $\supp{\Sigma}$. Hence, working modulo isomorphisms, we can talk about \emph{the open tropical modification of\/ $\Sigma$ along $\Delta$}. We denote it by $\tropmod{\Delta}{\Sigma}$. The same applies to the extended setting and leads for instance to what we call \emph{the closed extended tropical modification of\/ $\comp\Sigma$ along $\comp\Delta$} that we denote by $\ctropmod{\comp\Delta}{\comp\Sigma}$.

\subsubsection{Star fans of a tropical modification} \label{sec:star_fans_tropical_modifications}

Let $\Sigma$ be any fan and let $\sigma \in \Sigma$. Let $f$ be a conewise linear map on $\Sigma$. Then $f$ induces a conewise linear map $f^\sigma$ on the star fan $\Sigma^\sigma$ defined as follows. Let $\ell \in M_\R$ be a linear map that coincide with $f$ on $\sigma$. Then $f-\ell$ is zero on $\sigma$ and we set $f^\sigma = \pi_*(f-\ell)$ where $\pi\colon N_\R \to N_\R^\sigma$ is the natural projection. Then $f^\sigma$ is a conewise linear map on $\Sigma^\sigma$. Note that although $f^\sigma$ depends on the choice of $\ell$, it is well-defined up to an element of $M_{\sigma, \R} =(N^\sigma_\R)^\dual$. This is enough for our purpose. That is why, abusing the terminology, we sometimes call $f^\sigma$ \emph{the} conewise linear map induced by $f$ on $\Sigma^\sigma$. We also note that in the case $\Sigma$ is rational and $f$ is integral, by choosing $\ell \in M$ we can ensure that $f^\sigma$ is conewise integral linear on $\Sigma^\sigma$.

\begin{prop} \label{prop:star_fans_tropical_modifications}
Let $\Sigma$ be a tropical fan and let $f$ be a conewise integral linear map on $\Sigma$. Assume the divisor $\div(f)$ is reduced and denote it by $\Delta$. Set $\~\Sigma = \tropmod{f}{\Sigma}$. Then we have the following description of the star fans of\/ $\~\Sigma$.
\begin{itemize}
\item If $\delta \in \Delta$, then $\~\Sigma^{\uptm\delta} \simeq \Delta^\delta$.
\item If $\sigma \in \Delta$, then $\~\Sigma^{\basetm\sigma} \simeq \tropmod{f^\sigma}{\Sigma^\sigma}$ where $f^\sigma$ is the conewise integral linear map induced by $f$ on $\Sigma^\sigma$.
\item If $\sigma \in \Sigma\setminus\Delta$, then once again we have $\~\Sigma^{\basetm\sigma} \simeq \tropmod{f^\sigma}{\Sigma^\sigma}$. However, this time the tropical modification is degenerate. In particular, if\/ $\Sigma$ is div-faithful (at $\sigma$), then $\~\Sigma^{\basetm\sigma}$ is isomorphic to $\Sigma^\sigma$.
\end{itemize}
\end{prop}

\begin{proof}
The proof is a direct verification.
\end{proof}

\subsection{Star subdivision}

Let $\Sigma$ be any rational fan. Let $\sigma \in \Sigma$ be any cone of $\Sigma$ of dimension at least one. Let $\rho$ be a rational ray generated by a vector in the relative interior of $\sigma$. The \emph{blow-up of\/ $\Sigma$ along $\rho$}, also called the \emph{star subdivision of\/ $\Sigma$ along $\rho$}, consists in replacing any face $\sigma\vee\eta'$ for $\eta'\sim\sigma$ by the faces of the form $\tau\vee\rho \vee \eta'$ for any proper subface $\tau$ of $\sigma$. We obtain a new fan with the same support which we denote by $\Sigma_{(\rho)}$. On the contrary, if $\Sigma$ is the blow-up of some fan $\Sigma'$ along a ray $\rho$, then $\Sigma'$ is called the \emph{blow-down of\/ $\Sigma$ along $\rho$}. By an abuse of the terminology we say a fan $\~\Sigma$ is obtained \emph{by a blow-up of\/ $\Sigma$ along $\sigma$} if there exists a ray $\rho$ in the relative interior of $\sigma$ as above so that $\~\Sigma$ coincides with $\Sigma_{(\rho)}$.

\smallskip
If $\Sigma$ is unimodular, a blow-up of $\Sigma$ along $\sigma$ is called \emph{unimodular} if $\Sigma_{(\rho)}$ is still unimodular. Such a blow-up is in fact unique. Indeed, for any face $\sigma$ with rays $\varrho_1, \dots, \varrho_k$, there is exactly one ray $\rho$ such that the blow-up along $\rho$ is unimodular. This ray is generated by $\e_{\varrho_1} + \dots + \e_{\varrho_k}$. Via the link to toric geometry, the unimodular blow-up $\~\Sigma$ of $\Sigma$ along $\sigma$ corresponds to the blow-up of the toric variety $\P_\Sigma$ along the closure $D^\sigma$ of the torus orbit $T^\sigma$ associated to $\sigma$. For this reason, we denote this blow-up by $\Bl{\Sigma}{\sigma}$.

\subsection{$\Csh$-shellability}

In the previous sections, we defined two types of operations on tropical fans: tropical modifications and star subdivisions. In this section, we define a notion of shellability using these operations. Roughly speaking, a tropical fan is shellable if it can be obtained from a collection of \emph{basic tropical fans} by only using the above operations.

The idea behind the definition is that if a property $P$ holds for our basic tropical fans and if, in addition, this property happens to be preserved by the above operations, we then obtain a wide collection of tropical fans all verifying the property $P$. This happens in practice for many examples of properties which will be discussed later in the paper.

\smallskip
Recall that by our convention from Section~\ref{sec:prel}, we work with fans modulo isomorphisms, so that we can talk about the set of isomorphism classes of rational fans.

\begin{defi}[Tropical shellability] \label{def:shellability}
Let $\Csh$ be a class of tropical fans (or more precisely, a set of isomorphism classes of tropical fans). Let $\Bsh$ be a subset of $\Csh$ that we call the \emph{base set}. The set of \emph{tropically $\Csh$-shellable fans over $\Bsh$} denoted by $\Sh_\Csh(\Bsh)$ is defined as the smallest class of tropical fans all in $\Csh$ which verifies the following properties.

\begin{itemize}
\item (Any element in the base set is tropically $\Csh$-shellable) We have $\Bsh \subseteq \Sh_\Csh(\Bsh)$.
\item (Closeness under products) If $\Sigma, \Sigma' \in \Sh_\Csh(\Bsh)$ and the product $\Sigma \times \Sigma'$ belongs to $\Csh$, then $\Sigma \times \Sigma' \in \Sh_\Csh(\Bsh)$.
\item (Closeness under tropical modifications along a tropically shellable divisor) If $\Sigma \in \Sh_\Csh(\Bsh)$ and if $f$ is a conewise integral linear function on $\Sigma$ such that either the divisor $\div(f)$ is reduced and $\div(f)\in \Sh_\Csh(\Bsh)$, or $\div(f)$ is trivial, and if in addition $\tropmod{f}{\Sigma} \in \Csh$, then $\tropmod{f}{\Sigma} \in \Sh_\Csh(\Bsh)$.
\item (Closeness under blow-ups and blow-downs with shellable center) If $\Sigma \in \Csh$, for any cone $\sigma \in \Sigma$ and any ray $\rho$ in the relative interior of $\sigma$ which verify $\Sigma^\sigma \in \Sh_\Csh(\Bsh)$ and $\Sigma_{(\rho)} \in \Csh$, we have
\[ \Sigma \in \Sh_\Csh(\Bsh) \Longleftrightarrow \Sigma_{(\rho)} \in \Sh_\Csh(\Bsh). \]
\end{itemize}
If $\Csh$ is the set of all tropical fans, we just write $\Sh(\Bsh)$.
\end{defi}

We note that the category of tropical fans is stable by the above list of operations. The main purpose of restricting the operations to be performed inside $\Csh$ is to forbid arbitrary blow-ups and blow-downs. The main examples of classes $\Csh$ of tropical fans which are of interest to us are \emph{all}, \resp \emph{simplicial}, \resp \emph{unimodular}, \resp \emph{unimodular quasi-projective}, tropical fans that we denote by $\tropf, \simp, \unim$, and $\uqproj$, respectively. We call them the \emph{standard classes} of tropical fans. These classes only constrain blow-ups and blow-downs since they are all stable by products and by tropical modifications (see Theorem \ref{thm:standard_properties}).

\smallskip
An important example of the base set is the set $\Bsh_0=\{\{\conezero\}, \Lambda\}$ where $\Lambda$ is the complete fan in $\R$ with three cones $\conezero, \R_{\geq 0}$, and $\R_{\leq 0}$. We will later see that $\Sh(\Bsh_0)$ contains many interesting fans.

\begin{defi}[Shellable fans]
A fan $\Sigma$ is called \emph{shellable in $\Csh$} if $\Sigma \in \Sh_\Csh(\Bsh_0)$. If $\Csh = \tropf$, we simply say that $\Sigma$ is shellable.
\end{defi}

\begin{defi}[Shellable properties] Let $\Csh$ be a class of tropical fans.
\begin{itemize}
\item A subclass $\Ssh$ of $\Csh$ is called \emph{stable by shellability in $\Csh$} if $\Sh_\Csh(\Ssh) = \Ssh$.

\item By extension, if $P$ is a predicate on tropical fans of $\Csh$, then $P$ is called \emph{stable by shellability in $\Csh$}, or simply \emph{shellable}, if the subclass of fans of $\Csh$ verifying $P$ is stable by shellability. Equivalently, $P$ is shellable if for any subclass $\Bsh$ of $\Csh$ such that all elements of $\Bsh$ verifies $P$, then all elements of $\Sh_\Csh(\Bsh)$ verifies $P$. \qedhere
\end{itemize}
\end{defi}

\subsection{Stellar-stability} A class $\Ssh$ of tropical fans is called \emph{stellar-stable} if for any $\Sigma \in \Ssh$ and any $\sigma \in \Sigma$, the star fan $\Sigma^\sigma$ also belongs to $\Ssh$. For instance, the four standard classes $\tropf, \simp, \unim$, and $\uqproj$, as well as $\Bsh_0$ are stellar-stable (\cf Theorem \ref{thm:standard_properties}). Recall that a predicate $P$ on fans is called local or stellar-stable if $P(\Sigma) \Rightarrow P(\Sigma^\sigma)$ for any $\sigma \in \Sigma$.

There is a natural way to construct a stellar-stable predicate from an arbitrary one. If $P$ is any predicate on tropical fans, then we denote by $\Ploc P$ the predicate
\[ \Ploc P(\Sigma)\colon\quad \forall\sigma\in\Sigma,\ P(\Sigma^\sigma). \]
For instance, local irreducibility, div-faithfulness, and principality are all defined in this way relative to the corresponding property at $\conezero$. The same holds for the smoothness property where $P$ is the Poincaré duality of the fan at cone $\conezero$. As we will see later in Lemma \ref{lem:shellability_meta_lemma}, viewing them this way allows to considerably simplify the proofs we give of their shellability.

\subsection{Properties of standard classes} The following theorem summarizes several nice properties enjoyed by the four standard classes of tropical fans that we introduced in the previous section.

\begin{thm} \label{thm:standard_properties}
Let $\Csh$ be one of the four standard classes of tropical fans. Then $\Csh$ verifies the following.
\begin{itemize}
\item \emph{(Closeness under products)} If\/ $\Sigma$ and $\Sigma'$ belong to $\Csh$, then we have $\Sigma \times \Sigma' \in \Csh$.
\item \emph{(Closeness under tropical modifications)} If\/ $\Sigma \in \Csh$ and if $f$ is a conewise integral linear function on $\Sigma$ such that either $\div(f)$ is reduced, or it is trivial, then $\tropmod{f}{\Sigma} \in \Csh$.
\item \emph{(Containment of the basic fans)} Both fans of $\Bsh_0$ are in $\Csh$.
\item \emph{(Stellar-stability)} The class $\Csh$ is stellar-stable.
\item \emph{(Existence of unimodular quasi-projective subdivisions)} Any fan in $\Csh$ has a subdivision in $\Csh$ which is unimodular and quasi-projective.
\item \emph{(Weak factorization)} Let $\Sigma$ and $\Sigma'$ be two fans in $\Csh$ with the same support. Then there exists a sequence of fans $\Sigma=\Sigma^0, \Sigma^1, \Sigma^2, \dots, \Sigma^k=\Sigma'$ in $\Csh$ such that for any $i\leq k-1$, $\Sigma_{i+1}$ is obtained from $\Sigma_i$ by performing a blow-up or a blow-down.
\end{itemize}
\end{thm}

\begin{proof}
We only sketch the proof here.
\begin{itemize}[leftmargin=0pt, itemindent=1em]
\item (Closeness under products) For $\tropf$ this is stated by Proposition \ref{prop:trop_normal_product}. Moreover, simpliciality and unimodularity are preserved by taking products. Let us justify that the product of two quasi-projective fans $\Sigma$ and $\Sigma'$ is quasi-projective. Let $\pi\colon \Sigma \times \Sigma' \to \Sigma$ and $\pi'\colon \Sigma \times \Sigma' \to \Sigma'$ be the two natural projections. Let $f$ and $f'$ be two convex conewise linear functions on $\Sigma$ and $\Sigma'$, respectively. Then a direct calculation proves that $\pi^*(f)+\pi'^*(f')$ is convex on $\Sigma \times \Sigma'$, and so the product remains quasi-projective.

\item (Closeness under tropical modifications) For $\tropf$, this is Proposition \ref{prop:balancing_trop_modif}. The new cones of the form $\uptm\delta$ with $\delta$ in the divisor are simplicial, \resp unimodular, provided $\delta$ is simplicial, \resp unimodular. Hence, simpliciality and unimodularity are preserved by tropical modifications. For quasi-projectivity, if $g$ is a convex conewise linear function on a fan $\Sigma$ and if $\prtm$ is the projection associated to a tropical modification $\~\Sigma$ of $\Sigma$, then we claim that $\prtm^*(g)$ is convex on $\~\Sigma$. Take $\delta \in \Delta$ and let us prove that $\prtm^*(g)$ is convex around $\basetm\delta$. By the very definition, we know that there exists a linear form $\ell$ such that $g-\ell$ is zero on $\delta$ and is strictly positive on rays $\rho\sim\delta$. Moreover, let $\uptm{l}$ be a linear form which takes values one on $\etm$ and $0$ on $\basetm\delta$. Then, for a small enough positive real number $\varepsilon$, $\prtm^*(g)-\prtm^*(\ell)+\varepsilon\uptm{l}$ is zero on $\basetm\delta$, takes value $\varepsilon>0$ on $\etm$, and is strictly positive on rays $\basetm\rho$ for $\rho\sim\delta$. This proves that $\prtm^*(g)$ is convex around $\basetm\delta$. One can prove similarly that $g$ is convex around other faces of $\~\Sigma$. Hence $\~\Sigma$ is quasi-projective.

\item (Containment of the basic fans) This is trivial.

\item (Stellar-stability) For $\tropf$, stellar-stability follows from Proposition \ref{prop:trop_normal_support_local}. Simpliciality and unimodularity are local properties. For quasi-projectivity, a conewise linear convex function on a fan induces convex functions on the star fans around its faces.

\item (Existence of unimodular quasi-projective subdivisions) This is a well-known fact. We refer to Section 4 of \cite{AP-tht} for more details.

\item (Weak factorization) This last property is far from trivial. For $\simp$ and $\unim$, this is Theorem A of~\cite{Wlo97}, proved independently by Morelli~\cite{Mor96} and expanded by Abramovich-Matsuki-Rashid, see~\cite{AMR}.

For $\uqproj$, this can be obtained from relevant parts of~\cites{Wlo97, Mor96, AKMW} as discussed and generalized by Abramovich and Temkin in~\cite{AT19}*{Section 3}.  \qedhere
\end{itemize}
\end{proof}

\subsection{Shellability, support and factorization}

The following shows that in some cases of interest, shellability is only a property of the support and of the ambient lattice.

\begin{thm} \label{thm:shellability_factorization}
Let $\Csh$ be a standard class. Let $\Ssh$ be subclass of $\Csh$ which is both stellar-stable and stable by shellability. Then a tropical fan $\Sigma$ of\/ $\Csh$ is in $\Ssh$ if and only if any tropical fan of\/ $\Csh$ with the same support $\supp{\Sigma}$ and considered with the same lattice is in $\Ssh$.
\end{thm}

Notice that here we remember the ambient lattice. This is weaker than being a property of the support. For instance, if $\Sigma$ is not saturated in $N_\R$ and if $\Sigma'$ is the same fan considered with a different lattice $N' = \SF_1(\Sigma)$, then the above theorem does not imply that $\Sigma\in\Ssh$ if and only if $\Sigma'\in\Ssh$, even though $\Sigma\cap N = \Sigma'\cap N'$. For many applications, \eg, for properties depending on $\Sigma\cap N$ in the sense of Remark \ref{rem:independent_exterior_lattice}, this remark has no importance. However, for the study of the Chow ring, for instance for the shellability statement in Theorem \ref{thm:examples_shellable} \ref{thm:examples_shellable:PD_Chow_ring} in the next section, some counter-examples can be easily constructed (see Example \ref{ex:F1_vs_N}).

\begin{proof}
Let $\Sigma\in\Ssh$ and let $\Sigma'$ be another fan of $\Csh$ with the same support and the same lattice. By the weak factorization property of Theorem \ref{thm:standard_properties}, there exists a sequence of fans $\Sigma = \Sigma^0, \Sigma^1, \dots, \Sigma^{k-1}, \Sigma^k = \Sigma'$ all belonging to $\Csh$ such that $\Sigma^{i+1}$ is obtained from $\Sigma^i$ by performing a blow-up or a blow-down.

Let us prove that $\Sigma^1 \in \Ssh$. If $\Sigma^1$ is obtained from $\Sigma$ by blowing up the face $\eta \in \Sigma$, then we have $\Sigma^\eta\in\Ssh$ by stellar-stability. Hence $\Ssh$ is closed by blow-up along $\eta$, and we have $\Sigma^1 \in \Ssh$. Otherwise, this is $\Sigma$ which is obtained from $\Sigma^1$ by blowing up along a ray $\rho$ which is in the relative interior of some face $\eta'\in\Sigma^1$. In this case, notice that $(\Sigma^1)^{\eta'} = \Sigma^{\tau\vee\rho}$ where $\tau$ is any face of codimension one in $\eta'$. Once again, $(\Sigma^1)^{\eta'}$ belongs to $\Ssh$ and, by closeness under blow-downs, we get $\Sigma^1 \in \Ssh$.

Proceeding step by step and using the same argument for the next fans in the sequence, we obtain that $\Sigma' = \Sigma^k \in \Ssh$.
\end{proof}

\begin{remark} \label{rem:shellability_factorization}
Note that the proof is still valid if we just require $\Ssh$ to be stellar-stable and closed by blow-ups and blow-downs along faces whose star fan belongs to $\Ssh$.
\end{remark}

\subsection{Examples of shellable properties} The following theorem gives important examples of properties which are shellable in relevant classes of tropical fans.

\begin{thm} \label{thm:examples_shellable}
We have the following.
\begin{enumerate}
\item \label{thm:examples_shellable:normal} Normality is shellable.
\item \label{thm:examples_shellable:irreducible} Local irreducibility is shellable.
\item \label{thm:examples_shellable:div-faithful} Div-faithfulness is shellable.
\item \label{thm:examples_shellable:principal} Principality is shellable in the class of locally irreducible fans. More precisely, the property of being both principal and locally irreducible is shellable.
\item \label{thm:examples_shellable:PD_Chow_ring} Poincaré duality for the Chow ring is shellable in the class of div-faithful unimodular fans. More precisely, the set of div-faithful fans whose Chow rings verify the Poincaré duality is stable by shellability in $\unim$.
\item \label{thm:examples_shellable:HR_Chow_ring} Hodge-Riemann and Hard Lefschetz for the Chow ring are both shellable in $\uqproj$.
\item \label{thm:examples_shellable:smooth} Smoothness is shellable.
\end{enumerate}
\end{thm}

\begin{proof}
In order to prove this theorem, we will introduce in Section \ref{subsec:shellability_meta_lemma} a useful tool called \emph{shellability meta lemma} which allows to simplify the proofs of the different points in the definition of shellability. Then, we will be able to prove points \ref{thm:examples_shellable:normal} and \ref{thm:examples_shellable:irreducible} in Section \ref{subsec:normal_irreducible_shellable}.

The other points are treated later: \ref{thm:examples_shellable:div-faithful} and \ref{thm:examples_shellable:principal} are the contents of Theorems~\ref{thm:shellability_div-faithful} and \ref{thm:shellability_principal}, both proved in Section~\ref{sec:chow_ring}. Property \ref{thm:examples_shellable:PD_Chow_ring} is Theorem~\ref{thm:chow_shellable}, which is also established in the same section.

Point~\ref{thm:examples_shellable:HR_Chow_ring} is proved in~\cite{AP-tht}. The last point is the content of Theorem~\ref{thm:smooth_shellable} and will be proved in Section~\ref{sec:homology_tropical_modification}, using tools we will introduce in Section~\ref{sec:deligne}.
\end{proof}

All the properties of Theorem \ref{thm:examples_shellable} are verified by elements of $\Bsh_0$. Hence, all these properties are true in $\Sh_\uqproj(\Bsh_0)$. In Section \ref{subsec:Bergman_fans_shellable}, we prove that this class contains an interesting class of fans, namely generalized Bergman fans.

\subsection{A tool to prove shellability} \label{subsec:shellability_meta_lemma}

Checking all the axioms of shellability might be somehow tedious in general. The following meta lemma helps in practice to simplify the verification of these different points.

\begin{lemma}[Shellability meta lemma] \label{lem:shellability_meta_lemma}
Let $\Csh$ be one of the four standard classes. Let $P$ be a predicate on elements of\/ $\Csh$. Assume that the elements of $\Bsh_0$ verify $P$.

Let $\Sigma$ be an arbitrary fan in $\Csh$ such that for any face $\sigma \neq \conezero$ in $\Sigma$, the star fan $\Sigma^\sigma$ verifies $\Ploc P$, and such that at least one of the following points is verified.
\begin{itemize}
\item $\Sigma$ is the product of two unimodular quasi-projective fans verifying $\Ploc P$.
\item $\Sigma$ is the tropical modification of a unimodular quasi-projective fan verifying $\Ploc P$ with respect to some conewise integral linear function $f$ such that either $\div(f)$ vanishes, or $\div(f)$ is reduced, unimodular, quasi-projective, and verifies $\Ploc P$.
\item $\Sigma$ is the blow-up along some ray of a fan verifying $\Ploc P$.
\item $\Sigma$ is the blow-down along some ray of a fan verifying $\Ploc P$.
\end{itemize}
If for any fan $\Sigma$ as above, the property $P$ is verified, then $\Ploc P$ is shellable in $\Csh$.

\smallskip
Moreover, if $P$ is a predicate only depending on the support of the fan, we can restrict ourselves to fans $\Sigma$ verifying one of the two first points.
\end{lemma}

Using this lemma, one can now just focus on the following goal: proving that $\Sigma$ verifies $P$ for each of the four cases. For instance, we have to prove that if $\Sigma$ is a product of two fans verifying $P^\star$, then $\Sigma$ verifies $P$. To do so, we can assume without loss of generality that the two factors (and thus $\Sigma$ itself) are unimodular and quasi-projective, and that every proper star fan of $\Sigma$ verifies $P^\star$. Most of the time, it is not really needed to assume so many properties to be verified by $\Sigma$, and we can proceed in more generality. This being said, we will often assume $\Sigma$ is unimodular for instance.

\smallskip
The rest of this section is devoted to the proof of the above lemma. First, notice that the very last part of the proposition is clear since blow-ups and blow-downs do not change the support of the fan.

Let $\Ssh \subseteq \Csh$ be the class of fans verifying $\Ploc P$. Since $\Ploc P$ is stellar-shellable, so is $\Ssh$. For any integer $n$, we use the notation $\Csh^{<n}$, \resp $\Csh^{\leq n}$, to denote the subset of $\Csh$ of fans of dimension less than $n$, \resp at most $n$. We define $\Ssh^{<n}$ and $\Ssh^{\leq n}$ similarly.

Let us prove that, under the assumption of the lemma, $\Ssh$ is stable by shellability. We prove it by induction on $n$. Assume that, for some integer $n$, $\Ssh^{< n}$ is stable by shellability in $\Csh^{< n}$. This is clear for $n=1$. We prove that $\Ssh^{\leq n}$ is stable by shellability inside $\Csh^{\leq n}$.

\subsubsection{Closeness under blow-ups and blow-downs} We verify that $\Ssh^{\leq n}$ is closed under blow-ups and blow-downs. Let $\Sigma$ be a fan in $\Csh^{\leq n}$. Let $\eta \in \Sigma$ be a face and let $\rho$ be the ray in the relative interior of $\eta$. Assume that $\Sigma_{(\rho)} \in \Csh^{\leq n}$. We need to show the equivalence
\[ \Sigma \in \Ssh^{\leq n} \Longleftrightarrow \Sigma_{(\rho)} \in \Ssh^{\leq n}. \]

We first prove if $\Sigma_{(\rho)} \in \Ssh^{\leq n}$, then $\Sigma \in \Ssh^{\leq n}$. For this, we compare the star fans of $\Sigma$ and $\Sigma_{(\rho)}$ as follows. Consider a face $\sigma$ of $\Sigma$ different from $\conezero$. There are three cases.

\begin{itemize}
\item First, assume that $\sigma$ is not comparable with $\eta$. In this case, the two star fans $\Sigma_{(\rho)}^\sigma$ and $\Sigma^\sigma$ are identical. Since $\Ssh$ is stellar-stable, the first one is in $\Ssh$ by assumption. Hence $\Sigma^\sigma\in\Ssh$.

\item Second, we assume that $\sigma \supface \eta$. In this case, we have $\Sigma^\sigma \simeq \Sigma_{(\rho)}^{\tau\vee\rho}$ where $\tau$ is any face of codimension one in $\sigma$ such that $\tau\wedge\eta \ssubface \eta$. Once again $\Sigma^\sigma \in \Ssh$.

\item Finally, assume that $\sigma$ and $\eta$ are comparable but $\sigma \not\supface \eta$. Denote by $\eta^\sigma$, \resp $\rho^\sigma$, the face corresponding to $\eta$, \resp to $\rho$, in $\Sigma^\sigma$. Then $\rho^\sigma$ is a ray in the relative interior $\eta^\sigma$. The star fan $\Sigma_{(\rho)}^\sigma$ is then naturally isomorphic to the blow-up star fan $\bigl(\Sigma^\sigma\bigr)_{(\rho^\sigma)}$. By assumption, $\Sigma_{(\rho)}^\sigma \in \Ssh^{<n}$. Moreover $\bigl(\Sigma^\sigma\bigr)^{\eta^\sigma} = \Sigma^{\sigma\vee\eta}$ which is in $\Ssh^{<n}$ by the second point above. Hence, applying the closeness by blow-down of $\Ssh^{<n}$, we deduce that $\Sigma^\sigma$ is also in $\Ssh^{<n} \subseteq \Ssh$.
\end{itemize}

In any case, $\Sigma^\sigma \in \Ssh$ for any face $\sigma\neq\conezero$. One can apply the assumption of the lemma to deduce that $\Sigma$ verifies $P$. Hence, $\Sigma$ verifies $\Ploc P$, thus $\Sigma \in \Ssh^{\leq n}$. This proves the direction $\Leftarrow$.

\smallskip
To prove the direction $\Rightarrow$, assume $\Sigma \in \Ssh^{\leq n}$. Take a cone $\zeta$ in $\Sigma_{(\rho)}$ different from $\conezero$. Apart from the faces which already appeared in the above case analysis for which the reversed argument applies, it remains to consider those faces $\zeta$ with $\rho \subface \zeta$. Denote by $\zeta-\rho$ the face of $\zeta$ of codimension one which does not contain $\rho$. Set $\sigma=(\zeta-\rho)\vee\eta$. Then we get
\[ \supp{\Sigma_{(\rho)}^\zeta} \simeq \supp{\Sigma^\sigma} \times \R^k, \]
where $k = \dims{\sigma}-\dims{\zeta} \in \{0, \dots, n-1\}$. Since $\Ssh^{<n}$ is stellar-stable and closed by blow-ups and blow-downs, we can apply Remark \ref{rem:shellability_factorization} and Theorem \ref{thm:shellability_factorization}: $\Sigma_{(\rho)}^\zeta$ is in $\Ssh^{<n}$ if and only if there exists a fan with the same support in $\Ssh^{<n}$. This is the case. Indeed, $\Sigma^\sigma \in \Ssh^{<n}$. Moreover, $\Lambda$ and $\{\conezero\}$ belong to $\Ssh$ by assumption. Since $\Ssh^{<n}$ is stable by shellability, we deduce that $\Sigma^\sigma \times \Lambda^k \in \Ssh^{<n}$. This last fan has the same support as $\Sigma_{(\rho)}^\zeta$. Hence, we infer that $\Sigma_{(\rho)}^\zeta \in \Ssh^{<n} \subseteq \Ssh$.

We have proved that $\Sigma_{(\rho)}^\zeta \in \Ssh$ for any nontrivial face $\zeta \in \Sigma_{(\rho)}$. As before, we apply the assumption of the lemma to deduce that $\Sigma_{(\rho)} \in \Ssh^{\leq n}$. Hence, $\Ssh^{\leq n}$ is stable by blow-ups and blow-downs. In particular, we can apply Theorem \ref{thm:shellability_factorization} in $\Ssh^{\leq n}$ for the rest of this proof.

\subsubsection{Closeness under products} We whish to prove that $\Ssh^{\leq n}$ is stable by products inside $\Csh^{\leq n}$. Let $\~\Sigma\in\Csh^{\leq n}$ be the product of two fans in $\Ssh^{\leq n}$. Let us prove that $\~\Sigma\in\Ssh^{\leq n}$.

Denote this two factors by $\~\Sigma^1$ and $\~\Sigma^2$. Let $\Sigma^1$, \resp $\Sigma^2$, be a unimodular quasi-projective subdivision in $\Csh$ of $\Sigma^1$, \resp of $\Sigma^2$, which exists by Theorem \ref{thm:standard_properties}. Then, clearly $\Sigma^1$ has the same support as $\~\Sigma^1$, and by Theorem \ref{thm:shellability_factorization}, we get $\Sigma^1\in\Ssh^{\leq n}$. In the same way, we obtain $\Sigma^2\in\Ssh^{\leq n}$. Set $\Sigma=\Sigma^1\times\Sigma^2$. Let $\sigma^1\times\sigma^2$ be a nontrivial face of $\Sigma$. Then
\[ \Sigma^{\sigma^1\times\sigma^2} \simeq (\Sigma^1)^{\sigma^1} \times (\Sigma^2)^{\sigma^2}. \]
By stellar-stability of $\Ssh$, both factors belong to $\Ssh^{<n}$. Hence, the stability by shellability of $\Ssh^{<n}$ implies that the product belongs to $\Ssh^{<n}$.

Therefore, for any nontrivial face $\sigma$ of $\Sigma$, we get $\Sigma^\sigma\in\Ssh$. Applying the assumption of the lemma, we deduce that $\Sigma$ verifies $P$ and thus $\Sigma \in \Ssh^{\leq n}$. Since $\Sigma$ and $\~\Sigma$ have the same support, we can apply Theorem \ref{thm:shellability_factorization} to deduce that $\~\Sigma\in\Ssh^{\leq n}$ as well. Thus, $\Ssh^{\leq n}$ is stable by products.

\subsubsection{Closeness under tropical modifications} Let $\Sigma'$ be a fan in $\Ssh^{\leq n}$. Let $f$ be a conewise integral linear function on $\Sigma'$ such that $\div(f)$ is reduced. Set $\~\Sigma' = \tropmod{f}{\Sigma'}$. Set $\Delta'=\div(f)$ and assume that $\Delta'$ is in $\Ssh^{\leq n}$ (by convention, we assume in this proof that $\emptyset \in \Ssh^{\leq n}$). We wish to prove that $\~\Sigma'$ is in $\Ssh^{\leq n}$.

Let $\Sigma$ be a unimodular quasi-projective subdivision of $\Sigma'$. Then $f$ is still conewise integral linear on $\Sigma$. Moreover, $\Delta:=\div(f)$, taken in $\Sigma$, is a unimodular quasi-projective subdivision of $\Delta'$. As for the case of the product, Theorem \ref{thm:shellability_factorization} implies that both $\Sigma$ and $\Delta$ are in $\Ssh^{\leq n}$. Set $\~\Sigma = \tropmod{f}{\Sigma}$. A face $\~\sigma$ of dimension $k>0$ in $\~\Sigma$ is of two kinds, either it is equal to $\uptm\delta$ for $\delta \in \Delta_{k-1}$ or it coincides with $\basetm\sigma$ for $\sigma \in \Sigma_{k}$. By Proposition \ref{prop:star_fans_tropical_modifications}, in the first case, the star fan $\~\Sigma^{\uptm\delta}$ is isomorphic to $\Delta^\delta$ and so belongs to $\Ssh$. So we can now assume that $\~\sigma = \basetm\sigma$ for a cone $\sigma \in \Sigma$. Then, $\~\Sigma^{\basetm\sigma}$ is the tropical modification of $\Sigma^\sigma$ along $\Delta^\sigma$ with respect to the function $f^\sigma$. By convention here we set $\Delta^\sigma = \emptyset$ if $\sigma\not\in\Delta$. Note that $\Sigma^\sigma$ and $\Delta^\sigma$ are in $\Ssh^{<n}$. Since $\Ssh^{<n}$ is stable by shellability, $\tropmod{f^\sigma}{\Sigma^\sigma} \in \Ssh^{<n}$. We infer again that $\~\Sigma^\sigma \in \Ssh$, as desired.

At this point we have verified that for any nontrivial cone $\~\sigma$ in $\~\Sigma$, the star fan $\~\Sigma^{\~\sigma}$ is in $\Ssh$. Using the assumption of the lemma, we deduce that $\~\Sigma$ verifies $P$ and so $\~\Sigma \in \Ssh^{\leq n}$. By Theorem \ref{thm:shellability_factorization}, we deduce that $\~\Sigma' \in \Ssh^{\leq n}$. Therefore, $\Ssh^{\leq n}$ is closed by tropical modifications.

\subsubsection{End of the proof} Finally, we have proved that $\Ssh^{\leq n}$ is stable by shellability in $\Csh^{\leq n}$. By induction, we deduce that $\Ssh$ is stable by shellability, \ie, $\Ploc P$ is shellable in $\Csh$. \qed

\subsection{Normality and local irreducibility are shellable} \label{subsec:normal_irreducible_shellable}

In this section, we prove Points \ref{thm:examples_shellable:normal} and \ref{thm:examples_shellable:irreducible} of Theorem \ref{thm:examples_shellable}, namely that normality and local irreducibility are shellable. Let us start with normality.

\smallskip
First recall that a tropical fan $\Sigma$ is shellable if and only if for any face $\eta$ of codimension one in $\Sigma$, $\Sigma^\eta$ is isomorphic to the tropical line in some $\R^k$, \ie, if $\Sigma^\eta$ has rays $\rho_0, \dots, \rho_k$, then the only linear relation of dependence between $\e_{\rho_0}, \dots, \e_{\rho_k}$ is (up to a scalar) $\e_{\rho_0}+\dots+\e_{\rho_k} = 0$. In particular, normality is stellar-stable and only depends on the support (see Proposition \ref{prop:trop_normal_support_local}).

Denote by $P$ the predicate for a tropical fan to be normal. Hence $P = \Ploc P$. Note that elements of $\Bsh_0$ are normal. Let $\Sigma$ be a fan verifying the condition of Lemma \ref{lem:shellability_meta_lemma}. In particular, if the dimension of $\Sigma$ is at least two, then for any face $\eta$ of codimension one in $\Sigma$, the star fan $\Sigma^\eta$ is normal. As a consequence $\Sigma$ verifies $P$.

It remains to treat the case where $\Sigma$ has dimension one. In this special situation, the only nontrivial case is when $\Sigma$ is the tropical modification of a normal fan, \ie, of the tropical line in $\R^k$ for some integer $k$. There are two cases.
\begin{itemize}
\item First, assume $\Sigma$ is a non-degenerate tropical modification of the tropical line in $\R^k$. Let $\e_0, \dots, \e_k, \etm$ be the unit vectors of the rays of $\Sigma$ with $\etm$ the special one, \ie, the one such that $\prtm(\etm)=0$ for the projection map $\prtm$ associated to the tropical modification. Let $a_0, \dots, a_k, \uptm a$ be some coefficients such that $a_0\e_0 + \dots + \uptm a\etm = 0$. Applying $\prtm$ we deduce that $a_0\prtm(\e_0) + \dots + a_k\prtm(\e_k)=0$. The normality implies that $a_0 = \dots = a_k$. Moreover, by the balancing condition on $\Sigma$, we get $\e_0 + \dots + \e_k = -\etm$. Thus, $(-a_0 + \uptm a)\etm = 0$. Hence, we should have $a_0 = \dots = a_k = \uptm a$, which proves that $\Sigma$ is normal.

\item Otherwise, $\Sigma$ is a degenerate tropical modification of a tropical line. One can prove in the same way that $\Sigma$ is normal in this case.
\end{itemize}
We have proved that normality fulfills the conditions of Lemma \ref{lem:shellability_meta_lemma}. Hence, normality is shellable.

\smallskip
Denote by $Q$ the property of being connected in codimension one. Recall that a fan $\Sigma$ is locally irreducible if and only if it is normal and verifies $\Ploc Q$ (this is Proposition \ref{thm:characterization_irreducible}). Thanks to Lemma \ref{lem:shellability_meta_lemma}, one can easily prove that $\Ploc Q$ is shellable. Since normality is also shellable, it is clear that local irreducibility is shellable. \qed

\subsection{Shellability of the generalized Bergman fans} \label{subsec:Bergman_fans_shellable}

Recall that we say a tropical fan $\Sigma$ is a generalized Bergman fan if it is isomorphic to a tropical fan with the same support as the Bergman fan of a matroid. In this section we prove the following theorem.

\begin{thm} \label{thm:bergman_shellable} The Bergman fan of a matroid is shellable in $\uqproj$. More generally, any unimodular quasi-projective generalized Bergman fan is shellable in $\uqproj$.
\end{thm}
Note that the Bergman fan of a matroid is quasi-projective. Indeed, the Bergman fan of any matroid is a subfan of the Bergman fan of a free matroid, which is projective. This last observation also implies that any complete fan is a generalized Bergman fan. The above theorem thus implies that projective fans are shellable.

\smallskip
The rest of this section is devoted to recalling basic definitions and properties regarding matroids and their Bergman fans, and then giving the proof of this theorem.

\subsubsection{Matroids}
We start by recalling basic definitions involving matroids and refer to relevant part of \cite{Oxl06} for more details.

A matroid can be defined in different equivalent ways, for example by specifying what is called its \emph{collection of independent sets}, or its \emph{collection of bases}, or its \emph{collection of flats}, or its \emph{collection of circuits}, or still by giving its \emph{rank function}. The data of any of these collections determines all the others. The definition of a matroid with respect to independent sets is the following.

\begin{defi}[Matroid: definition with respect to the family of independent sets] \label{defi:matroid}
A \emph{matroid} $\Ma$ is a pair $(E, \ind)$ consisting of a finite set $E$ called the \emph{ground set} and a collection $\ind$ of subsets of $E$ called the \emph{family of independent sets of\/ $\Ma$} which verifies the following axiomatic properties:
\begin{enumerate}
\item The empty set is an independent set: $\emptyset\in\ind$.
\item (Hereditary property) $\ind$ is stable under inclusion: if $J\subset I$ and $I\in\ind$, then $J\in\ind$.
\item (Augmentation property) for two elements $I, J\in\ind$, if $\Card J<\Card I$, then one can find an element $i$ in $I\setminus J$ such that $J+i\in\ind$. \label{defi:matroid:augmentation_property} \qedhere
\end{enumerate}
\end{defi}

An example of a such a pair $(E, \ind)$ is given by a collection of vectors $v_1, \dots, v_m$ in a finite dimensional vector space $H$ over a field $\k$. The corresponding matroid $\Ma$ has ground set $E = [m]$ and has independent sets $\ind$ consisting of all subsets $I \subseteq [m]$ verifying that the corresponding vectors $v_i$ for $i\in I$ are linearly independent. A matroid $\Ma$ of this form is said to be \emph{representable} (over $\k$). A recent result of Peter Nelson~\cite{Nel18} shows that \emph{almost any matroid is non-representable over any field}: denoting by $\epsilon_m$ the proportion of representable matroids on the ground set $[m]$ among all the matroids on this ground set, the numbers $\epsilon_m$ converge to zero when $m$ goes to infinity.

\smallskip
To a given matroid $\Ma=(E,\ind)$ we can associate the so-called \emph{rank function} $\rkm\colon 2^E \to \Z_{\geq 0}$ which is defined as follows. Consider a subset $A \subseteq E$. The \emph{rank of $A$} is then defined by taking the maximum of $\Card I$ over all independent sets $I$ which are included in $A$. The \emph{rank of\/ $\Ma$} is then the integer $\rkm(E)$. The rank function satisfies what is called the \emph{submodularity property} which is
\[\forall A, A' \subseteq E, \qquad \rkm(A\cap A') + \rkm(A \cup A') \leq \rkm(A) + \rkm(A').\]
The rank function $\rkm$ determines $\ind$ as the family of all subsets $I \subseteq E$ which verify $\rkm(I) = \Card{I}$.

A \emph{basis} of $\Ma$ by definition is a maximal independent set. The collection of bases of $\Ma$ is denoted by $\bases(\Ma)$. On the contrary, a \emph{circuit} of $\Ma$ is a minimal dependent set, where a set is \emph{dependent} if it is not independent. The collection of circuits of $\Ma$ is denoted by $\crct(\Ma)$.

The \emph{closure} $\cl(A)$ of a subset $A$ in $E$ is defined as
\[ \cl(A):=\{a \in E \textrm{ such that } \rkm(A + a)=\rkm(A)\}. \]
A \emph{flat} of $\Ma$ is a subset $F \subseteq E$ with $\cl(F) = F$. Flats are also equally called \emph{closed sets} and the collection of flats of the matroid $\Ma$ is denoted by $\Cl(\Ma)$.

For a representable matroid, given by a collection of vectors $v_1, \dots, v_m$ in a vector space $ H \simeq \k^n$, the above definitions coincide with the usual ones from linear algebra. For example, a flat $F$ of $\Ma$ is a subset of $[m]$ such that $\{v_i\}_{i\in F}$ is the set of all the vectors in the collection which live in a vector subspace of $H$.

An element $e$ in the ground set is a \emph{loop} if $\{e\}$ is not in $\ind$, that is, if $\rkm(\{e\})=0$. Two elements $e$ and $e'$ of $E$ are called \emph{parallel} if $\rkm(\{e, e'\}) = \rkm(\{e\}) = \rkm(\{e'\}) = 1$. A \emph{simple} matroid is a one which neither contains a loop nor parallel elements.

An element $e\in E$ is called a \emph{coloop} if $\rkm(E-e)=\rkm(E)-1$. This is equivalent to $E-e$ being a flat.

Finally, a \emph{proper element} in a matroid is an element which is neither a loop nor a coloop.

\subsubsection{Bergman fans} For a simple matroid $\Ma$ on a ground set $E$, the Bergman fan of $\Ma$ denoted by $\Sigma_\Ma$ is defined as follows.

Recall that for a subset $A \subseteq E$, we denote by $\e_A$ the sum $\sum_{i\in A} \e_i$. Here $\{\e_i\}_{i\in E}$ is the standard basis of $\R^E$. Let $N = \rquot{\Z^E}{\Z \e_E}$ and denote by $M$ the dual of $N$. By an abuse of the notation, we denote by the same notation $\e_A$ the projection in $N_\R$ of $\sum_{i\in A}\e_i$. Note in particular that $\e_E =0$.

The Bergman fan $\Sigma_\Ma$ of $\Ma$ is a rational fan in $N_\R$ of dimension $\rkm(\Ma) -1$ defined as follows. First, a \emph{flag of proper flats} $\Fl$ of $\Ma$ is a collection
\[\Fl\colon \quad \emptyset \neq F_1 \subsetneq F_2 \subsetneq \dots \subsetneq F_{\ell} \neq E\]
consisting of flats $F_1, \dots, F_\ell$ of $\Ma$. The number $\ell$ is called the \emph{length} of $\Fl$.

To such a flag $\Fl$, we associate the rational cone $\sigma_\Fl \subseteq N_\R$ generated by the vectors $\e_{F_1}, \e_{F_2}, \dots, \e_{F_\ell}$, that is
\[\sigma_\Fl := \Bigl\{ \lambda_1 \e_{F_1} + \dots + \lambda_\ell \e_{F_\ell} \Bigst \lambda_1, \dots, \lambda_\ell \geq0\Bigr\}.\]
Obviously, the dimension of $\sigma_\Fl$ is equal to the length of $\Fl$.

The \emph{Bergman fan of $\Ma$} is the fan consisting of all the cones $\sigma_\Fl$, $\Fl$ a flag of proper flats of $\Ma$, \ie,
\[\Sigma_\Ma := \Bigl\{\, \sigma_\Fl \Bigst \Fl \textrm{ flag of proper flats of }\Ma\,\Bigr\}.\]
It is straightforward using the properties of matroids to check that $\Sigma_\Ma$ is a tropical fan of pure dimension $\rkm(\Ma)-1$. A \emph{generalized Bergman fan} is any fan isomorphic to a fan of support $\supp{\Sigma_\Ma}$ (considered with the same lattice) for some matroid $\Ma$.

\smallskip
We now mention some basic properties and examples of generalized Bergman fans. Let $\Ma$ be a matroid on the ground set $E=[m]$ and denote by $\crct(\Ma)$ the circuits of $\Ma$. The support of $\Sigma_\Ma$ can be described as follows, see~\cite{AK06} for more details. Let $\hat\Sigma_{\Ma}$ be the set of points $x=(x_i)_{i \in E}\in \R^E$ such that for every circuit $C\in\crct(\Ma)$, the minimum of $x_i$ for $i\in C$ is achieved at least twice. Note that if $x\in\hat\Sigma_{\Ma}$, then $x+\lambda\ones\in\hat\Sigma_\Ma$ for all $\lambda\in\R$.

\begin{prop}[Ardila-Klivans~\cite{AK06}] \label{prop:support_bergman} The support of\/ $\Sigma_\Ma$ coincides with the projection of\/ $\hat\Sigma_\Ma$ in $\rquot{\R^E\!}{\R\ones}$.
\end{prop}

As we previously said, any complete rational fan $\Sigma$ is a generalized Bergman fan. To see this, let $U_{r+1}$ be the uniform matroid on the ground set $E=[r+1]$, with collection of independent sets $\ind = 2^{[r+1]}$. In this case, the support of $\Sigma_{U_{r+1}}$ is the full space $\R^r \simeq \rquot{\R^{r+1}}{\R \e_{E}}$.

\smallskip
Let now $\Ma$ be a realizable matroid on the ground set $E=[m]$ defined by a collection of non-zero vectors $v_1, \dots, v_m$ in a vector space $V$ over the field $\k$. We can view $v_1, \dots, v_m$ as elements of the dual vector space of $V^\dual$ and in this case, we get a collection of hyperplanes $H_1, \dots, H_m \subset V^\dual$ with $H_j = \ker\Bigl( v_j \colon V^\dual \to \k\Bigr)$, $j\in[m]$. This leads to a hyperplane arrangement in the projective space $\P(V^\dual)$ again denoted by an abuse of the notation $H_1, \dots, H_m \subseteq \P(V^\dual)$.

Let now $X$ be the complement in $\P(V^\dual)$ of the union $H_1 \cup \dots \cup H_m$. The collection of linear forms $v_j$ gives an embedding of $X\hookrightarrow T = \spec(\k[M]) \simeq \rquot{\G_m^{m}}{\G_m}$ with $\G_m = \spec(\k[\x^{\pm1}])$ the algebraic torus over $\k$. The quotient $\rquot{\G_m^{m}}{\G_m}$ here refers to the diagonal action of $\G_m$ on $\G_m^{m}$. In this case, the tropicalization of $X \hookrightarrow T$ coincides with the Bergman fan $\Sigma_\Ma$, see~\cites{AK06, MS15} for more details.

\smallskip
As mentioned above, we have the following proposition.

\begin{prop}\label{prop:productBergman} The product of two generalized Bergman fans is again a generalized Bergman fan.
\end{prop}

More precisely, for two matroids $\Ma$ and $\Ma'$, we have
\[\supp{\Sigma_{\Ma}} \times \supp{\Sigma_{\Ma'}} \simeq \supp{\Sigma_{\Ma\vee\Ma'}}\]
where $\Ma\vee \Ma'$ is any \emph{parallel connection} of $\Ma$ and $\Ma'$. We give the definition below and refer to~\cite{Oxl06}*{Chapter 7} for more details.

A \emph{pointed matroid} is a matroid $\Ma$ on a ground set $E$ with a choice of a distinguished element $\distel = \distel_\Ma$ in $E$. Let now $\Ma$ and $\Ma'$ be two pointed matroids on the ground sets $E$ and $E'$ with distinguished elements $\distel_\Ma \in E$ and $\distel_{\Ma'} \in E'$, respectively. The parallel connection of $\Ma$ and $\Ma'$ denoted by $\Ma \vee \Ma'$ is by definition the pointed matroid on the ground set $E \vee E' = \rquot{E \sqcup E'}{(\distel_\Ma=\distel_{\Ma'})}$, the wedge sum of the two pointed sets $E$ and $E'$, with distinguished element $\distel = \distel_\Ma = \distel_{\Ma'}$, and with the following collection of bases:
\begin{align*}
\bases( \Ma \vee \Ma') =& \Bigl\{ \,B \cup B' \st B \in \bases(\Ma), B' \in \bases(\Ma') \,\textrm{ with }\,\distel_\Ma\in B\, \textrm{ and }\,  \distel_{\Ma'}\in B'\Bigr \} \\
& \cup  \Bigl\{ \,B \cup B' \setminus\{\distel_\Ma\} \st B \in \bases (\Ma), B' \in \bases(\Ma') \, \textrm{ with }\, \distel_\Ma\in B \, \textrm{ and }\, \distel_{\Ma'}\notin B'\Bigr \}\\
& \cup  \Bigl\{ \,B \cup B' \setminus\{\distel_{\Ma'}\}\st B \in \bases(\Ma), B' \in \bases(\Ma')\, \textrm{ with }\, \distel_\Ma \notin B\, \textrm{ and }\,  \distel_{\Ma'}\in B'\Bigr \}.
\end{align*}

We note that the circuits of $\Ma \vee \Ma'$ are given by
\begin{align*}
\crct(\Ma\vee \Ma') = &\ \crct(\Ma) \cup \crct(\Ma') \\
& \ \cup \Bigl\{ (C \cup C')\setminus \{\distel\} \st C \in\crct(M) \textrm{ with } \distel_\Ma \in C \myand C'\in\crct(M')\textrm{ with } \distel_{\Ma'} \in C' \Bigr\}.
\end{align*}

A parallel connection of two matroids $\Ma$ and $\Ma'$ is the wedge sum of the form $(\Ma, e) \vee (\Ma', e')$ for choices elements $e\in \Ma$ and $e'\in \Ma'$ which make pointed matroids out of them.

The analogous operation for graphs consists in gluing two different graphs along distinguished oriented edges as illustrated below. One can check that the relations between the bases, \resp the circuits, of three matroids involved in a parallel connection mimics the relations between the spanning trees, \resp the circuits, of the three graphs.

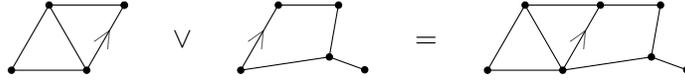
\[ \tikz[vcenter]{\draw (0,0)[dot] --++(60:1)[dot] --+(0:1)[dot] ++(0,0) --++(-60:1)[dot] --(0,0); \draw (0:1)--++(60:1)[dot] node[midway, sloped] {>};}
\quad \vee \quad
\tikz[vcenter]\draw (0,0)[dot]--++(60:1)[dot] node[midway, sloped] {>} --++(0:.8)[dot] --++(80:-.7)[dot] --+(-20:.5)[dot] ++(0,0)--(0,0);
\quad = \quad
\tikz[vcenter]{\draw (0,0)[dot] --++ (60:1)[dot] --+(0:1) ++(0,0) --++(-60:1)[dot] --(0,0); \draw (0:1)--++(60:1)[dot] node[midway, sloped] {>} --++(0:.8)[dot] --++(80:-.7)[dot] --+(-20:.5)[dot] +(0,0) --(0:1);} \]

\begin{proof}[Proof of Proposition~\ref{prop:productBergman}]
Let $\Ma$ and $\Ma'$ be two matroids on ground sets $E=[m]$ and $E'=[m']$, respectively. We show that the support of $\Sigma_{\Ma}\times\Sigma_{\Ma'}$ is isomorphic, by an integral linear isomorphism on the ambient spaces, to $\Sigma_{\Ma\vee \Ma'}$.

Consider the following maps
\[ \R^{E\vee E'} \xrightarrow{\ \phi\times \phi'\ } \R^{E} \times \R^{E'} \xrightarrow{\ \pi \times \pi' \ } \rquot{\R^{E}\!}{\R\ones} \times \rquot{\R^{E'}\!}{\R\ones}. \]
Here, $\pi$, $\pi'$, $\phi$ and $\phi'$ are the natural projections.

We follow the notations introduced in Proposition~\ref{prop:support_bergman} and consider the subsets $\hat \Sigma_{\Ma \vee \Ma'}$, $\hat\Sigma_{\Ma}$ and $\hat\Sigma_{\Ma'}$ of $\R^{E\vee E'}$, $\R^E$ and $\R^{E'}$, respectively.

Since $\crct(\Ma\vee \Ma')$ contains $\crct(\Ma)$ and $\crct(\Ma')$, it is clear that $\phi\times \phi'$ restricts to a map $\hat\phi\times \hat \phi'$ from $\hat\Sigma_{\Ma\vee \Ma'}$ to $\hat\Sigma_{\Ma} \times \hat\Sigma_{\Ma'}$. Let us denote by $\psi\colon \hat\Sigma_{\Ma\vee\Ma'} \to \supp{\Sigma_{\Ma}}\times \supp{\Sigma_{\Ma'}}$ the composition of this map with the projection $\pi\times\pi'$. It is clear that $\psi$ is a linear map and that $\ker(\psi)=\R\ones$. Hence, it remains to prove that $\psi$ is surjective. But the surjectivity can be checked directly by using the description of the circuits in $\Ma \vee \Ma'$ given above and by applying Proposition~\ref{prop:support_bergman}.
\end{proof}

\subsubsection{Proof of the shellability of quasi-projective Bergman fans} We are now in position to give the proof.

\begin{proof}[Proof of Theorem~\ref{thm:bergman_shellable}] Let $\Ma$ be a simple matroid with Bergman fan $\Sigma_\Ma$. We need to show that any quasi-projective fan $\Sigma$ with $\supp{\Sigma} = \supp{\Sigma_\Ma}$ in $\uqproj$ is shellable.

By Theorem~\ref{thm:shellability_factorization}, it will be enough to produce one such fan $\Sigma$ and show it is shellable.

First, we observe that since $\Lambda \in \Bsh_0$, the product $\Lambda^n$ is shellable. This implies that complete fans in $\uqproj$ are all shellable.

It follows that projective fans are shellable. We can therefore assume that $\Ma$ is not a free matroid. We now proceed by induction on the size of $\Ma$.

There exists a proper element of $\Ma$ that is an element $i$ in the ground set $E$ of $M$ such that $\Ma \del \{i\}$ has the same rank as $\Ma$. It is well-known that in this case, the tropical modification of $\Sigma_{\Ma \del \{i\}}$ along the divisor $\Sigma_{\Ma \contr \{i\}} \subset \Sigma_{\Ma \del \{i\}}$ is a quasi-projective fan of support $\supp{\Sigma_\Ma}$, see~\cite{Sha13a}. By induction, both $\Sigma_{\Ma \del \{i\}}$ and $\Sigma_{\Ma \contr \{i\}}$ are shellable in $\uqproj$. Since $\Sigma_\Ma$ is quasi-projective, it is itself shellable in $\uqproj$ and the result follows.
\end{proof}

Example \ref{ex:shellable_not_Bergman} shows that the property of \emph{being a generalized Bergman fan} is not shellable, and at the same time provides an example of a shellable tropical fan which is not a generalized Bergman fan.

\section{Chow rings of tropical fans}\label{sec:chow_ring}

In this section, we study Chow rings of tropical fans. We define a natural cycle class map from the Chow groups to the tropical homology groups of the compactification $\comp \Sigma$, and moreover, study the behavior of Chow groups under tropical modifications and blow-ups. From this description, we deduce the shellability of local irreducibility and div-faithfulness stated in Theorem \ref{thm:examples_shellable}.

Unless otherwise explicitly stated, $\Sigma$ will be a simplicial tropical fan; we will assume unimodularity for the cycle class map and for the Poincaré duality.

\subsection{The additive structure of the Chow ring of a simplicial fan} Consider a rational simplicial fan $\Sigma$ in $N_\R$. We do not assume in this subsection that $\Sigma$ is tropical.

Recall that the Chow ring $A^\bul(\Sigma)$ is defined as the quotient
\[A^\bul(\Sigma) = \rquot{\Z[\x_\rho \mid \rho\in \Sigma_1]}{(I + J)}\]
where
\begin{itemize}
\item $I$ is the ideal generated by the products $\x_{\rho_1}\!\cdots \x_{\rho_k}$, for $k\in \N$, such that $\rho_1, \dots, \rho_k$ do not form a cone in $\Sigma$, and
\item $J$ is the ideal generated by the elements of the form
\[\sum_{\rho\in \Sigma_1} \ell(\e_\rho)\x_\rho\]
for $\ell \in M$ where the primitive vector $\e_\rho$ is the generator of $N\cap\rho$.
\end{itemize}
In order to distinguish to which fan we refer, sometimes we denote by $N_\Sigma$ and $M_\Sigma$ the lattices underlying the definition of $\Sigma$, $M_\Sigma = N_\Sigma^\dual$, and denote by $I_{\Sigma}$ and $J_{\Sigma}$ the ideals used above in the definition of $A^\bul(\Sigma)$.

In the following, if $\sigma$ is a cone of $\Sigma$, we set
\[\x_\sigma := \prod_{\rho\in\Sigma_1\\ \rho\subface\sigma}\x_\rho.\]

Let
\[Z^k(\Sigma) := \bigoplus_{\sigma\in\Sigma_k} \Z\x_\sigma \hookrightarrow \Z[\x_\rho \mid \rho\in \Sigma_1]\]
so that we get a well-defined map $Z^k(\Sigma) \to A^k(\Sigma)$.

\begin{lemma}[Localization Lemma] \label{lem:kernel_Z^k} The map $Z^k(\Sigma) \to A^k(\Sigma)$ is surjective and its kernel is generated by elements of the form
\[ \sum_{\rho\in\Sigma_1 \\ \tau\sim\rho} \ell(\e_\rho)\x_\rho\x_\tau, \]
where $\tau \in \Sigma_{k-1}$ and $\ell$ is an element in $M$ that is zero on $\tau$.
\end{lemma}

\begin{proof}
This fact is a special case of a more general result stated in~\cite{FMSS}*{Theorem 1}, which establishes an isomorphism between the Chow groups of a variety coming with an action of a solvable linear algebraic group with the specific Chow groups defined by cycles and relations which are all stable under the action of the group. A generalization of that theorem can be found in~\cite{Tot14} and subsequent works~\cites{Fra06, Josh01, Pay06}.

\smallskip
We will give an elementary combinatorial proof in Section \ref{sec:kernel_Z^k}.
\end{proof}

\begin{remark} \label{rem:chow_to_mw} In the case the fan $\Sigma$ is unimodular, the above lemma provides an isomorphism $\MW_k(\Sigma) \simeq A^k(\Sigma)^\dual$. Namely, we have a pairing
\[\langle \cdot\,, \cdot \rangle \colon Z^k(\Sigma) \times \MW_k(\Sigma) \to \Z, \qquad \langle \x_\sigma, (\Sigma_k, w)\rangle \to w(\sigma)\]
for any $\sigma\in \Sigma_k$. By the balancing condition for the element $(\Sigma_k, w)$ and Localization Lemma, the pairing vanishes on the kernel of the map $Z^k(\Sigma) \to A^k(\Sigma)$ and we get the desired map $\MW_k(\Sigma) \to A^k(\Sigma)^\dual$. A direct application of Localization Lemma shows that this map is actually an isomorphism, see \cite{AHK}*{Proposition 5.6}. However, $A^k(\Sigma)$ might have torsion (see Examples \ref{ex:F1_vs_N} and \ref{ex:non_saturated_smooth_fan}), in which case $A^k(\Sigma)$ is not isomorphic to $\MW_k(\Sigma)^\dual$. We will discuss torsion-freeness of the Chow ring of tropical fans further in Section~\ref{sec:hodge_isomorphism}.

For $k=d$ and if $\Sigma$ is tropical, this gives the degree map $\deg\colon A^d(\Sigma) \to \Z$ which is obtained by evaluating at the canonical element $[\Sigma] \in \MW_d(\Sigma)$. More precisely, the degree map is obtained by sending the element $\x_\sigma$ of $Z^d(\Sigma)$ to 1, for any $\sigma \in \Sigma_d$.
\end{remark}

\begin{prop}[Künneth formula]
Let $\Sigma$ and $\Sigma'$ be two simplicial fans. Then
\[ A^\bul(\Sigma \times \Sigma') \simeq A^\bul(\Sigma) \otimes A^\bul(\Sigma') \quad\text{and}\quad \MW_\bul(\Sigma \times \Sigma') \simeq \MW_\bul(\Sigma) \otimes \MW_\bul(\Sigma'). \]
Moreover, the first isomorphism is a ring isomorphism.
\end{prop}

In this section, every tensor product is over $\Z$.

\begin{proof}
Notice that for two simplicial fans $\Sigma$ and $\Sigma'$, we have $M_{\Sigma\times\Sigma'}\simeq M_{\Sigma} \times M_{\Sigma'}$. Then one can check that
\begin{gather*}
I_{\Sigma}\otimes \Z[\x_\varrho \mid \varrho \in \Sigma'_1]\,+\,\Z[\x_\varrho \mid \varrho \in \Sigma_1]\otimes I_{\Sigma'} = I_{\Sigma\times\Sigma'},\quad\text{and} \\
J_{\Sigma}\otimes \Z[\x_\varrho \mid \varrho \in \Sigma'_1]\,+\,\Z[\x_\varrho \mid \varrho \in \Sigma_1]\otimes J_{\Sigma'} = J_{\Sigma\times\Sigma'}.
\end{gather*}
The first part of the statement then follows because the tensor product is right-exact.

\smallskip
For Minkowski weights, there is a perfect pairing $Z^k(\Sigma) \times \W_k(\Sigma) \to \Z$ where $\W_k(\Sigma)\simeq \Z^{\Sigma_k}$ denotes the group of weights of dimension $k$ on $\Sigma$. The orthogonal sum of the different degrees leads to a perfect pairing between $Z^\bul(\Sigma)$ and $\W_\bul(\Sigma)$. Denote by $Z_0^k(\Sigma)$ the subgroup of $Z^k(\Sigma)$ generated by elements of the form
\[ \sum_{\sigma \ssupface\tau} \ell(\nvect_{\sigma/\tau})\x_\sigma, \]
where $\tau \in \Sigma_{k-1}$ and $\ell$ is an element in $M_{\Sigma}$ that is zero on $\tau$. From the definition of Minkowski weights we infer that $\MW_k(\Sigma)$ is equal to $Z_0^k(\Sigma)^\perp$. One can check that
\[ Z_0^\bul(\Sigma) \otimes Z^\bul(\Sigma') \,+\, Z^\bul(\Sigma) \otimes Z_0^\bul(\Sigma') = Z_0^\bul(\Sigma\times\Sigma'). \]
Taking the orthogonal of each member in $\W_\bul(\Sigma)\otimes \W_\bul(\Sigma') \simeq \W_\bul(\Sigma \times \Sigma')$, we get the desired isomorphism $\MW_\bul(\Sigma)\otimes\MW_\bul(\Sigma')\simeq\MW_\bul(\Sigma\times\Sigma')$.
\end{proof}

\subsection{The cycle class map from the Chow ring to the tropical homology} \label{subsec:cycle_class_map} Let now $\Sigma$ be a \emph{unimodular} tropical fan. In this section, we prove the following proposition.

\begin{prop}\label{prop:map-chow-homology}
For any integer $k$, there exists a well-defined map $\psi\colon A^k(\Sigma) \to H_{d-k,d-k}(\comp\Sigma)$ defined by sending $\x_\sigma$ to the canonical element $\nu_{\comp\Sigma_\infty^\sigma}$ associated to $\comp\Sigma_\infty^\sigma$.
\end{prop}

For any $\sigma$, the support of $\nu_{\comp\Sigma_\infty^\sigma}$ is the compactification of a tropical fan. Thus, it is a tropical cycle of $\comp\Sigma$. One can thus think of the image of $A^\bul(\Sigma)$ as a part of the homology generated by tropical Chow groups of $\comp\Sigma$, endowed with the intersection product.

\smallskip
Recall from Section~\ref{sec:prel} that for simplicial $\Sigma$, and $\tau \subface \eta$, the face $\comp C^\tau_\eta$ of $\comp \Sigma$ is denoted by $\cube^\tau_\eta$.

\begin{proof}
The fact that $\nu_{\comp\Sigma_\infty^\sigma}$ is a cycle comes from the balancing condition. By Localization Lemma~\ref{lem:kernel_Z^k}, it remains to check that, for any $\tau \in \Sigma_{k-1}$ and any $\ell \in M$ orthogonal to $\tau$, $\sum_{\rho\sim\tau}\ell(\e_\rho)\nu_{\comp\Sigma_\infty^{\tau\vee\rho}}$ is a boundary. This then shows the natural map $Z^k \to H_{d-k,d-k}(\comp\Sigma)$ passes to the quotient and defines the map stated in the proposition.

Since $\ell$ is orthogonal to $\tau$, it induces a linear map on $N^\tau$ which we also denote by $\ell$. To each facet $\eta \supface \tau$, we associate the contraction of $\nu^\tau_\eta$ by $\ell$ denoted $\iota_\ell(\nu^\tau_\eta) \in \SF_{d-k}(\cube_\eta^\tau)$: if $\rho_1, \dots, \rho_{d-k+1}$ are the rays generating $\cube_\eta^\tau$ in $N_{\infty,\R}^\tau$ such that $\nu^\tau_\eta = \e_{\rho_1} \wedge \dots \wedge \e_{\rho_{d-k+1}}$, then
\[ \iota_\ell(\nu^\tau_\eta) = \sum_{i = 1}^{d-k+1} (-1)^{i-1} \ell(\e_{\rho_i})\e_{\rho_1}\wedge\dots\wedge\hat{\e_{\rho_i}}\wedge\dots\wedge\e_{\rho_{d-k+1}}. \]
Set $\iota_\ell(\nu_{\comp\Sigma_\infty^\tau}) := \sum_{\eta\in\Sigma_d \\ \eta\supface\tau} \iota_\ell(\nu^\tau_\eta)$. A direct computation shows that
\[ \partial\bigl(-\iota_\ell(\nu_{\comp\Sigma_\infty^\tau})\bigr) = \sum_{\rho \in \Sigma_1 \\ \rho\sim\tau} \ell(\e_\rho) \nu_{\comp\Sigma_\infty^{\tau\vee\rho}}, \]
which proves the claim, and concludes the proof.
\end{proof}

There is another natural cycle class map. It is given by the following composition.
\[ \phi\colon A^k(\Sigma) \to A^{d-k}(\Sigma)^\dual \simto \MW_{d-k}(\Sigma) \simto H_{d-k,d-k}(\comp\Sigma). \]
The first map sends $\alpha$ to $\bigl(\beta \mapsto \deg(\alpha\cdot\beta)\bigr)$ and the last map sends a Minkowski weight $(\Sigma_{(d-k)},w)$ to $\bigl(w(\sigma)\nu_\sigma\bigr)_{\sigma\in\Sigma_{d-k}} \in C_{d-k,d-k}(\comp\Sigma)$. This last map is studied in detail in Section \ref{sec:hodge_isomorphism} where we prove that it is an isomorphism.

We claim the above two cycle class maps coincide. Let us briefly justify this fact.

Let $\sigma$ be a face of dimension $k$ in $\Sigma$. There is a natural ring morphism $\i^*\colon A^\bul(\Sigma) \to A^\bul(\Sigma_\infty^\sigma)$ described for instance in \cites{AHK, AP-tht}. This map is entirely determined by the property that, for a ray $\rho\not\subface\sigma$, it maps $\x_\rho$ to 0 if $\rho$ is not comparable with $\sigma$, and to $\x_{\rho^\sigma}$ if $\rho\sim\sigma$ where $\rho^\sigma$ is the corresponding ray of $\Sigma_\infty^\sigma$. Using the arguments in Section \ref{subsec:Ap_to_Hpp}, one can check that the map $A^p(\Sigma) \to A^p(\Sigma_\infty^\sigma)$ commutes with the restriction map $\i^*\colon H^{p,p}(\comp\Sigma) \to H^{p,p}(\comp\Sigma_\infty^\sigma)$. (This is the reason the map between the Chow rings is denoted $\i^*$.) With this in mind, we get the following commutative diagram.
\[ \begin{tikzcd}
A^0(\Sigma_\infty^\sigma) \rar\dar& A^{d-k}(\Sigma_\infty^\sigma)^\dual \dar\rar{\sim}& \MW_{d-k}(\Sigma_\infty^\sigma) \rar{\sim}\dar& H_{d-k, d-k}(\comp\Sigma_\infty^\sigma) \dar{\i_*} \\
A^k(\Sigma)            \rar& A^{d-k}(\Sigma)^\dual            \rar{\sim}& \MW_{d-k}(\Sigma)            \rar{\sim}& H_{d-k,d-k}(\comp\Sigma)
\end{tikzcd} \]
The second and last vertical maps are the dual of the maps $\i^*$ described above. The first vertical map sends 1 to $\x_\sigma$. In the first row, 1 is mapped to the Minkowski weight $[\Sigma_\infty^\sigma] \in \MW_{d-k}(\Sigma_\infty^\sigma)$ and then to $\nu_{\comp\Sigma_\infty^\sigma}$. Hence, the last vertical map sends $\nu_{\comp\Sigma_\infty^\sigma}=\psi(\x_\sigma)$ to $\phi(\x_\sigma)$.

\subsection{Chow ring and tropical modifications}\label{sec:chow_tropmod} In this section, we prove that the Chow ring and the cycle class map behave well with respect to tropical modifications.

Let $\Sigma$ be a simplicial tropical fan of dimension $d$ and let $\Delta = \div(f)$ be the principal divisor associated to a conewise integral linear function on $\Sigma$. Assume that $\Delta$ is effective reduced, and let $\~\Sigma$ be the fan obtained by the open tropical modification of $\Sigma$ with respect to $\Sigma$.

\begin{prop} \label{prop:chow_ring_tropical_modification} Notations as above, we have a natural surjective ring morphism
\[ A^\bul(\Sigma) \twoheadrightarrow A^\bul(\~\Sigma). \]
\end{prop}

\begin{proof}
Following the notations of Section \ref{subsec:tropical_modification}, for any ray $\varrho \in \Sigma_1$, we denotes by $\basetm\varrho$ the corresponding ray in $\~\Sigma_1$. We assume for now that $\Delta$ is nontrivial and denote by $\rho$ the only ray of $\~\Sigma$ that does not come from $\Sigma$, that is, $\rho = \uptm{\conezero\,}$ where $\conezero$ refers to the zero cone $\conezero_\Sigma$ of $\Sigma$.

Consider the map $\psi\colon \Z[\x_\varrho \mid \varrho\in\Sigma_1] \to \Z[\x_{\varrho} \mid \varrho\in\~\Sigma_1]$ defined by sending $\x_\varrho$ to $\x_{\varrho_0}$. Let $I_\Sigma$, $J_\Sigma$, $I_{\~\Sigma}$ and $J_{\~\Sigma}$ be the ideals appearing in the definition of the Chow rings of $\Sigma$ and $\~\Sigma$. We clearly have $\psi(I_\Sigma) \subseteq I_{\~\Sigma}$. Moreover, if $\ell$ is a linear form on the ambient space of $\Sigma$ and if $\prtm$ denotes the projection associated to the tropical modification, since $\prtm^*\ell$ is zero on the special ray $\rho$, we get that
\[ \psi\Bigl(\sum_{\varrho \in \Sigma_1} \ell(\e_\varrho)\x_\varrho\Bigr) = \sum_{\varrho \in \Sigma_1} \prtm^*\ell(\e_{\basetm\varrho}) \x_{\basetm\varrho} = \sum_{\varrho \in \~\Sigma_1} \prtm^*\ell(\e_{\varrho}) \x_{\varrho}. \]
This implies that $\psi(J_\Sigma) \subseteq J_{\~\Sigma}$. All together, we thus get a well-defined map $\psi$ from $A^\bul(\Sigma)$ to $A^\bul(\~\Sigma)$.

\smallskip
To prove the surjectivity, by the definition of $\psi$, we just need to find a preimage for $\x_\rho \in A^\bul(\~\Sigma)$. For this, take any integral linear form $\ell$ on $\~N$ that takes value one on the primitive vector $\e_\rho$ of the special ray $\rho$. In the Chow ring $A^\bul(\~\Sigma)$, using $\ell$, we can express $\x_\rho$ as a linear combination of the $\x_{\varrho}$ for $\varrho \in \~\Sigma_1 - \rho$. This linear combination is in the image of $\psi$, and the surjectivity of $\psi$ follows. In the case $\Delta=\div(f)$ is trivial, the proof is similar.
\end{proof}

For the sake of completeness, we state the following proposition and refer to Sections~\ref{sec:hodge_isomorphism}~and~\ref{sec:homology_tropical_modification} for the definition of the relevant maps used in the statement.
\begin{prop}
Let $\Sigma$ be a unimodular tropical fan. The morphism $A^\bul(\Sigma) \to A^\bul(\~\Sigma)$ commutes via the Hodge isomorphism theorem with the map $H^{\bul,\bul}(\Sigma) \to H^{\bul,\bul}(\~\Sigma)$ studied in Section~\ref{sec:homology_tropical_modification}. Hence, it also commutes with cycle class maps and we get the following diagram.
\[ \begin{tikzcd}
  A^k(\~\Sigma)  \rar                        &  H_{d-k,d-k}(\comp{{\~\Sigma}}) \dar \\
  A^k(\Sigma)    \uar[twoheadrightarrow]\rar &  H_{d-k,d-k}(\comp\Sigma)
\end{tikzcd} \]
\end{prop}

\begin{proof}
The proof is straightforward once the different maps have been defined.
\end{proof}

\subsection{Shellability of the principality} In this section we prove the following theorem.

\begin{thm}[Shellability of the principality for locally irreducible tropical fans]\label{thm:shellability_principal} The property for a tropical fan to be both principal and locally irreducible is shellable.
\end{thm}

Local irreducibility assumption is needed here as shown by Example \ref{ex:div-faithful_not_irreducible}. To prove the theorem, we need to check the different cases of Lemma \ref{lem:shellability_meta_lemma}. The fact that elements of $\Bsh_0$ are locally irreducible and principal is trivial.

\subsubsection{Closeness under products} Consider two locally irreducible principal tropical fans $\Sigma, \Sigma'$ of dimension $d$ and $d'$, respectively. By Theorem \ref{thm:examples_shellable} \ref{thm:examples_shellable:irreducible}, local irreducibility is shellable. Hence $\Sigma\times\Sigma'$ is locally irreducible. Let us prove that $\Sigma\times\Sigma'$ is principal at $\conezero$. By Künneth formula for Minkowski weights, we get
\[\MW_{d+d'-1}(\Sigma\times \Sigma') = \MW_d(\Sigma) \otimes \MW_{d'-1}(\Sigma')\ \oplus\ \MW_{d-1}(\Sigma) \otimes \MW_{d'}(\Sigma'). \]
Using the above decomposition, and arguing by symmetry, we only need to verify that divisors of the form $D \times C'$ for $D \in \MW_{d-1}(\Sigma)$ and $C'\in \MW_{d'}(\Sigma')$ are principal.

Since $\Sigma'$ is locally irreducible, we can suppose that $C'$ is equal to the element $[\Sigma']$ in $\MW_{d'}(\Sigma')$ given by $\Sigma'$. Since $\Sigma$ is principal, we have $D = \div(f)$ for a conewise integral linear function $f$ on $\Sigma$. Denote by $\pi\colon \Sigma\times\Sigma' \to \Sigma$ the natural projection. We have $\div(\pi^*f) = \div(f)\times[\Sigma'] = D\times C'$, from which the result follows.

\subsubsection{Closeness under tropical modifications} Let $\Sigma$ be a locally irreducible principal tropical fan, and let $\~\Sigma$ be a tropical modification of $\Sigma$. Once again $\~\Sigma$ is locally irreducible. By Lemma \ref{lem:shellability_meta_lemma}, we can assume without loss of generality that $\Sigma$ is unimodular.

By Proposition~\ref{prop:chow_ring_tropical_modification}, the applications $A^1(\Sigma) \to A^1(\~\Sigma)$ and $A^{d-1}(\Sigma) \to A^{d-1}(\~\Sigma)$ are both surjective. It follows that the dual map $A^{d-1}(\~\Sigma)^* \to A^{d-1}(\Sigma)^*$ is injective, and from the explicit description of the map $\cycl^1$ given in Proposition~\ref{prop:map_cl}, we get the following commutative diagram
\[ \begin{tikzcd}
  A^1(\~\Sigma)  \rar                        &  A^{d-1}(\~\Sigma)^* \dar[hook]\rar{\sim} & \MW_{d-1}({\~\Sigma}) \dar[hook] \\
  A^1(\Sigma)    \uar[twoheadrightarrow]\rar &  A^{d-1}(\Sigma)^*  \rar{\sim}            & \MW_{d-1}(\Sigma)
\end{tikzcd}
\]
We now show that $\~\Sigma$ is principal at $\conezero$. By Theorem~\ref{thm:char_div-faithful} this is equivalent to the surjectivity of the map $A^1(\~\Sigma)\to \MW_{d-1}(\~\Sigma)$ (saturation is not needed for this part of the theorem). Moreover, since $\Sigma$ is principal, $A^1(\Sigma) \to \MW_{d-1}(\Sigma)$ is surjective. Using the above diagram, we obtain the surjectivity of $A^1(\~\Sigma)\to \MW_{d-1}(\~\Sigma)$. This in turn implies that $\~\Sigma$ is principal at $\conezero$, as desired.

\subsubsection{Closeness under blow-ups and blow-downs}

Once again, closeness of local irreducibility has already been checked. Consider a tropical fan $\Sigma$ of dimension $d$ and let $\sigma$ be a cone of $\Sigma$ and $\rho$ a ray inside $\sigma$. Set $\Sigma' = \Sigma_{(\rho)}$. Assume that $\Sigma$ and $\Sigma'$ are principal at any nontrivial face. We need to prove $\Sigma$ is principal at $\conezero_{\Sigma}$ if and only if $\Sigma'$ is principal at $\conezero_{\Sigma'}$.

Assume first that $\Sigma$ is principal. Let $D$ be a divisor in $\Sigma'$. Let $D^\rho$ be the induced divisor on $\Sigma'^\rho$. We infer the existence of a conewise integral linear function $f^\rho$ on $\Sigma'^\rho$ such that $\div(f^\rho) = D^\rho$. Via the projection map $N_\R \to N^\rho_\R$, we view $f^\rho$ as a conewise integral linear function on the fan $\Sigma'' \subset \Sigma'$ consisting of the cones $\tau$ comparable with $\rho$, and extend it to a conewise integral linear function on full $\Sigma'$. By an abuse of the notation, we denote this function by $f^\rho$. The divisor $D - \div(f^\rho)$ on $\Sigma'$ does not have any of the faces $\tau \supface \rho$ with $\tau \in \Sigma'_{d-1}$ in its support. This means it can be viewed as a divisor in $\Sigma$. Since $\Sigma$ is principal, we can find a conewise integral linear function $f$ on $\Sigma$ with $\div(f) = D - \div(f^\rho)$. We infer that $D$ is the divisor of the conewise integral linear function $f + f^\rho$ on $\Sigma'$, which shows that $D$ is principal, as desired.

Assume now that $\Sigma'$ is principal. Let $D$ be divisor in $\Sigma$. Viewing $D$ as a divisor in $\Sigma'$, we find a conewise integral linear function $f$ on $\Sigma'$ such that $\div(f) = D$. Take a facet $\eta\supface\sigma$ in $\Sigma$. Denote by $\Sigma'\rest\eta$ the fan of support $\eta$ induced by $\Sigma'$. If $f$ is not linear on $\eta$, then one can find a face $\tau$ of dimension $d-1$ in $\Sigma'\rest\eta \setminus \Sigma$ such that $\ord_f(\tau)\neq0$. This is impossible since $D$ is supported on faces of $\Sigma$. Hence $f$ is linear on all the facets of $\Sigma$ containing $\sigma$. Therefore $f$ is conewise integral linear on $\Sigma$ which proves that $D$ is principal in $\Sigma$.

\subsubsection{Proof of Theorem~\ref{thm:shellability_principal}} At this point, we have verified all the cases of Lemma \ref{lem:shellability_meta_lemma} and this concludes the proof of Theorem~\ref{thm:shellability_principal}. \qed

\subsection{Div-faithfulness and stability of the Chow ring under tropical modifications} The following theorem implies the invariance of the Chow ring under tropical modifications under the assumption that the underlying tropical fan $\Sigma$ is div-faithful.

\smallskip
We follow the preceding notations and denote by $\~\Sigma$ the tropical modification of $\Sigma$ along $\Delta =\div(f)$.

\begin{thm}[Stability of the Chow ring under tropical modifications] \label{thm:invariant_chow_tropical_modification} Notations as above, let $\~\Sigma$ be the tropical modification of the tropical fan $\Sigma$ along $\Delta$. Assume furthermore that $\Sigma$ is div-faithful. Then we get an isomorphism
\[ A^\bul(\~\Sigma) \simeq A^\bul(\Sigma) \]
between the Chow rings. Moreover, this isomorphism is compatible with the Poincaré duality for Chow rings.
\end{thm}

\begin{proof}
If $\Delta$ is empty, since $\Sigma$ is div-faithful, $f$ is linear and $\~\Sigma \simeq \Sigma$ from which the result follows. Hence, in what follows, we assume $\Delta\neq\emptyset$, and we denote by $\rho$ the special ray of $\~\Sigma$. We already know there exists a surjective map $\psi\colon A^\bul(\Sigma) \twoheadrightarrow A^\bul(\~\Sigma)$ given in Proposition \ref{prop:chow_ring_tropical_modification}. We now construct a surjective map $\phi\colon A^\bul(\~\Sigma) \twoheadrightarrow A^\bul(\Sigma)$ and show that $\phi\circ\psi=\id$, from which we get the result.

\smallskip
Consider first the surjective map $\phi \colon \Z[\x_\varrho \mid \varrho\in\~\Sigma_1] \to \Z[\x_{\varrho} \mid \varrho\in \Sigma_1]$, which appeared implicitly in the proof of Proposition~\ref{prop:chow_ring_tropical_modification}. This is defined on the level of generators as follows. Take the linear form $\ell$ on $\~N$ that takes value one on the primitive vector $\e_\rho$ of the special ray $\rho$ and that is zero on $N\times\{0\}$. For each ray of $\~\Sigma$ of the form $\basetm\varrho$ for $\varrho \in \Sigma_1$, set $\phi(\x_{\basetm \varrho}) := \x_\varrho$, and define
\[\phi(\x_\rho) = - \sum_{\varrho \in\Sigma} \ell(\e_{\basetm\varrho})\x_\varrho.\]
Notice that from the equality $\e_{\basetm\varrho}=(\e_\varrho,f(\e_\varrho))$, we get $\ell(\e_{\basetm\varrho})=f(\e_\varrho)$.

We claim that $\phi$ passes to the quotient. In this case, it is clear that $\phi$ induces an inverse to $\psi$, and the isomorphism between $A^\bul(\Sigma)$ and $A^\bul(\~\Sigma)$ stated in the theorem follows.

Denote by $I_\Sigma$, $J_\Sigma$, $I_{\~\Sigma}$ and $J_{\~\Sigma}$ the ideals appearing in the definition of the Chow rings of $\Sigma$ and $\~\Sigma$, as previously.

We show that we have $\phi(J_{\~\Sigma}) \subseteq J_\Sigma$ and $\phi(I_{\~\Sigma}) \subseteq I_\Sigma +J_\Sigma$, which allow to pass to the quotient and get the above map on the level of Chow rings.

\smallskip
To show that $\phi(J_{\~\Sigma}) \subseteq J_\Sigma$, consider an integral linear form $l$ on $\~N$. Let $l' := l - l(\e_\rho)\ell$. We have $l'(\e_\rho)=0$ and so $l'$ gives an integral linear form on $N \simeq \rquot{\~N}{\Z \e_\rho}$. We have
\[ \phi\Bigl(\sum_{\varrho\in\~\Sigma_1} l(\e_\varrho)\x_\varrho\Bigr) = l(\e_\rho)\phi(\x_{\rho}) + \sum_{\varrho\in \Sigma_1} l(\e_{\basetm \varrho})\phi(\x_{\basetm\varrho})=\sum_{\varrho\in \Sigma_1} l'(\e_\varrho)\x_\varrho. \]
This shows that $\phi(J_{\~\Sigma}) \subseteq J_\Sigma$.

\smallskip
We now consider the image of $I_{\~\Sigma}$. Consider a collection of distinct rays $\rho_1, \dots, \rho_k$ of $\~\Sigma$, $k\in \N$, and suppose they are not comparable so that we have $\x_{\rho_1} \dots \x_{\rho_k} \in I_{\~\Sigma}$. Two cases can happen.
\begin{enumerate}
\item \label{enum:invariant_chow:different_rho} Either $\rho_1, \dots, \rho_k$ are different from $\rho$.
\item \label{enum:invariant_chow:contains_rho} Or one of the rays, say $\rho_1$, is equal to $\rho$.
\end{enumerate}
Consider first the case \ref{enum:invariant_chow:different_rho}. Then we have $\rho_j = \basetm{(\varrho_j)}$ for $j=1, \dots, k$ and rays $\varrho_j$ in $\Sigma$. Moreover, the rays $\varrho_1, \dots, \varrho_k$ do not form a cone in $\Sigma$ and we get $\phi(\x_{\rho_1} \dots \x_{\rho_k}) = \x_{\varrho_1} \dots \x_{\varrho_k} \in I_\Sigma$, as desired.

\smallskip
We now study the case \ref{enum:invariant_chow:contains_rho}. Let as in the previous case, $\rho_j = \basetm{(\varrho_j)}$ for $j=2, \dots, k$ and rays $\varrho_j$ in $\Sigma$. At this point, two cases can happen. Either, these rays do not form a cone in $\Sigma$ in which case we get
\[\phi(\x_{\rho_1}\x_{\rho_2} \dots \x_{\rho_k}) = \phi(\x_{\rho}) \x_{\varrho_2} \dots \x_{\varrho_k} \in I_\Sigma.
\]

Or, $\varrho_2, \dots, \varrho_k$ form a cone $\tau$ in $\Sigma$. Note that $\tau \not\in\Delta$: otherwise $\basetm\tau$ and $\rho$ would be comparable which would contradict our assumption. This implies that the divisor induced by $\div(f)$ on $\Sigma^\tau$ is trivial. In particular, taking an integral linear function $l$ on $N_\R$ which has restriction to $\tau$ equal to $f$, the difference $f^\tau:=f-l$ induces a conewise integral linear function $\pi_*(f^\tau)$ on $\Sigma^\tau$, where $\pi\colon N \to N^\tau$ is the projection, and $\div(\pi_*(f^\tau)) = 0$. Since $\Sigma$ is div-faithful at $\tau$, this implies that $\pi_*(f^\tau)$ coincides with an integral linear function $l^\tau$ on $N^\tau_\R$. All together, this means $f$ coincides with the restriction of the linear function $l_\tau:=l+\pi^*(l^\tau)$ on all faces $\sigma$ of $\Sigma$ with $\sigma \supface \tau$. Set $\~f = f-l_\tau$. In particular, $\~f$ is zero on every ray comparable with $\tau$.

We now observe that
\[ \phi(\x_\rho) = -\sum_{\varrho\in\Sigma_1} f(\e_{\varrho}) \x_\varrho = -\sum_{\varrho\in\Sigma_1} \~f(\e_{\varrho}) \x_\varrho - \sum_{\varrho\in\Sigma_1} l_\tau(\e_\varrho) \x_\varrho, \]
and so, using the notation $\x_\tau = \x_{\varrho_2} \cdots \x_{\varrho_k}$,
\[ \phi(\x_{\rho_1}\x_{\rho_2} \dots \x_{\rho_k}) = \phi(\x_{\rho})\x_{\tau}
  = -\sum_{\varrho\in\Sigma_1} \~f(\e_\varrho)\x_\varrho\x_{\tau} - \Bigl(\sum_{\varrho\in\Sigma_1} l_\tau(\e_\varrho) \x_\varrho\Bigr) \x_{\tau}. \]
Since $\~f(\e_\varrho)$ is trivial if $\varrho$ is comparable with $\tau$, the first sum is in $I_\Sigma$. Since $l_\tau \in M$, the second term is in $J_\Sigma$.

In any case, we conclude that $\phi(\x_{\rho_1}\x_{\rho_2} \cdots \x_{\rho_k}) \in I_\Sigma + J_\Sigma$, as desired. This finishes the proof of the theorem.
\end{proof}

\subsection{Shellability of div-faithfulness}
In this section we prove the following.

\begin{thm} \label{thm:shellability_div-faithful}
The property for a tropical fan to be div-faithful is shellable.
\end{thm}

Once again, we study the different cases of Lemma \ref{lem:shellability_meta_lemma}.

\subsubsection{Closeness under products} Consider two tropical fans $\Sigma$ and $\Sigma'$ that are div-faithful. It will be enough to show that $\Sigma \times \Sigma'$ is div-faithful at $\conezero$. We can assume without loss of generality that $\Sigma$ and $\Sigma'$ are unimodular (by Lemma \ref{lem:shellability_meta_lemma}) and saturated at $\conezero$ (see Remark \ref{rem:independent_exterior_lattice}).

We use the decomposition
\[\MW_{d+d'-1}(\Sigma\times \Sigma') = \MW_{d-1}(\Sigma) \otimes \MW_{d'}(\Sigma')\ \oplus\ \MW_d(\Sigma) \otimes \MW_{d'-1}(\Sigma'). \]
We have injections $\Z=A^0(\Sigma)\hookrightarrow \MW_d(\Sigma)$ and $A^0(\Sigma') \hookrightarrow \MW_{d'}(\Sigma')$. We have moreover a decomposition
\[A^1(\Sigma \times \Sigma') \simeq A^1(\Sigma)\otimes A^0(\Sigma')\ \oplus\ A^0(\Sigma)\otimes A^1(\Sigma'). \]

The injectivity of the map $A^1(\Sigma \times \Sigma')\to\MW_{d+d'-1}(\Sigma\times \Sigma')$ now follows from the injectivity of the corresponding maps for $\Sigma$ and $\Sigma'$, given by Theorem~\ref{thm:char_div-faithful}, which shows that $\Sigma\times \Sigma'$ is div-faithful at $\conezero$, as desired.

\subsubsection{Closeness under tropical modifications} Let $\Sigma$ be a tropical fan that is div-faithful. Let $\~\Sigma$ be a tropical modification of $\Sigma$ and let us prove that $\~\Sigma$ is div-faithful at $\conezero$. Once again we assume without loss of generality that $\Sigma$ is unimodular and saturated at $\conezero$.

By Theorem~\ref{thm:invariant_chow_tropical_modification}, we get isomorphisms $A^1(\Sigma) \simeq A^1(\~\Sigma)$ and $A^{d-1}(\Sigma) \simeq A^{d-1}(\~\Sigma)$. The injectivity of the map $A^1(\~\Sigma)\to\MW_{d-1}(\~\Sigma)$ then follows from the injectivity of $A^1(\Sigma)\to\MW_{d-1}(\Sigma)$. This implies that $\~\Sigma$ is div-faithful at $\conezero$.

\subsubsection{Closeness under blow-ups and blow-downs} Consider a fan $\Sigma$, and let $\sigma$ be a cone of $\Sigma$ and $\rho$ a ray inside $\sigma$. Set $\Sigma' = \Sigma_{(\rho)}$. Assume that $\Sigma$ and $\Sigma'$ are div-faithful at any nontrivial faces. Let us prove that $\Sigma$ is div-faithful at $\conezero$ if and only if $\Sigma'$ is div-faithful at $\conezero$.

First assume that $\Sigma'$ is div-faithful. Let $f$ be a conewise integral linear function on $\Sigma$ such that $\div(f)=0$. Then $f$ is also a conewise linear function on $\Sigma'$, and since its divisor is trivial, it is linear. This proves the direction $\Leftarrow$.

In the other direction, assume $f$ is a conewise integral linear function on $\Sigma'$ whose divisor is trivial. Then the induced conewise integral linear function $f^\rho$ on $\Sigma^\rho$ also verifies $\div(f^\rho)=0$. By assumption, $\Sigma^\rho$ is div-faithful. Thus, $f^\rho$ is integral linear, and we deduce that $f$ is integral linear on each face $\eta\supface\sigma$ of $\Sigma$. Thus, $f$ is conewise integral linear on $\Sigma$. We infer that $f$ is globally linear which concludes this part of the proof.

\subsubsection{Proof of Theorem~\ref{thm:shellability_div-faithful}} We have studied all the cases of Lemma \ref{lem:shellability_meta_lemma}. This concludes the proof of Theorem~\ref{thm:shellability_div-faithful}. \qed

\subsection{Shellability of the Poincaré duality for the Chow ring of div-faithful unimodular fans} Let $\Sigma$ be a unimodular tropical fan of dimension $d$. In this section, we denote by $\PD_\chow(\Sigma)$ the Poincaré duality for the Chow ring of $\Sigma$. This is the statement that the map
\[A^\bul(\Sigma) \times A^{d-\bul}(\Sigma) \to \Z, \qquad (\alpha, \beta) \to \deg(\alpha\cdot \beta),\]
is a perfect pairing. Recall that the degree map $\deg\colon A^d(\Sigma) \to \Z$ is given by $[\Sigma]\in\MW_d(\Sigma)\simeq A^d(\Sigma)^\dual$. In particular, it maps any element $\x_\sigma$, $\sigma$ a facet of $\Sigma$, to 1.

In case of Poincaré duality for the Chow ring, $A^d(\Sigma)\simeq\Z$ has a natural generator of degree $1$ that we denote $\omega_\Sigma$. We have $\omega_\Sigma=\x_\sigma$ for any facet $\sigma\in\Sigma_d$, this last element being independent of the choice of the facet.

\begin{thm}\label{thm:chow_shellable}
The set of div-faithful fans whose Chow rings verify the Poincaré duality is stable by shellability in $\unim$.
\end{thm}

As a corollary, we get the following theorem.

\begin{thm}\label{thm:chow-Poincare_shellable}
The property $\Ploc{\PD_\chow}$ is shellable in $\unim$.
\end{thm}
\begin{proof}
If $\Sigma$ verifies $\PD_\chow$, then it is div-faithful at $\conezero$. Indeed, in this case, we have $A^1(\Sigma^\sigma) \simeq \MW_{d-1}(\Sigma^\sigma)$. Hence $\Ploc{\PD_\chow}$ implies div-faithfulness. The theorem now follows from Theorem \ref{thm:chow_shellable} via Lemma \ref{lem:shellability_meta_lemma}.
\end{proof}

\begin{remark}
By Remark~\ref{rem:torsion_duality_chow} in the next section, the Poincaré duality for Chow ring is implied by the cohomological Poincaré duality for $\comp\Sigma$. However, the two notions could be a priori different since the Poincaré duality for Chow groups does not imply the vanishing of the cohomology groups $H^{p,q}$ of canonical compactifications for bidegrees $p>q$. This is for instance the case for Example \ref{ex:cube} with rational coefficients. The authors are not aware of any such example for integral coefficients.

Example \ref{ex:PD_chow_not_smooth} describes a fan $\Sigma$ which verifies $\PD_\chow$ and such that $\comp\Sigma$ verifies the Poincaré duality, but $\Sigma$ is not smooth and does not verify $\Ploc{\PD_\chow}$. This shows that $\PD_\chow$ is not a local property.
\end{remark}

The rest of this section is devoted to the proof of this theorem. This time, we cannot use Lemma \ref{lem:shellability_meta_lemma} since $\PD_\chow$ is not local as we discussed in the previous remark. We have already proved that div-faithfulness is shellable. Hence in what follows we only focus on the shellability of $\PD_\chow$.

\subsubsection{Closeness under products} Let $\Sigma$ and $\Sigma'$ be two unimodular tropical fans of dimension $d$ and $d'$, respectively. We know that we have a ring isomorphism
\[ A^\bul(\Sigma \times \Sigma') \simeq A^\bul(\Sigma) \otimes A^\bul(\Sigma'). \]
Moreover, the dual isomorphism sends $\deg_{\Sigma\times\Sigma'}$ to $\deg_\Sigma\otimes\deg_{\Sigma'}$. We infer that $\PD_\chow(\Sigma \times \Sigma')$ holds provided that $\PD_\chow(\Sigma)$ and $\PD_\chow(\Sigma')$ both hold.

\subsubsection{Closeness under tropical modifications} Let $\Sigma$ be a div-faithful tropical fan and let $\~\Sigma$ be a tropical modification of $\Sigma$. Applying Theorem~\ref{thm:invariant_chow_tropical_modification}, we obtain the Chow stability of tropical modification, namely that the natural ring morphism $A^\bul(\Sigma) \to A^\bul(\~\Sigma)$ is an isomorphism. By compatibility of the degree maps, we conclude that $\PD_\chow(\Sigma)$ implies $\PD_\chow(\~\Sigma)$.

\subsubsection{Closeness under blow-ups and blow-downs} \label{sec:Keel} Consider a unimodular fan $\Sigma$ of dimension $d$ and let $\sigma$ be a cone in $\Sigma$. Let $\Sigma' = \Bl\Sigma\sigma$ be the unimodular star subdivision of $\Sigma$ obtained by star subdividing the cone $\sigma$ and denote by $\rho$ the new ray in $\Sigma'$, that is, $\rho=\R_+(\e_1+\dots+\e_k)$ with $\e_1, \dots, \e_k$ the primitive vectors of the rays of $\sigma$. Recall that for any $\sigma\in\Sigma$, we have a surjective map $\i^*_{\conezero\subface\sigma}\colon A^\bul(\Sigma) \to A^\bul(\Sigma^\sigma)$ described in Section \ref{subsec:cycle_class_map}. We state the following useful result.

{ \newcommand{\J}{\mathfrak{J}}
\begin{thm}[Keel's lemma] \label{thm:keel} Let $\J$ be the kernel of the surjective map $\i^*_{\conezero\subface\sigma}\colon A^\bul(\Sigma)\to A^\bul(\Sigma^\sigma)$ and let
\[ P(T):=\prod_{\varrho\subface \sigma \\ \dims\varrho=1}(\x_\varrho+T). \]
We have
\[ A^\bul(\Sigma')\simeq \rquot{A^\bul(\Sigma)[T]}{(\J T+P(T))}. \]
The isomorphism is given by the map
\[\chi\colon \rquot{A^\bul(\Sigma)[T]}{(\J T+P(T))} \simto A^\bul(\Sigma')\]
which sends $T$ to $-\x_\rho$ and which verifies
\[ \forall \varrho \in \Sigma_1, \qquad
\chi(\x_\varrho) = \begin{cases}
\x_\varrho+\x_\rho & \text{if $\varrho \subface \sigma$,}\\
\x_\varrho & \text{otherwise.}
\end{cases}
\]
In particular this gives a vector space decomposition of $A^\bul(\Sigma')$ as
\begin{equation}\label{eq:keel}
A^\bul(\Sigma')\simeq A^\bul(\Sigma)\oplus A^{\bul-1}(\Sigma^\sigma)T \oplus \dots \oplus A^{\bul-\dims{\sigma}+1}(\Sigma^\sigma)T^{\dims{\sigma}-1}.
\end{equation}
\end{thm}

\begin{proof} This follows from~\cite{Kee92}*{Theorem 1 in the appendix} for the map of toric varieties $\P_{\Sigma'} \to \P_{\Sigma}$. Here $P(T)$ is the polynomial in $A^\bul(\P_\Sigma)$ whose restriction in $A^\bul(\P_{\Sigma^\sigma})$ is the Chern polynomial of the normal bundle for the inclusion of toric varieties $\P_{\Sigma^\sigma} \hookrightarrow \P_{\Sigma}$.

For an elementary proof, one can follow the proof of an analogous statement in \cite{AHK}.
\end{proof}

Assume now $\PD_\chow(\Sigma^\sigma)$. We have to show the equivalence of $\PD_\chow(\Sigma)$ and $\PD_\chow(\Sigma')$. We follow the proof of an analogue result in \cite{AHK}.

By Keel's lemma, we get $A^d(\Sigma) \simeq A^d(\Sigma')$, and under this isomorphism, the canonical element $\omega_\Sigma$ gets identified with the canonical element $\omega_{\Sigma'}$.

Let $\alpha$ be a lifting of the canonical element $\omega_{\Sigma^\sigma}$ in $A^{d-\dims\sigma}(\Sigma)$, for the restriction map $\i^*_\sigma\colon A^\bul(\Sigma) \to A^\bul(\Sigma^{\sigma})$, so that we have $\omega_{\Sigma^\sigma} =\i^*_\sigma(\alpha)$. We can obtain such a lifting by fixing a top dimensional cone in $\Sigma^{\sigma}$, \ie, a cone of the form $\eta^\sigma$ for $\eta\in \Sigma_d$ and $\eta \supface \sigma$. Denoting by $\tau$ the face of $\eta$ generated by all the rays of $\eta$ not included in $\sigma$, we get
\[\omega_{\Sigma^\sigma} = \x_{\eta^\sigma} = \prod_{\zeta \subface \eta^\sigma \\ \dims\zeta=1} \x_\zeta\]
and a lifting of $\omega_{\Sigma^\sigma}$ in $A^{d-\dims\sigma}(\Sigma)$ is given by
\[\alpha=\x_\tau= \prod_{\varrho \subface \tau \\ \dims\varrho=1}\x_\varrho.\]

We now claim that the canonical element $\omega_{\Sigma'}$ can be identified with the element $-T^{\dims \sigma} \alpha$ in $A^\bul(\Sigma') \simeq \rquot{A^\bul(\Sigma)[T]}{(\J T+P(T))}$. Indeed, using the identity $P(T) = 0$, we get
\[-T^{\dims \sigma} = \sum_{j=1}^{\dims \sigma} S_{j} T^{\dims \sigma -j}\]
with $S_j$ referring to the $j$-th symmetric function in the variables $\x_\varrho\in A^1(\Sigma')$ for $\varrho$ a ray in $\sigma$. Therefore,
\[ -T^{\dims \sigma} \alpha = \sum_{j=1}^{\dims \sigma} S_{j} \alpha T^{\dims \sigma -j} = \prod_{\varrho \subface \sigma\\ \dims\varrho =1} \x_\varrho \cdot \alpha + \sum_{j=1}^{\dims \sigma -1} S_{j} \alpha T^{\dims \sigma -j}. \]
Since $\omega_{\Sigma^{\sigma}}=\i^*_\sigma(\alpha)$ lives in the top-degree part of $A^\bul(\Sigma^\sigma)$, the products $S_j\alpha$ all belong to $\J$. We infer that the terms of the sum are all vanishing for $j=1,\dots, \dims \sigma -1$, and
\[ -T^{\dims \sigma} \alpha = \prod_{\varrho \subface \sigma\\ \dims\varrho =1} \x_\varrho\cdot \alpha = \x_\sigma \cdot \x_\tau = \omega_{\Sigma} = \omega_{\Sigma'}. \]

We thus obtain the following commutative diagram.
\[ \begin{tikzcd}
A^{d-\dims\sigma}(\Sigma^\sigma) \arrow[r, "\cdot \left(-T^{\dims{\sigma}}\right)", "\sim"'] \arrow[rd, "\sim"{sloped, above}, "\deg"'] & A^d(\Sigma') \arrow[r, "\sim"'] \arrow[d, "\deg", "\sim"{sloped, below}] & A^d(\Sigma) \arrow[ld, "\deg", "\sim"{sloped, above}] \\
& \Z
\end{tikzcd} \]
Let us denote by $\Psi_k(\Sigma')$ the pairing $A^k(\Sigma') \times A^{d-k}(\Sigma') \to \Z$. We use similar notations for $\Sigma$ and $\Sigma^\sigma$. The above diagram allows us to describe the pairing induced by $\Psi_k(\Sigma')$ on the different parts of the decomposition given by the Keel's lemma.
\begin{itemize}
\item Between $A^k(\Sigma)$ and $A^{d-k}(\Sigma)$ we get the pairing $\Psi_k(\Sigma)$.
\item For any integer number $i$, between $A^{k-i}(\Sigma^\sigma)T^i$ and $A^{d-\dims\sigma-k+i}(\Sigma^\sigma)\,T^{\dims\sigma-i}$ we get the pairing $-\Psi_{k-i}(\Sigma^\sigma)$.
\item For integers $i<j$, between $A^{k-i}(\Sigma^\sigma)T^i$ and $A^{d-\dims\sigma-k+j}(\Sigma^\sigma)\,T^{\dims\sigma-j}$ the pairing is trivial since $A^{d-\dims\sigma+j-i}(\Sigma^\sigma)$ is trivial.
\item For any integer $i<\dims\sigma$, the pairing between $A^k(\Sigma)$ and $A^{d-k-i}(\Sigma^\sigma)T^i$ is trivial. Indeed, the product lives in $A^{d-i}(\Sigma^\sigma)T^i$ which is trivial since $d-i > d-\dims\sigma$.
\end{itemize}
We do not need to compute the pairing between the other parts. The bilinear form $\Psi_k(\Sigma')$ can be written in the form of a block matrix consisting of bilinear maps, as follows.
\[ \begin{tikzpicture}
\matrix[matrix of math nodes] {
  & \quad & A^{d-k}(\Sigma) & A^{d-k-\dims\sigma+1}(\Sigma^\sigma) T^{\dims\sigma-1} & \quad\cdots\quad & A^{d-k-1}(\Sigma^\sigma)T \\[.5ex]
A^k(\Sigma) && \Psi_k(\Sigma) & 0 && \\
A^{k-1}(\Sigma^\sigma)T && 0 & -\Psi_{k-1}(\Sigma^\sigma) & \node (zero) {}; & \\
\vdots && \vdots & \node (ast) {}; & \rotatebox{-30}{$\cdots$} & \\
A^{k-\dims\sigma+1}(\Sigma^\sigma) T^{\dims\sigma-1} && 0 &&& -\Psi_{k-\dims\sigma+1}(\Sigma^\sigma) \\};
\draw[gray, thin, dotted] ($(zero)+(-1.7,1)$) --++ (4.2,0) --++ (0,-2.1) --++ (-1,0) -- cycle;
\path (zero) node[shift={(1.4,0)}] {\Large0};
\draw[gray, thin, dotted] ($(ast)+(2.6,-1)$) --++ (-3.4,0) --++ (0,1.4) --++ (1.2,0) -- cycle;
\path (ast) node[shift={(.2,-.5)}] {\LARGE*};
\path (-2.8,-.4) node {$\left(\rule[-1.4cm]{0pt}{2.8cm}\right.$};
\path (7.2,-.4) node {$\left.\rule[-1.4cm]{0pt}{2.8cm}\right)$};
\end{tikzpicture} \]
The matrix is lower triangular. On the diagonal, all the bilinear maps apart from the first one are non-degenerate, since we have $\PD_\chow(\Sigma^\sigma)$. Thus, $\Psi_k(\Sigma')$ is non-degenerate if and only if $\Psi_k(\Sigma)$ is non-degenerate. This shows the properties $\PD_\chow(\Sigma)$ and $\PD_\chow(\Sigma')$ are equivalent provided that $\PD_\chow(\Sigma^\sigma)$ holds.
}

\subsubsection{Proof of Theorem~\ref{thm:chow_shellable}} At this point, we have verified all the needed properties for the shellability of the property $\PD_{\chow}$ in $\unim$, and this concludes the proof of Theorem~\ref{thm:chow_shellable}. \qed

\subsection{Elementary proof of Localization Lemma~\ref{lem:kernel_Z^k}} \label{sec:kernel_Z^k}

In this section, we give a combinatorial self-contained proof of Lemma~\ref{lem:kernel_Z^k}. Although the proof is somehow technical, we hope it could be of independent interest for several reasons. First, the result is used by Adiprasito-Huh-Katz~\cite{AHK} in order to prove the duality between the Chow groups and the Minkowski weights of unimodular simplicial fans. While almost all the proofs in that paper are written to be accessible to a combinatorial audience, this one refers to a more general result taken from~\cite{FMSS} on Chow rings of algebraic varieties admitting a solvable group action. The proof we give is elementary and arguably more transparent. Second, combined with a similar in spirit theorem proved in~\cite{Ami20} for the combinatorial Chow rings of products of graphs, in connection with semi-stable reduction theorem, it suggests a more general phenomenon concerning the combinatorial part of the Chow rings of algebraic varieties that certainly merits a further study.

\subsubsection{Admissible expansions} Let $\Sigma$ be any simplicial rational fan. It will be convenient in the course of the proof to introduce the following notation. For an element $\ell \in M =N^\dual$, let
\[ \x_\ell := \sum_{\rho \in \Sigma_1} \ell(\e_\rho) \x_\rho. \]

Let $B^k$ be the group of homogeneous polynomials of degree $k$ in $\Z[\x_\rho \mid \rho\in\Sigma_1]$. Recall that $I$ is the ideal generated by the products $\x_{\rho_1}\cdots \x_{\rho_s}$, for $s\in \N$, such that $\rho_1, \dots, \rho_s$ do not belong to a same cone in $\Sigma$, and that $J$ is the ideal generated by the elements $\x_\ell$ for $\ell \in M$.

Recall that $Z^k = \bigoplus_{\sigma \in \Sigma_k} \Z\x_\sigma$, which is thus a subgroup of $B^k$. We need to show that $Z^k \to \rquot{B^k}{B^k \cap(I + J)} $ is surjective and to determine its kernel. We will prove the latter, the former becomes clear in the course of the proof; one can also find a proof of the surjectivity in \cite{AHK}.

\smallskip
Let $K^k$ be the subgroup of $B^k$ generated by elements of the form
\[ \sum_{\rho \in \Sigma_1 \\ \rho\sim\tau} \ell(\e_\rho) \x_\rho \x_\tau, \]
for $\tau \in \Sigma_{k-1}$ and $\ell$ an element in $M$ which is orthogonal to $\tau$. Consider an element $a$ of $B^k \cap \bigl(I+J\bigr)$ which belongs to $Z^k$. We want to prove that $a \in K^k$.

\smallskip
First, notice that we can write $a$ as an element of $I$ plus a sum of monomials of the form $\x_\ell \x_{\varrho_2} \cdots \x_{\varrho_k}$ for some $\ell \in M=N^\dual$ and rays $\varrho_2, \dots, \varrho_k$. The first step is to prove the following result.

\begin{claim}\label{claim:admissible_expansion}
Each element $a \in B^k\cap (I + J) $ can be written as a sum consisting of an element in $I$ and a sum consisting of elements of the form $\x_{\ell_1} \x_{\ell_2} \cdots \x_{\ell_s} \x_{\tau}$ for $s\in \N$, $\tau \in \Sigma_{k-s}$, and $\ell_1, \dots, \ell_s \in M$.
\end{claim}

An expansion of the form described in the claim for an element $a \in B^k\cap (I + J)$ will be called \emph{admissible} in the sequel.

\begin{proof}
To prove the above claim, it will be enough to prove the statement for an element $a$ of the form $\x_{\ell_1} \cdots \x_{\ell_s} \x_{\varrho_{s+1}}^{\kappa_{s+1}} \cdots \x_{\varrho_r}^{\kappa_r}$ with $\kappa_{s+1} \geq \dots \geq \kappa_{r}$ and $\varrho_{s+1}, \dots, \varrho_r$ distinct and comparable, \ie, they form the rays of some cone $\tau$ of dimension $r-s$. Actually, the case $s=1$ will already give the result, but considering arbitrary values for $s$ allows to proceed by induction on the lexicographical order of the $k$-tuples of non-negative integer numbers.

To a term of the form $\x_{\ell_1} \cdots \x_{\ell_s} \x_{\varrho_{s+1}}^{\kappa_{s+1}} \cdots \x_{\varrho_r}^{\kappa_r}$ with $\kappa_{s+1} \geq \dots \geq \kappa_r$ we associate the $k$-tuple of integers $(\kappa_{s+1}, \dots, \kappa_r, 0, \dots, 0)$. We now show that an element as above with $\kappa_{s+1} > 1$ can be rewritten as the sum of an element of $I$ and a sum of terms of the above form having a $k$-tuple with strictly smaller lexicographical order.

So suppose $\kappa_{s+1} > 1$. We take a linear form $l$ which takes the value one on $\e_{\varrho_{s+1}}$ and $0$ on the other rays $\varrho_{s+2}, \dots, \varrho_r$. This is doable since $\tau$ is simplicial. Using
\[\x_l = \x_{\varrho_{s+1}} + \sum_{\rho \in \Sigma_1 \\ \rho\not\subface \tau} l(\e_\rho)\x_\rho,\]
we get
\[ \begin{split} \x_{\ell_1}\! \cdots \x_{\ell_s} \x_{\varrho_{s+1}}^{\kappa_{s+1}} \cdots \x_{\varrho_r}^{\kappa_r}
  =&\: \x_{\ell_1}\! \cdots \x_{\ell_s} \x_{l} \x_{\varrho_{s+1}}^{\kappa_{s+1}-1} \x_{\varrho_{s+2}}^{\kappa_{s+2}} \cdots \x_{\varrho_r}^{\kappa_r} \\
  &\quad- \sum_{\rho \in \Sigma_1 \\ \rho \not\subface \tau} l(\e_\rho) \x_{\ell_1} \cdots \x_{\ell_s} \x_{\rho} \x_{\varrho_{s+1}}^{\kappa_{s+1}-1} \x_{\varrho_{s+2}}^{\kappa_{s+2}} \cdots \x_{\varrho_r}^{\kappa_r}.
\end{split} \]
Each term in the right hand side is either in $I$ or is of the form described above with a lower lexicographic order. Proceeding by induction, we get the claim.
\end{proof}

\smallskip
So each element $a\in B^k \cap (I +J)$ is the sum of an element of $I$ plus a sum of terms of the form $\x_{\ell_1}\!\cdots\x_{\ell_s}\x_\tau$ with $1\leq s < k$ and for a cone $\tau \in \Sigma_{k-s}$.

\subsubsection{An auxilary filtration} Using admissible expansions, we now introduce an increasing filtration $\F_\bul$ on $Z^k \cap (I +J)$ as follows.

First, for each $s>0$, denote by $\filt_s$ the group generated by the elements of $B^k$ of the form $\x_{\ell_1} \cdots \x_{\ell_s} \x_\tau$ for some linear forms $\ell_1, \dots, \ell_s \in M$ and for some $\tau \in \Sigma_{k-s}$. Moreover, for $s=0$, define the subgroup $\filt_0 \subseteq \filt_1$ as the one generated by the elements of the form $\x_\ell \x_\tau$ for some $\tau \in \Sigma_{k-1}$ and some linear form $\ell$ orthogonal to $\tau$.

For any $t\geq 0$, let $\H_t := \filt_0 + \dots + \filt_t$ and define $\F_t := Z^k \cap (\H_t + I)$. Note that $\H_t$ is the subgroup generated by the elements which admit an admissible expansion having only terms $\x_{\ell_1} \x_{\ell_2} \cdots \x_{\ell_s} \x_{\tau}$ with $s\leq t$. In this way, we get a filtration
\[ \F_0 \subseteq \F_1 \subseteq \cdots \subseteq \F_{k-1} \subseteq \F_k \subseteq Z^k \cap (I + J). \]

We now prove that all these inclusions are equalities.
\begin{claim} \label{claim:equality}
We have $\F_0 = \F_1 = \cdots = \F_{k-1} = \F_k= Z^k \cap (I + J)$.
\end{claim}
\begin{proof}
By Claim~\ref{claim:admissible_expansion}, we have $\F_k = Z^k \cap (I + J)$, which proves the last equality. We prove all the other equalities.

So fix $s > 0$ and let $a$ be an element of $\F_s = Z^k \cap (\H_s + I)$ admitting an admissible expansion consisting of an element of $I$ plus a sum of terms each in $\filt_1, \dots, \filt_{s-1}$ or $\filt_s$. We need to show that $a \in \F_{s-1}$.

We can assume there is a term of $a$ in this admissible expansion which lies in $\filt_s$, \ie, of the form $\x_{\ell_1} \cdots \x_{\ell_s} \x_\tau$, with $\tau \in \Sigma_{k-s}$. Otherwise, the statement $a \in \F_{s-1}$ holds trivially. Consider all the terms in the admissible expansion of $a$ which are of the form $\x_{\ell'_1} \cdots \x_{\ell'_s} \x_\tau$ with the same cone $\tau$, but with possibly different linear forms $\ell_1', \dots, \ell_s'$. We will prove that the sum of those terms that we denote by $a_\tau$ belongs to $\filt_{s-1} + I$. Applying this to any $\tau \in \Sigma_{k-s}$, we obtain $a \in \F_{s-1}$ and the claim follows.

\smallskip
For each ray $\rho$ of $\tau$, choose a linear form $\ell_{\rho, \tau}$ which takes value one on $\e_\rho$ and zero on other rays of $\tau$. This is again possible since $\Sigma$ is simplicial. Then $\ell_1$, for instance, can be decomposed as the sum
\[ \ell_1 = l_1 + \sum_{\rho \subface \tau \\
\rho \in \Sigma_1} \ell_1(\e_\rho) \ell_{\rho, \tau} \]
with $l_1$ vanishing on $\tau$.

The observation now is that the term $\x_{l_1} \x_{\ell_2} \cdots \x_{\ell_s} \x_\tau$ belongs to $\filt_{s-1} + I$: this is by definition if $s = 1$, and for $s>1$ is obtained by expanding the product
\[\x_{l_1} \x_{\ell_2} \cdots \x_{\ell_s} \x_\tau = \sum_{\rho \in \Sigma_1\\
\rho \not \subface \tau} l_1(\e_\rho) \x_\rho \x_{\ell_2} \cdots \x_{\ell_s} \x_\tau \]
and by observing that each term in the right hand side is either in $I$ or is in $\filt_{s-1}$. (Note that the sum is on rays $\rho \not \subface \tau$ because $l_1$ vanishes on $\tau$.) Hence, in proving that $a_\tau$ belongs to $\filt_{s-1}+I$, we can ignore this term.

Decomposing in the same way each $\x_{\ell_i}$ which appear in a term in the initial admissible expansion of $a_\tau$, we can rewrite $a_\tau$ as a sum of terms which already belong to $\H_{s-1}+I$ plus a sum of terms of the form
\begin{equation}\label{eq:atau}\prod_{\rho \subface \tau\\\rho\in \Sigma_1}\x_{\ell_{\rho, \tau}}^{\kappa_\rho} \x_\tau, \qquad \kappa_\rho \in \Z_{\geq 0}.
\end{equation}
These all together give a new admissible expansion of $a$ in which $a_\tau$ (defined as before in the new admissible expansion) is a sum of the terms of the form in~\eqref{eq:atau}. We will work from now on with this admissible expansion of $a$.

Fix non-negative integers $\kappa_\rho$ for rays $\rho$ of $\tau$ whose sum is $s$. Denote by $\bigl[\prod_{\rho \in \tau} \x_\rho^{\kappa_\rho + 1} \bigr] a$ the coefficient of this monomial in $a$ written as sum of monomials in $B^k$. Since $a$ belongs to $Z^k$, all the monomials of $a$ are square-free. But the product $\prod_{\rho \in \tau} \x_\rho^{\kappa_\rho + 1}$ is not square-free since $s>0$. Hence, the corresponding coefficient is zero.

Consider now the monomial $\prod_{\rho \in \tau} \x_\rho^{\kappa_\rho + 1}$ and let us see in which terms in the admissible expansion of $a$ it can occur. Such a monomial cannot appear in a term of the admissible expansion of $a$ which belongs to $I$. It cannot neither be in a term of $a$ of the form $\x_{\lambda_1} \cdots \x_{\lambda_{s'}} \x_\sigma$ with $\sigma \in \Sigma_{k-s'}$ with $s'<s$ nor with $\sigma \in \Sigma_{k-s}$ and $\sigma \neq \tau$.

It follows that the monomial $\prod_{\rho \in \tau} \x_\rho^{\kappa_\rho + 1}$ can only appear in the terms of $a_\tau$. More precisely, by the definition of $\ell_{\rho,\tau}$, it can only appear in each term of the form $\prod_{\rho \subface \tau}\x_{\ell_{\rho, \tau}}^{\kappa_\rho} \x_\tau$ for the chosen $\kappa_\rho$, and with coefficient one in each term. Hence, the sum of these terms have to cancel out. This proves that $a_\tau$ is zero in the new admissible expansion, which shows that $a \in \F_{s-1}$ and the claim follows.
\end{proof}

\subsubsection{End of the proof} We can now finish the proof of the lemma.

\begin{proof}[Proof of Lemma~\ref{lem:kernel_Z^k}]
Recall that $K^k$ is the set generated by elements of the form
\[ \sum_{\rho \in \Sigma_1 \\ \rho\sim\tau} \ell(\e_\rho) \x_\rho \x_\tau, \]
where $\tau \in \Sigma_{k-1}$ and $\ell$ is orthogonal to $\tau$. The following facts are then clear:
\[ \F_0 \subseteq K^k + I, \quad K^k \subseteq Z^k, \quad \text{and } I \cap Z^k = \{0\}. \]
Since $\F_0 \subseteq Z^k$, we deduce that $\F_0 \subseteq K^k$. Applying Claim~\ref{claim:equality}, we get $Z^k \cap (I + J) = \F_0$ which implies $Z^k \cap (I + J) \subseteq K^k$. The inclusion $K^k \subseteq Z^k \cap (I + J)$ is obvious. So we get $Z^k \cap (I + J) = K^k$ which concludes the proof.
\end{proof}

\section{Hodge isomorphism for unimodular fans}\label{sec:hodge_isomorphism}

In this section we prove the Hodge isomorphism theorem which relates the Chow group of unimodular fans to their cohomology groups. As a byproduct of the methods we develop in this section, we obtain an alternate characterization of tropical smoothness.

\smallskip
Unless otherwise stated, the fans which appear in this section are not assumed to be tropical.

\subsection{Statement of the main theorems}

The main theorem we prove in this section is the following.

\begin{thm}[Hodge isomorphism for unimodular fans] \label{thm:ring_morphism}
Let $\Sigma$ be a saturated unimodular fan in $N_\R$. For any integer $p$, there is a well-defined isomorphism
\[ \begin{array}{ccc}
\Psi\colon H^{p,p}(\comp\Sigma) & \overset\sim\longrightarrow & A^p(\Sigma) \\
\alpha               & \longmapsto & \displaystyle\sum_{\sigma \in \Sigma_p} \alpha_\sigma(\nu_\sigma)\x_\sigma
\end{array} \]
which induces a ring morphism $H^{\bul,\bul}(\comp\Sigma) \to A^\bul(\Sigma)$ by mapping $H^{p,q}(\comp\Sigma)$ to zero in bidegrees $p\neq q$. Moreover, $H^{p,q}(\comp\Sigma)$ is trivial in bidegrees $p < q$ and bidegrees $(p,0)$ for $p > 0$.
\end{thm}

\smallskip
As an immediate corollary, we get the following theorem.

\begin{thm}[Hodge isomorphism for smooth unimodular tropical fans]\label{thm:Hodge_isomorphism}
Let $\Sigma$ be a saturated unimodular tropical fan in $N_\R$. Suppose in addition that $\Sigma$ is tropically smooth. Then we get an isomorphism of rings $A^\bul(\Sigma) \simto H^{\bul, \bul}(\comp\Sigma)$. Moreover, $A^\bul(\Sigma)$ verifies the Poincaré duality. In particular, all these groups are torsion-free.
\end{thm}

\begin{proof}
This follows from the discussion in Remark \ref{rem:torsion_duality_chow}.
\end{proof}

This saturation hypothesis is not always possible to obtain (see Example \ref{ex:non_saturated_smooth_fan}). The following dual version of the above theorem, which does not use the saturation hypothesis, is sometimes more useful for applications.

\begin{thm}
Let $\Sigma$ be any unimodular fan in $N_\R$. Then the natural map $\MW_p(\Sigma) \to H_{p,p}(\comp\Sigma)$ is an isomorphism. Moreover, $H_{p,q}(\comp\Sigma)$ is trivial in bidegrees $p < q$ and bidegrees $(p,0)$ for $p > 0$.
\end{thm}

\begin{proof}
Once again, this follows from the discussion in Remark \ref{rem:torsion_duality_chow} (and from a direct computation for the vanishing of $H^{p,0}(\comp\Sigma)$ for $p > 0$).
\end{proof}

In the smooth case, we get the following corollary.

\begin{thm}[Hodge conjecture for tropical fans] \label{thm:Hodge_conjecture}
If\/ $\Sigma$ is a smooth unimodular tropical fan, $H_{p,q}(\comp\Sigma)$ and $H^{p,q}(\comp\Sigma)$ are trivial for $p \neq q$. Moreover, to each element $\alpha$ in $H^{p,p}(\comp\Sigma)$ one can associate a tropical cycle whose class in $H_{d-p,d-p}(\comp\Sigma)$ is the Poincaré dual of $\alpha$.
\end{thm}

\begin{proof}
The tropical cycle is obtained via the isomorphism
\[ H^{p,p}(\comp\Sigma) \simeq H_{d-p, d-p}(\comp\Sigma) \simeq \MW_{d-p}(\Sigma). \]
The details can be found in Remark \ref{rem:torsion_duality_chow}.
\end{proof}

Before going through the proof of Theorem~\ref{thm:ring_morphism}, we state some clarifying remarks.

\begin{remark}[Non saturated case] \label{rem:ring_morphism_non_saturated}
The proof of Theorem \ref{thm:ring_morphism} gives as well results in the case $\Sigma$ is not saturated. In this situation, we still have a surjective ring morphism $A^\bul(\Sigma) \to H^{\bul,\bul}(\Sigma)$, but this is never injective. Indeed, the kernel is a nontrivial torsion group. In particular, we get an isomorphism for rational coefficients, as well as a dual isomorphism for integral coefficients: $A^\bul(\Sigma)^\dual \simeq H^{\bul,\bul}(\comp\Sigma)^\dual$. In any case, the last part of the theorem concerning the vanishing are still valid.
\end{remark}

\begin{remark}[Torsion-freeness and duality for Chow rings] \label{rem:torsion_duality_chow} Let $\Sigma$ be a unimodular fan. By the universal coefficient theorem applied to the tropical chain and cochain complexes of $\comp \Sigma$, for each pair of non-negative integers $p,q$, we get the following exact sequences
\begin{align}
\label{eq:univ_coeff_homology}
0\to \Ext{H^{p, q+1}(\comp\Sigma)}{\Z} \to H_{p,q}(\comp \Sigma) \to H^{p,q}(\comp \Sigma)^\dual \to 0,& \quad\text{and} \\
\label{eq:univ_coeff_cohomology}
0\to \Ext{H_{p, q-1}(\comp\Sigma)}{\Z} \to H^{p,q}(\comp \Sigma) \to H_{p,q}(\comp \Sigma)^\dual \to 0.&
\end{align}
From the first exact sequence, using the vanishing result stated in Theorem~\ref{thm:ring_morphism}, we obtain an isomorphism $H_{p,q}(\comp\Sigma) \simeq H^{p,q}(\comp\Sigma)^\dual$ for $p\leq q$. This shows the vanishing of $H_{p,q}(\comp \Sigma)$ for $p<q$ and an isomorphism $H_{p,p}(\comp \Sigma) \simeq H^{p,p}(\comp \Sigma)^\dual$. This in turn implies in particular that $H_{p,p}(\Sigma)$ is torsion-free.

\smallskip
In the case $\Sigma$ is saturated, the dual of the map from $H^{p,p}(\comp \Sigma)$ to $A^p(\Sigma)$ depicted in the statement of Theorem~\ref{thm:ring_morphism} is the natural morphism from the group of Minkowski weights $\MW_p(\Sigma)$ to the tropical homology group $H_{p,p}(\comp\Sigma)$. This map is obtained by sending a Minkowski weight to its closure in $\comp\Sigma$, which is a tropical cycle of dimension $p$, and then taking the corresponding homology class in $H_{p,p}(\comp\Sigma)$, as described in Section \ref{subsec:cycle_class_map}; for the cycle class map, see for example~\cites{MZ14, Sha-thesis, JRS18, GS19}. Theorem~\ref{thm:ring_morphism} implies that this map is an isomorphism. By Remark \ref{rem:ring_morphism_non_saturated}, this isomorphism still exists in the non saturated case.

\smallskip
In the case $\Sigma$ is a smooth saturated unimodular tropical fan, we get an isomorphism $H^{p,q}(\comp \Sigma) \simeq H_{d-p,d-q}^\BM(\comp\Sigma) = H_{d-p,d-q}(\comp \Sigma)$, $p,q \in \Z_{\geq 0}$, which implies as well the vanishing of $H_{p,q}(\comp \Sigma)$ for $p>q$. Applying the exact sequence \eqref{eq:univ_coeff_cohomology}, we infer the isomorphism $H^{p,q}(\comp \Sigma) \simeq H_{p,q}(\comp \Sigma)^\dual$. This in particular implies that the cohomology groups $H^{p,p}(\comp \Sigma)$ and so $A^p(\Sigma)$ are all torsion-free, we have $A^p(\Sigma) \simeq \MW_p(\Sigma)^\dual$, and the Chow ring verifies the Poincaré duality.
\end{remark}

\begin{remark}
The smoothness assumption in the statement of Theorem~\ref{thm:Hodge_isomorphism} is needed, and the result is not true, in general, for any saturated unimodular tropical fan. Example \ref{subsec:cube} is a saturated unimodular tropical fan $\Sigma$ with $H^{2,1}(\comp\Sigma)$ of rank two.
\end{remark}

\smallskip
In order to prove Theorem \ref{thm:ring_morphism}, we establish an alternative way of computing the cohomology of $\comp\Sigma$ by using its decomposition into \emph{hypercubes}. The proof of the first half then follows quite directly from this description by using the localization Lemma~\ref{lem:kernel_Z^k} proved in the previous section. To show the second half, concerning the products, we provide an explicit description of the inverse of the map $\Psi$ from $A^p(\Sigma)$ to $H^{p,p}(\comp \Sigma)$, a result which might be of independent interest.

\smallskip
For non-negative integers $p,q$, we define
\[ \cubC^{p,q}(\comp\Sigma) := \bigoplus_{\dims\sigma = q} \SF^{p-q}(\infty_\sigma). \]
We get a cochain complex for each integer $p$
\[ \cubC^{p,\bullet}\colon \quad \dots\longrightarrow \cubC^{p, q-2}(\comp\Sigma) \xrightarrow{\ \cubd^{q-2}\ } \cubC^{p,q-1}(\comp\Sigma) \xrightarrow{\ \cubd^{q-1}\ } \cubC^{p,q}(\comp\Sigma) \longrightarrow \cdots \]
where the differential is given by the sum of the following maps for $\tau \ssubface \sigma$
\[ \begin{array}{ccc}
  \SF^{k}(\infty_\tau) & \longrightarrow & \SF^{k-1}(\infty_\sigma), \\
  \alpha               & \longmapsto     & \pi_*(\iota_{\e^\tau_\sigma}(\alpha)).
\end{array} \]
Here $\pi\colon N^\tau \to N^\sigma$ is the natural projection.

\begin{thm} \label{thm:cohomology_compactification}
The cohomology of\/ $\bigl(\cubC^{p,\bul}, \cubd \bigr)$ is $H^{p,\bul}(\comp\Sigma)$.
\end{thm}

The theorem is obtained by exploiting the decomposition of $\comp\Sigma$ into hypercubes. This explains the notation for this new complex. Note that one can also use the dual complex which compute the homology of $\comp\Sigma$.

\smallskip
The proof of Theorem \ref{thm:cohomology_compactification} uses an interesting double complex. A deeper exploitation of this double complex leads to the proof of the following alternate characterization of smoothness.

\begin{thm}[Alternate characterization of tropical smoothness] \label{thm:smoothness_alternate}
A unimodular tropical fan $\Sigma$ is smooth if and only if, for any face $\sigma\in\Sigma$, $\comp\Sigma^\sigma$ verifies the Poincaré duality.
\end{thm}

Note that, in particular, this theorem gives a new proof of Theorem \ref{thm:PD_canonical_compactification} which states that if $\Sigma$ is smooth, then $\comp\Sigma$ verifies the Poincaré duality.

\smallskip
The rest of this section is devoted to the proof of the above results.

\subsection{Proof of Theorem \ref{thm:cohomology_compactification}}

We start with the proof of Theorem \ref{thm:cohomology_compactification}. We fix a non-negative integer $p$.

\subsubsection{The fine double complex} The idea behind the proof is quite natural. The original tropical cochain complex $C^{p,\bul}(\comp\Sigma)$ has a finer decomposition. Namely, for a face $\cube^\tau_\sigma = \comp C^\tau_\sigma$ of $\comp\Sigma$, instead of using a grading just following its dimension $\dims\sigma-\dims\tau$, which is what is done in $C^{p,\bul}(\comp\Sigma)$, one can remember both $\dims\sigma$ and $\dims\tau$. In this way, the tropical cochain complex $C^{p,\bul}(\comp\Sigma)$ can be unfolded and identified with the total complex of the following double complex:
\[ \Ep{p}^{a,b} := \bigoplus_{\tau\subface\sigma \\ \dims\sigma = a \\ \dims\tau = -b} \SF^p(\cube^\tau_\sigma) \]
where the differentials are given by the ones between the coefficient groups $\SF^p$ of faces of $\comp\Sigma$ which define the differentials of $C^{p,\bul}(\comp\Sigma)$. So we have $C^{p,\bul}(\comp\Sigma) = \Tot^\bul(\Ep{p}^{\bul, \bul})$.

The cohomology of $C^{p,\bul}(\comp\Sigma)$ can be calculated by using the spectral sequence associated to the filtration given by the columns of $\Ep{p}^{\bul, \bul}$. (We will use below the other filtration, the one given by the rows, in the proof the other stated results.)

 We denote by $\EpI{p}^{\bul,\bul}_0$ the $0$-th page of this spectral sequence which has abutment
\[ \EpI{p}_0^{\bul,\bul} \Longrightarrow H^{p,\bul}(\comp\Sigma). \]

\subsubsection{Computation of the first page} The proof of Theorem~\ref{thm:cohomology_compactification} is given by the following claim.
\begin{claim} \label{claim:degeneration_spectral_sequence}
On the first page of the spectral sequence, we have
\[\EpI{p}_1^{a,b} =\begin{cases} \cubC^{p,a}(\comp\Sigma) & \text{for $b = 0$,} \\
0 & \text{otherwise}.
\end{cases}\]
Moreover, the differential on the first page for $b=0$ coincides with the differential $\cubd^\bul$ of $\cubC^{p,\bul}(\comp\Sigma)$.
\end{claim}

We prove the claim. Fix an integer $a$. Then the $a$-th column in page zero of the spectral sequence is
\[ \EpI{p}_0^{a, \bul}\colon \quad \dots \longrightarrow \bigoplus_{\tau\subface\sigma \\ \dims\sigma = a \\ \dims\tau = -b+1} \SF^p(\cube^\tau_\sigma) \longrightarrow \bigoplus_{\tau\subface\sigma \\ \dims\sigma = a \\ \dims\tau = -b} \SF^p(\cube^\tau_\sigma) \longrightarrow \bigoplus_{\tau\subface\sigma \\ \dims\sigma = a \\ \dims\tau = -b-1} \SF^p(\cube^\tau_\sigma) \longrightarrow \cdots \]
This complex can be decomposed as a direct sum running on the faces $\sigma \in \Sigma_a$ of the complexes
\[ \dots \longrightarrow \bigoplus_{\tau\subface\sigma \\ \dims\tau = -b+1} \SF^p(\cube^\tau_\sigma) \longrightarrow \bigoplus_{\tau\subface\sigma \\ \dims\tau = -b} \SF^p(\cube^\tau_\sigma) \longrightarrow \bigoplus_{\tau\subface\sigma \\ \dims\tau = -b-1} \SF^p(\cube^\tau_\sigma) \longrightarrow \cdots \]
Denote this last complex by $\EpI{p}_0^{\sigma,\bul}$, $\sigma \in \Sigma$, and let $\EpI{p}_1^{\sigma,\bul}$ be its cohomology. We thus have
\begin{equation}\label{eq:decomposition}
\EpI{p}_0^{a,\bul} = \bigoplus_{\sigma \in \Sigma_a} \EpI{p}_0^{\sigma,\bul}, \qquad \EpI{p}_1^{a,\bul} = \bigoplus_{\sigma \in \Sigma_a} \EpI{p}_1^{\sigma,\bul}.
\end{equation}
These latter cohomology groups are given by the following lemma.

\begin{lemma}[Hypercube vanishing lemma] \label{lem:hypercube_vanishing}
For any face $\sigma$ of\/ $\Sigma$, the cohomology $\EpI{p}_1^{\sigma,\bul}$ of the complex $\EpI{p}_0^{\sigma,\bul}$ is given by
\[ \EpI{p}_1^{\sigma,b} \simeq \begin{cases} \SF^{p-\dims\sigma}(\infty_\sigma) & \text{ if $b=0$,} \\
0 & \text{otherwise.}
\end{cases} \]
This isomorphism is moreover induced by the natural map
\[ \begin{array}{ccccc}
  \EpI{p}_0^{\sigma,0} &\simeq& \SF^p(\comp\sigma) &\longrightarrow& \SF^{p-\dims\sigma}(\infty_\sigma), \\
  && \alpha & \longmapsto & \pi_*(\iota_{\nu_\sigma}(\alpha))
\end{array} \]
with $\pi$ referring to the projection $N\to N^\sigma$, extended to the exterior algebra and restricted to the corresponding subspaces $\SF^p$.
\end{lemma}

Let us assume the lemma for the moment and finish the proof of the claim.

\begin{proof}[Proof of Claim~\ref{claim:degeneration_spectral_sequence}]By the above lemma, and in view of the decomposition given in~\ref{eq:decomposition}, the first page of the spectral sequence is concentrated in the 0-th row and is given by
\[ \EpI{p}_1^{a,0} = \bigoplus_{\dims\sigma = a} \SF^{p-a}(\infty_\sigma), \]
which is exactly $\cubC^{p,a}(\comp\Sigma)$. It thus remains to check that the differentials coincide. Let $\tau\ssubface\sigma$ and denote by $\d_{\tau\ssubface\sigma}\colon \SF^{p-\dims\tau}(\infty_\tau) \to \SF^{p-\dims\sigma}(\infty_\sigma)$ the corresponding part of the differential in $\EpI{p}^{\bul,0}$. By Lemma~\ref{lem:hypercube_vanishing}, we have the following commutative diagram.
\[ \begin{tikzcd}
\SF^p(\comp\tau) \rar{\sign(\tau,\sigma)\i^*}\dar{\pi^\tau_*\circ\,\iota_{\nu_\tau}}& \SF^p(\comp\sigma) \dar{\pi^\sigma_*\circ\,\iota_{\nu_\sigma}} \\
\SF^{p-\dims\tau}(\infty_\tau) \rar{\d_{\tau\ssubface\sigma}}& \SF^{p-\dims\sigma}(\infty_\sigma)
\end{tikzcd} \]
We infer that the bottom map should be given by $\pi_*\circ\iota_{\e^\tau_\sigma}$ with $\pi\colon N^\tau \to N^\sigma$, \ie, by the differential of $\cubC^{p,\bul}(\comp\Sigma)$. This concludes the proof of the claim.
\end{proof}

\subsubsection{Proof of the hybercube vanishing lemma} It remains to prove Lemma \ref{lem:hypercube_vanishing}. For faces $\tau \subface \sigma$, we have isomorphisms
\[ \SF^p(\cube^\tau_\sigma) \simeq \bigoplus_{p_1+p_2=p} \bigwedge^{p_1}N^{\tau\,\dual}_\sigma \otimes \SF^{p_2}(\infty_\sigma). \]
These isomorphisms are not canonical, nevertheless, we can choose them in a compatible way. In fact, the terms appearing in the above decomposition of $\SF^p_\Sigma(\cube^\tau_\sigma)$ correspond to the graded pieces associated to a natural filtration on $\SF^p_\Sigma(\cube^\tau_\sigma)$, that we name the \emph{toric weight filtration}, which plays an important role in the results we prove in~\cite{AP-tht}. In order to simplify the presentation here, we do not detail this filtration here and refer to the relevant part of \cite{AP-tht} for more information.

We can decompose $\EpI{p}^{\sigma,\bul}_0$ as a direct sum of complexes
\[ \EpI{p}^{\sigma,\bul}_0 \simeq \bigoplus_{p_1+p_2=p} \Bigl( \dots
  \longrightarrow \bigoplus_{\tau\supface\sigma \\ \dims\tau = -b+1} \bigwedge^{p_1}N^{\tau\,\dual}_\sigma
  \longrightarrow \bigoplus_{\tau\supface\sigma \\ \dims\tau = -b} \bigwedge^{p_1}N^{\tau\,\dual}_\sigma
\longrightarrow \cdots \Bigr)\ \otimes \SF^{p_2}(\infty_\sigma). \]
Since $\sigma$ is unimodular, the cochain complex appearing in the above summand is isomorphic to $C_c^{p_1,\bul+\dims\sigma}(\eR^{\dims\sigma})$ for the natural subdivision $\{\eR, \infty\}^{\dims\sigma}$:
\[ \EpI{p}^{\sigma,\bul}_0 \simeq \bigoplus_{p_1+p_2=p} C_c^{p_1,\bul+\dims\sigma}(\eR^{\dims\sigma}) \otimes \SF^{p_2}(\infty_\sigma). \]
Hence, the cohomology is given by
\[ \EpI{p}^{\sigma,\bul}_1 \simeq \bigoplus_{p_1+p_2=p} H_c^{p_1,\bul+\dims\sigma}(\eR^{\dims\sigma}) \otimes \SF^{p_2}(\infty_\sigma). \]

The cohomology of $H_c^{\bul,\bul}(\eR^k)$ is easy to compute: this is only nontrivial in bidegree $(k,k)$ where it is isomorphic to $\Z$. It follows that $\EpI{p}_1^{\sigma,b}$ is trivial unless $b=0$ in which case it becomes equal to $\SF^{p-\dims\sigma}(\infty_\sigma)$. The last part of the lemma follows from a cautious study of the different isomorphisms, that we omit here. This concludes the proof of Lemma \ref{lem:hypercube_vanishing}. \qed

\subsubsection{Proof of Theorem \ref{thm:cohomology_compactification}} This is a direct consequence of Claim~\ref{claim:degeneration_spectral_sequence}, which shows that the spectral sequence degenerates at page two, and, moreover, we have
\[ H^{p,q}(\comp \Sigma) \simeq \Ep{p}_2^{q,0} = H^q\bigl(\cubC^{p,\bul}(\comp\Sigma)\bigr). \pushQED{\qed}\qedhere \]

\subsection*{Proof of Theorem \ref{thm:smoothness_alternate}} We have so far used only one of the two spectral sequences coming from the double complex. The spectral sequence $\EpII{p}^{\bul,\bul}$ given by the filtration by rows will allow to prove Theorem \ref{thm:smoothness_alternate}. We explain this now before going through the proof of the hodge isomorphism theorem.

\smallskip
We first prove the forward direction. Assume that $\Sigma$ is a smooth unimodular tropical fan. Recall that the $b$-th row of the double complex is
\[ \EpII{p}_0^{\bul, b}\colon \quad \dots \longrightarrow
  \bigoplus_{\tau\subface\sigma \\ \dims\sigma = a-1 \\ \dims\tau = -b} \SF^p(\cube^\tau_\sigma) \longrightarrow
  \bigoplus_{\tau\subface\sigma \\ \dims\sigma = a \\ \dims\tau = -b} \SF^p(\cube^\tau_\sigma) \longrightarrow
  \bigoplus_{\tau\subface\sigma \\ \dims\sigma = a+1 \\ \dims\tau = -b} \SF^p(\cube^\tau_\sigma) \longrightarrow
\cdots \]
We recognize the cochain complex for the cohomology with compact support of the star fans of codimension $-b$:
\[ \EpII{p}_0^{\bul, b} = \bigoplus_{\dims\tau = -b} C_c^{p,\bul+b}(\Sigma^\tau). \]
All the star fans are smooth by assumption, and so $H_c^{p,q}(\Sigma^\tau)$ is trivial unless $q=d-(-b)$, and for this $q$, we have $H_c^{p,d+b}(\Sigma^\tau) \simeq \SF_{d-p+b}(\infty_\tau)$. Hence, the first page is given by
\[ \EpII{p}_1^{a,b} = \begin{cases}
  \bigoplus_{\dims\tau = -b} \SF_{d-p+b}(\infty_\tau) & \text{if $a = d$,} \\
  0 & \text{otherwise.}
\end{cases} \]
Stated differently, $\EpII{p}^{a,b}$ is trivial unless $a=d$ and
\[ \EpII{p}^{d,\bul} = \hcubC_{d-p,-\bul}(\comp\Sigma), \]
where $\hcubC_{d-p,\bul}(\comp\Sigma)$ denotes the dual complex of $\cubC^{d-p,\bul}(\comp\Sigma)$. Theorem \ref{thm:cohomology_compactification} combined with the universal coefficient theorem imply that the homology of $\hcubC_{d-p,\bul}(\comp\Sigma)$ computes $H_{d-p,\bul}(\comp\Sigma)$. Hence, $\EpII{p}^{\bul,\bul}$ degenerates at page two, and we have
\[ \EpII{p}^{\bul,\bul} \Longrightarrow H_{d-p,d-\bul}(\comp\Sigma). \]

On the other hand, by the definition of the double complex, we know that
\[ \EpII{p}^{\bul,\bul} \Longrightarrow H^{p,\bul}(\comp\Sigma). \]
Therefore, we get the isomorphism of the Poincaré duality $H^{p,q}(\comp\Sigma) \simeq H_{d-p,d-q}(\comp\Sigma)$; more precisely, we infer that
\[ H^{p,p}(\comp\Sigma) \simeq \EpI{p}_\infty^{p,0} \simeq \EpII{p}_\infty^{d,p-d} \simeq H_{d-p,d-p}(\comp\Sigma), \]
and every other terms of $\EpI{p}^{\bul,\bul}_\infty$ and $\EpII{p}^{\bul,\bul}_\infty$ are trivial. This concludes the forward direction.

\begin{remark}
In the same way that $\EpI{p}_\infty^{p,0}$ is naturally isomorphic to $A^p(\Sigma)$ (in the saturated case), for any unimodular fan $\Sigma$ we have the following isomorphisms
\[ H_{d-p,d-p}(\comp\Sigma) \simeq \hcubC_{p,d-p}(\comp\Sigma) = \ker\Bigl(\bigoplus_{\dims\sigma=d-p}\SF_0(\infty_\sigma) \to \bigoplus_{\dims\tau=d-p-1}\SF_1(\infty_\tau)\Bigr) \simeq \MW_{d-p}(\Sigma). \]

Moreover, if $\Sigma$ is smooth as above, the study of the spectral sequences gives us a new natural cycle class map $H^{p,p}(\comp\Sigma) \simto \MW_{d-p}(\Sigma)$. To a cocycle $a$ in $C^{p,p}(\comp\Sigma)$, one can associate the Minkowski weight $(\Sigma_{(d-p)}, w)$ where $w(\sigma) = \langle a, \nu_{\comp\Sigma_\infty^\sigma} \rangle$ for any $\sigma \in \Sigma_{d-p}$.
\end{remark}

\smallskip
We now prove the other direction. Assume that $\comp\Sigma^\sigma$ verifies the Poincaré duality for any face $\sigma$ of $\Sigma$. By induction, we can assume that $\Sigma^\sigma$ is smooth for any proper face $\sigma$ of $\Sigma$.

We have a natural inclusion of complexes $C_c^{p,\bul}(\Sigma) \simeq \Ep{p}^{\bul,0} \hookrightarrow \Ep{p}^{\bul,\bul}$. Denote by $\Ep{p}^{\bul,\bul<0}$ the cokernel. We get a short exact sequence of complexes:
\begin{equation} \label{eqn:cut_E}
0 \longrightarrow C_c^{p,\bul}(\Sigma) \longrightarrow \Tot^{\bul}(\Ep{p}^{\bul,\bul}) \longrightarrow \Tot^{\bul}(\Ep{p}^{\bul,\bul<0}) \longrightarrow 0.
\end{equation}
We already know the cohomology of the first two terms. For the last one, we can compute the first page as before (since the proper star fans are smooth) and we get $\Ep{p}_1^{d,\bul<0} = \hcubC_{d-p,-\bul>0}(\comp\Sigma)$ where
\[ \hcubC_{d-p,\bul>0}(\comp\Sigma) := \coker\Bigl(F_{d-p}(\conezero) \hookrightarrow \hcubC_{d-p,\bul}(\comp\Sigma)\Bigr), \]
and every other terms are trivial. Assume for now that $p < d-1$. The case $p=d$ is trivial and the case $p=d-1$ is studied later below. Then we get
\[ H^k(\hcubC_{d-p,\bul>0}(\comp\Sigma)) = \begin{cases}
  H_{d-p,d-p}(\comp\Sigma) & \text{if $k = d-p$,} \\
  F_{d-p}(\conezero)       & \text{if $k = 1$,} \\
  0                        & \text{otherwise.}
\end{cases} \]

We are now ready to study the long exact sequence associated to \eqref{eqn:cut_E}:
\[ \cdots \longrightarrow H_c^{p,q}(\Sigma) \longrightarrow H^{p,q}(\comp\Sigma) \longrightarrow H^q(\hcubC_{d-p,(d-\bul)>0}(\comp\Sigma)) \longrightarrow H_c^{p,q+1}(\Sigma) \longrightarrow \cdots \]
The two middle terms are trivial unless $q = p$ or $q = d-1$. In the case $q=p$, we get a map $H^{p,p}(\comp\Sigma) \to H_{d-p,d-p}(\comp\Sigma)$. By assumption, this is an isomorphism, and we can just replace both terms by $0$ in the long exact sequence. For $q=d-1$, $H^{p,d-1}(\comp\Sigma) = H^{p,d}(\comp\Sigma) = 0$ and we get
\[ 0 \longrightarrow F_{d-p}(\conezero) \longrightarrow H_c^{p,d}(\Sigma) \longrightarrow 0. \]
All the other terms of $H_c^{p,\bul}(\Sigma)$ vanish. Using the universal coefficient theorem, we get $H^\BM_{p,d}(\Sigma) \simeq H_c^{p,d}(\Sigma)^\dual \simeq F^{d-p}(\Sigma)$, and $H^\BM_{p,q}(\Sigma) = 0$ for $q \neq d$. This is exactly the Poincaré duality for (the homology with coefficients in) $\SF_p$.

It remains to treat the case $p = d-1$. In this case, $\hcubC_{1,\bul>0}(\Sigma)$ contains only one nontrivial term: $\hcubC_{1,1>0}(\Sigma) = \bigoplus_{\varrho\in\Sigma_1} F_0(\infty_\varrho)$. Its cohomology is identical. In the above long exact sequence, the only interesting part is the first row of the following diagram.
\[ \begin{tikzcd}
0 \rar& H_c^{d-1,d-1}(\Sigma) \rar\dar& H^{d-1,d-1}(\comp\Sigma) \rar\dar{\vsim}& \bigoplus_{\varrho\in\Sigma_1}F_0(\infty_\varrho) \rar\dar{\vsim}& H_c^{d-1,d}(\Sigma) \rar\dar& 0 \\
 & 0 \rar& \MW_1(\Sigma) \rar& \W_1(\Sigma) \rar& \SF_1(\Sigma) \rar& 0
\end{tikzcd} \]
Let us describe the second row and the vertical maps. We know that $H^{d-1,d-1}(\comp\Sigma) \simeq H_{1,1}(\comp\Sigma) \simeq \MW_1(\Sigma)$. Here, $\W_1(\Sigma) := \Z^{\Sigma_1} \simeq \bigoplus_{\varrho\in\Sigma_1}F_0(\infty_\varrho)$ denotes the group of one dimensional weights of $\Sigma$. Note also that we have a natural map $H_c^{d-1,d}(\Sigma) \to \SF_1(\Sigma)$: it is induced by the dual of the map $\SF^1(\Sigma) \to H^\BM_{d-1,d}(\Sigma)$. By the definition of the Minkowski weights, the second row is a short exact sequence. The five lemma implies that the first and the last vertical maps are isomorphisms. Therefore, $H_c^{d-1,q}(\Sigma)$ is trivial for $q\neq d$, and $H_c^{d-1,q}(\Sigma) \simeq \SF_1(\Sigma)$. Using the universal coefficient theorem as above, and dualizing, we get the Poincaré duality for $H^\BM_{d-1,\bul}(\Sigma)$. All together we obtain the Poincaré duality for $H^\BM_{\bul,\bul}(\Sigma)$.  \qed

\subsection{Proof of Theorem \ref{thm:ring_morphism}} We conclude with the proof of Theorem \ref{thm:ring_morphism}, given in this section and the next one.

First, we note that by Theorem \ref{thm:cohomology_compactification}, we have $H^{p,q}(\comp\Sigma) \simeq H^q(\cubC^{p,\bul}(\comp\Sigma))$. Notice as well that for $p < q$, $\cubC^{p,q}(\comp\Sigma)$ is trivial. Hence $H^{p,q}(\comp\Sigma)$ is trivial for $p < q$.

\smallskip
For $q=0$ and $p>0$, we get
\[ H^{p,0}(\comp\Sigma) \simeq \ker\bigl(\SF^p(\conezero) \to \bigoplus_{\varrho \in \Sigma_1} \SF^{p-1}(\infty_\varrho)\bigr). \]
Let $\alpha \in \SF^p(\conezero)$ be an element of the kernel. Then $\iota_{\e_\varrho}(\alpha) = 0$ for any $\varrho \in \Sigma_1$. Let $V \subseteq N_\R$ be the subspace spanned by $\Sigma$. Then $\alpha \in \bigwedge^p V^\dual$. Since the vectors $\e_\varrho$, $\varrho \in \Sigma_1$, span $V$, any element of $\bigwedge^p V^\dual$ whose contraction by each $\e_\varrho$ is trivial must be trivial. Hence, we obtain $\alpha = 0$ and $H^{p,0}(\comp\Sigma) = 0$.

\smallskip
For $p=q$, we get
\[ H^{p,p}(\comp\Sigma) \simeq \coker\bigl(\cubC^{p,p-1}(\comp\Sigma) \to \cubC^{p,p}(\comp\Sigma)\bigr) \simeq \coker\Bigl(\bigoplus_{\dims\tau=p-1}\SF^1(\infty_\tau) \to \bigoplus_{\dims\sigma=p} \SF^0(\infty_\sigma)\Bigr). \]
We have a natural isomorphism $\bigoplus_{\dims\sigma=p} \SF^0(\infty_\sigma) \simeq Z^p(\Sigma)$. Moreover, since $\Sigma$ is saturated at any face $\tau\in\Sigma_{p-1}$, $\bigoplus_{\dims\tau=p-1}\SF^1(\infty_\tau)$ is naturally isomorphic to the kernel of $Z^p(\Sigma) \to A^p(\Sigma)$ as described in the Localization Lemma \ref{lem:kernel_Z^k}. Hence $H^{p,p}(\comp\Sigma) \simeq A^p(\Sigma)$. Moreover, by Lemma \ref{lem:hypercube_vanishing}, the isomorphism is induced by the map
\[ \begin{tikzcd}[column sep=small, row sep=0em]
C^{p,p}(\comp\Sigma) \rar& \cubC^{p,p}(\comp\Sigma) \rar& Z^p(\Sigma) \\
\alpha \ar[rr, mapsto] && \sum_{\sigma\in\Sigma_p}\alpha(\nu_\sigma)\x_\sigma.
\end{tikzcd} \]

It remains to prove that the above isomorphism respects the products. This can be done by an explicit calculation of the inverse of the map $\Psi$ in the statement of \ref{thm:ring_morphism}, given in the next section, which shows that the cup-product in cohomology corresponds to product in Chow ring. The theorem thus follows. \qed

\subsection{Explicit description of the inverse $\Psi^{-1} \colon A^p(\Sigma) \to H^{p,p}(\comp \Sigma)$} \label{subsec:Ap_to_Hpp}

In this section, we provide an explicit description of the inverse of $\Psi$.

\smallskip
Before going through the description, let us introduce some new notations which will help in working with the cohomology of $\comp\Sigma$. In the following, if $a$ is a cochain in $C^{p,q}(\comp\Sigma)$ with $q=\dims\sigma-\dims\tau$, we denote by $a^\tau_\sigma \in \SF^p(\cube^\tau_\sigma)$ the part of $a$ which lives on the face $\cube^\tau_\sigma$. If $\tau=\conezero$, we just write $a_\sigma$. Moreover, if $\delta$ is any face of dimension $q$ in $\comp\Sigma$ and $\alpha \in \SF^p(\delta)$, then we denote by $(\delta, \alpha)$ the cochain in $C^{p,q}(\comp\Sigma)$ whose restriction to $\delta$ is $\alpha$ and which vanishes everywhere else.

\subsubsection{The inverse in cohomology degree one} \label{sec:inverse_degree_one} We consider first the part $A^1(\Sigma)$. Let $\rho \in \Sigma_1$. We give an element of $H^{1,1}(\comp\Sigma)$ whose image by $\Psi$ coincides with $\x_\rho$. For this, we first define an element $a \in C^{1,1}(\comp\Sigma)$ as follows. For any cone $\sigma\sim \rho$ in $\Sigma$, we choose an element $\alpha_\sigma \in \SF^1(\cube^{\sigma}_{\sigma \vee \rho})$ which takes value one on the vector $\e^{\sigma}_{\sigma \vee \rho}$. We then set
\[a := \sum_{\sigma \sim \rho} (\cube^{\sigma}_{\sigma \vee \rho}, \alpha_\sigma).\]

Note that the only part of $a$ which has sedentarity $\conezero$ is $(\cube_{\rho}, \alpha_\rho)$, and we have $\Psi(a) =\alpha_\conezero(\e_\rho)\x_\rho = \x_\rho$. So the statement would have followed if $a$ were a cocycle. The idea now is to find an element $b$ in $C^{1,1}(\comp\Sigma)$ supported only on faces of non-zero sedentarity such that $a-b$ is a cocycle: in that case, since $\Psi(b)=0$, we still get $\Psi(a-b) =\x_\rho$ and the inverse image of $\x_\rho$ will be represented by the class of $a-b$.

\smallskip
Recall that $\e^{\sigma}_{\sigma \vee \rho}$ is the projection of $\e_\rho$ at infinity into $N^\sigma_\infty$. We now describe $\~a := \d a$. Consider a pair of faces $\tau \subface \eta$ in $\Sigma$ with $\dims\eta-\dims\tau = 2$. Then, we have $\~a^\tau_\eta =0$ if either $\rho \not\subface \eta$ or $\rho \subface \tau$. In the remaining cases, there must exist a ray $\rho'$ such that $\eta = \tau \vee \rho \vee \rho'$, and in this case, we get
\[ \~a^\tau_\eta = \pm(a^\tau_{\tau\vee\rho} - \pi^*(a^{\tau\vee\rho'}_{\tau\vee\rho'\vee\rho})) = \pm(\alpha_\tau - \pi^*(\alpha_{\tau\vee\rho'})), \]
where $\pi$ is the projection $N^\tau \to N^{\tau\vee\rho'}$. In particular, we notice that $\~a^\tau_\eta$ is zero on the projection of $\e_\rho$.

From the above discussion, we obtain a well-defined pushforward $\pi^\rho_*(\~a^\tau_\eta) \in \SF^1(\cube^{\tau\vee\rho}_\eta)$ where $\pi^\rho\colon N^\tau_\infty \to N^{\tau\vee\rho}_\infty$ is the natural projection.

\smallskip
Set
\[ b := \sum_{\tau,\eta} \sign(\cube^{\tau\vee\rho}_\eta, \cube^\tau_\eta) \Bigl(\cube^{\tau\vee\rho}_\eta,\pi^{\rho}_{*}(\~a^\tau_\eta)\Bigr) \]
where the sum is restricted to pairs $\tau\subface\eta$ for which the term is well-defined, that is those with $\rho \not \subface \tau$, $\rho \subface \eta$ and $\dims\eta-\dims\tau=2$. The element $b$ has been defined in order to get $\d b = \d a$ on every face which is not in $\comp\Sigma^\rho_\infty$. In particular, the coboundary $c := \d(a-b)$ of $a-b$ has support in $\comp\Sigma^\rho_\infty$. Since $\d c=0$, the positive-sedentarity-vanishing lemma~\ref{lem:psed_vanishing} below implies that $c =0$. This shows that $a-b$ is a cocycle. Since $b$ has support only on the faces of non-zero sedentarity, we get in addition $\Psi(a-b) = \x_\rho$, as required.

\begin{lemma}[Positive-sedentarity-vanishing lemma] \label{lem:psed_vanishing} Let $\Sigma$ be a unimodular fan. Let $\zeta$ be a non-zero face of\/ $\Sigma$ and let $c \in C^{p,q}(\comp \Sigma)$ be a cochain supported in $\Sigma^\zeta_\infty$. Assume that $\d c =0$. Then we have $c=0$.
\end{lemma}

\begin{proof} Actually this is a consequence of the vanishing of $\EpI{p}_\infty^{a,b}$ for $b<0$. Let us give a more elementary proof here. Since $c$ is supported in $\Sigma^\zeta_\infty$, we need to show the vanishing of $c^\tau_\sigma$ for $\tau, \sigma$ in $\Sigma$ with $\dims\sigma - \dims\tau =q$ and $\zeta\subface \tau$.

Let $\rho$ be a ray of $\zeta$ and set $\eta \ssubface \tau$ so that we get $\tau = \rho \vee \eta$. Consider the face $\delta = \cube^\eta_\sigma$ of $\comp \Sigma$ and note that $\delta$ is (a hypercube) of dimension $q+1$. The only face of dimension $q$ in $\delta$ which lies in $\comp\Sigma^\zeta_\infty$ is $\cube^\tau_\sigma$. Since $c$ has support in $\comp\Sigma^\zeta_\infty$, it follows that the components of $c$ on faces of $\delta$ different from $\cube^\tau_\sigma$ are all zero. Denoting by $\pi$ the projection $N^\eta_\infty \to N^\tau_\infty$, which induces a surjection $\pi_*\colon \SF_p(\cube^\eta_\sigma) \to \SF_p(\cube^\tau_\sigma)$ and in injection $\pi^*\colon \SF^p(\cube^\tau_\sigma)\to \SF^p(\cube^\eta_\sigma)$, we thus get
\[ \pi^*(c^\tau_\sigma) = \pm(\d c)^\eta_\sigma =0.\]
This implies that $c^\tau_\sigma =0$, and the lemma follows. \end{proof}

\subsubsection{Cup-product in cubical complexes}

To extend the description of the inverse of $\Psi$ to higher degrees, we need to recall the formula for the cup product in cubical complexes. Let $a \in C^{p,q}(\comp\Sigma)$ and $b \in C^{p',q'}(\comp\Sigma)$. Then for any pair of faces $\tau \subface \eta$ with $\dims\eta - \dims\tau = q+q'$,
\[ (a \smile b)^\tau_\eta = \sum_{\tau\subface\sigma\subface\eta \\ \dims\sigma-\dims\tau = q} \omega^\tau_\eta(\nu^\tau_\sigma \wedge \pi^*(\nu^\sigma_\eta))\ \cdot\ a^\tau_\sigma \wedge \pi^*(b^\sigma_\eta) \]
where, for each $\tau \subface \sigma$, $\pi\colon N^\tau \to N^\sigma$ is the projection. As usual, this cup product induces a cup product on cohomology $\smile\colon H^{p,q}(\comp\Sigma) \times H^{p',q'}(\comp\Sigma) \to H^{p+p',q+q'}(\comp\Sigma)$.

\subsubsection{The inverse in higher cohomological degrees}

{
\renewcommand{\a}[1]{{}_{#1}a}%
Let $\sigma \in \Sigma_p$. By Localization Lemma~\ref{lem:kernel_Z^k}, it suffices to find a preimage of $\x_\sigma$. Let $\rho_1, \dots, \rho_p$ be the rays of $\sigma$. For $i \in \{1,\dots,p\}$, we denote by $\a{i}$ the preimage of $\x_{\rho_i}$ as defined in Section \ref{sec:inverse_degree_one}. In particular, we note that $\a{i}$ is supported by faces of the form $\cube^\tau_\sigma$ with $\rho_i\subface\sigma$. We claim that the element
\[ a := \a1 \smile \a2 \smile \dots \smile \a{p} \]
is a preimage of $\x_\sigma$. We compute $a$ as follows. Consider a face $\sigma' \in \Sigma_p$ and denote by $\rho'_1, \dots, \rho'_p \in \Sigma_1$ the rays of $\sigma'$. In the following, for a permutation $s \in \mathfrak S_p$ of $[p]$ and $k\in [p]$, we set
\[ \sigma'(s, k) := \rho'_{s(1)} \vee \dots \vee \rho'_{s(k)}, \]
the face of $\sigma'$ with rays $\rho'_{s(1)}, \dots, \rho'_{s(k)}$.

\smallskip
Expanding the cup product using the formula stated in the previous section, we find
\[ a_{\sigma'} = \sum_{s \in \mathfrak S_p} \omega_{\sigma'}(\nu_{\rho'_{s(1)}} \wedge \dots \wedge \nu_{\rho'_{s(p)}}) \ \cdot\ \a1^{\sigma'(s,0)}_{\sigma'(s,1)} \wedge \a2^{\sigma'(s,1)}_{\sigma'(s,2)} \wedge \dots \wedge \a{p}^{\sigma'(s,p-1)}_{\sigma'(s,p)}, \]
where for the ease of reading, we omit to precise the pullback by different projections. Each term in the above sum is nontrivial only if $\rho_1 \subface \sigma'(s,1)$, $\rho_2 \subface \sigma'(s,2)$, \dots, and $\rho_p \subface \sigma'(s,p)$, \ie, if and only if $\sigma'=\sigma$ and the permutation $s$ is identity.

Since $\pi_{\sigma'(\id,j)}(\e_{\rho_i}) = 0$ for $j \geq i$, we get that
\begin{align*}
a_\sigma(\nu_\sigma)
  &= \omega_\sigma(\e_{\rho_1}\wedge\dots\wedge\e_{\rho_p}) \cdot a_\sigma(\e_{\rho_1}\wedge\dots\wedge\e_{\rho_p}) \\
  &= \bigl(\omega_\sigma(\e_{\rho_1}\wedge\dots\wedge\e_{\rho_p})\bigr)^2 \bigl(\a1^{\sigma(\id,0)}_{\sigma(\id,1)} \wedge \a2^{\sigma(\id,1)}_{\sigma(\id,2)} \wedge \dots \wedge \a{p}^{\sigma(\id,p-1)}_{\sigma(\id,p)}\bigr)\bigl(\e_{\rho_1}\wedge\dots\wedge\e_{\rho_p}\bigr) \\
  &= \a1^{\sigma(\id,0)}_{\sigma(\id,0)\vee\rho_1}(\e_{\rho_1})\ \cdots\ \a{p}^{\sigma(\id,p-1)}_{\sigma(\id,p-1)\vee\rho_p}(\e_{\rho_p}) \\
  &= 1.
\end{align*}
We thus infer that $a$ is a preimage of $\x_\sigma$. This achieves the description of the inverse of the map $\Psi$.
}

\section{Tropical Deligne resolution} \label{sec:deligne}

Let $\Sigma$ be a smooth unimodular tropical fan. We follow the notations of previous sections. In particular, for each cone $\sigma$, we denote by $\Sigma^\sigma$ the star fan of $\sigma$ in $\Sigma$, which is a smooth unimodular tropical fan in $N^\sigma_\R$, and by $\comp \Sigma^\sigma$, we denote its canonical compactification. For the canonical compactification of a tropical fan $\Sigma$, recall as well that we set
\[H^k(\comp \Sigma) := \bigoplus_{p+q=k} H^{p,q}(\comp \Sigma).\]
In the case of a smooth unimodular tropical fan $\Sigma$, by Theorem~\ref{thm:Hodge_conjecture}, $H^k(\comp \Sigma)$ is non-vanishing only in even degrees and in that case it is equal to $H^{k/2,k/2}(\comp \Sigma)$.

\smallskip
The aim of this section is to prove the following theorem.

\begin{thm}[Tropical Deligne resolution]\label{thm:deligne} Notations as above, we have the following long exact sequence
\[0 \rightarrow \SF^p(\conezero) \rightarrow \bigoplus_{\sigma \in \Sigma_p} H^0(\comp \Sigma_\infty^\sigma) \rightarrow \bigoplus_{\sigma \in \Sigma_{p-1}} H^2(\comp \Sigma_\infty^\sigma) \rightarrow \dots \rightarrow \bigoplus_{\sigma \in \Sigma_1} H^{2p-2}(\comp \Sigma_\infty^\sigma) \rightarrow H^{2p} (\comp \Sigma) \to 0. \]
\end{thm}

The maps between cohomology groups in the above sequence are given by the sum of Gysin maps. Namely, for a cone $\sigma \in \Sigma$, using the notations of Section~\ref{sec:prel}, the fan $\Sigma_\infty^\sigma$ is based at the point $\infty_\sigma$ of $\comp \Sigma$, is naturally isomorphic to $\Sigma^\sigma$, and has closure $\comp{\Sigma}_\infty^\sigma$ in $\comp\Sigma$ which can be identified with the canonical compactification $\comp \Sigma^\sigma$ of $\Sigma^\sigma$. This means the canonical compactification $\comp\Sigma_\infty^\sigma$ naturally lives in $\comp \Sigma$. Moreover, for an inclusion of cones $\tau \ssubface \sigma$, we get an inclusion $ \i_{\sigma \ssupface \tau}\colon \comp \Sigma_\infty^\sigma \hookrightarrow \comp \Sigma_\infty^\tau$. This induces a map on cohomology $\i^*_{\tau \ssubface \sigma} \colon H^\bul(\comp \Sigma_\infty^\sigma) \to H^\bul(\comp \Sigma_\infty^\tau)$, which, by applying the Poincaré duality for the smooth tropical varieties $\comp \Sigma_\infty^\sigma$ and $\comp \Sigma_\infty^\tau$, leads to the \emph{Gysin map} $\gys_{\sigma \ssupface \tau}\colon H^\bul(\comp \Sigma_\infty^\sigma) \to H^{\bul+2}(\comp \Sigma_\infty^\tau)$.

In what follows, in order to simplify the notations, sometimes we drop $\infty$ from indices and identify $\comp \Sigma^\sigma$ with $\comp{\Sigma}_\infty^\sigma$ living naturally in $\comp \Sigma$ as the closure of those points which have sedentarity equal to $\sigma$.

\smallskip
Apart from the last subsection, where we study what happens is the non-smooth case, the rest of this section is devoted to the proof of this theorem. As it becomes clear from the above discussion, using the Poincaré duality for canonical compactifications $\comp \Sigma^\sigma$, it will be enough to prove the exactness of the following complex for each $k$ (here $k=d-p$):

\begin{equation} \label{eqn:Deligne_dual}
0 \rightarrow H^{2k}(\comp \Sigma) \to \bigoplus_{\sigma \in \Sigma_1} H^{2k}(\comp\Sigma^\sigma) \to \bigoplus_{\sigma \in \Sigma_2} H^{2k}(\comp\Sigma^\sigma) \to \dots \to \bigoplus_{\sigma \in \Sigma_{d-k}} H^{2k}(\comp \Sigma^\sigma) \to \SF_{d-k}(\conezero) \to 0.
\end{equation}

\subsection{The sheaf of tropical holomorphic forms} Let $Z$ be a tropical variety with an extended polyhedral structure. We will be only interested here in $Z = \Sigma$ or $\comp \Sigma$, for a tropical fan $\Sigma$, so there is no harm in assuming this in what follows in this section.

For the tropical variety $Z$, we denote by $\Omega^k_Z$ the sheaf of \emph{tropical holomorphic $k$-forms on $Z$}. With real coefficients, this can be defined as the kernel of the second differential operator from Dolbeault $(k,0)$-forms to Dolbeault $(k,1)$-forms on $Z$, as in~\cites{JSS19, CLD, GK17}. We will give an alternative characterization of this sheaf which shows that it is actually defined over $\Z$ as the sheafification of the combinatorial sheaf $\SF^k$. As we are going to use some results from~\cite{JSS19}, it might be helpful to note that this sheaf is denoted $\mathcal L^k_Z$ in \emph{loc.\ cit.}

Let $Z$ be a tropical variety with an extended polyhedral structure. If $U$ is an open subset of $Z$, we say that $U$ is \emph{nice} if either $U$ is empty, or there exists a face $\gamma$ of $Z$ intersecting $U$ such that, for each face $\delta$ of $Z$, every connected components of $U\cap\delta$ contains $U\cap\gamma$. (Compare with the \emph{basic open sets} from~\cite{JSS19}.) This implies that $U$ is connected, and that for each face $\delta\in Z$ which intersects $U$, $\gamma$ is a face of $\delta$ and the intersection $\delta\cap U$ is connected. We call $\gamma$ the \emph{minimum face of\/ $U$}. We have the following elementary result.
\begin{prop} Nice open sets form a basis of open sets on $Z$.
\end{prop}

The sheaf $\Omega^k_Z$ is the unique sheaf on $Z$ such that for each nice open set $U$ of $Z$ with minimum face $\gamma$, we have
\[ \Omega^k_Z(U)=\SF^k(\gamma). \]

\medskip

Suppose now $\comp Z$ is a compactification of $Z$ and denote by $\i\colon Z\hookrightarrow\comp Z$ the inclusion. We denote by $\Omega^k_{Z,c}$ the \emph{sheaf of holomorphic $k$-forms on $\comp Z$ with compact support in $Z$} defined on connected open sets by
\[ \Omega^k_{Z,c}(U):=\begin{cases}
\Omega^k_Z(U) & \text{if $U\subseteq Z$,} \\
0             & \text{otherwise.}
\end{cases} \]

In other words, we have
\[ \Omega^k_{Z,c}=\i_!\Omega^k_Z \]
that is the \emph{direct image with compact support} of the sheaf $\Omega^k_Z$. In particular, in the case we study here, the cohomology with compact support of $\Omega^k_Z$ is computed by the usual sheaf cohomology of $\Omega^k_{Z,c}$, \ie, we have
\[ H_c^\bul(Z, \Omega^k_Z)=H^\bul(\comp Z, \Omega^k_{Z,c}). \qedhere \]

\subsection{Cohomology with coefficients in the sheaf of tropical holomorphic forms}
We now specify the above set-up for $Z = \Sigma$ and $\comp \Sigma$ the canonical compactification of $\Sigma$.

For each cone $\sigma \in \Sigma$, we get the sheaf $\Omega^k_{\comp \Sigma^\sigma_\infty}$ of holomorphic $k$-forms on $\comp \Sigma^\sigma_\infty \hookrightarrow \comp \Sigma$ which by extension by zero leads to a sheaf on $\comp \Sigma$. We denote this sheaf by $\Omega^k_\sigma$. The following proposition describes the cohomology of these sheaves.

\begin{prop} \label{prop:cohomology_Omega_fan} Notations as above, for each pair of non-negative integers $m, k$, we have
\[ H^m(\comp \Sigma, \Omega^k_\sigma) = \begin{cases}
  H^{2k}(\comp\Sigma^\sigma) = H^{k,k}(\comp\Sigma^\sigma) & \textrm{if $m = k$,} \\
  0  & \textrm{otherwise.}
\end{cases}\]
\end{prop}
\begin{proof}
We have
\[H^m(\comp \Sigma, \Omega^k_\sigma) \simeq H^m(\comp \Sigma^\sigma_\infty, \Omega^k_{\comp \Sigma^\sigma_\infty}) \simeq H^m(\comp \Sigma^\sigma, \Omega^k_{\comp \Sigma^\sigma}) \simeq H^{k,m}(\comp\Sigma^\sigma).\]
The result now follows from our Theorem~\ref{thm:Hodge_conjecture}.
\end{proof}

\subsection{The resolution $\Omega^k_\bul$ of the sheaf $\Omega^k_{\Sigma,c}$}
For a pair of faces $\tau\subface\sigma$ in $\Sigma$, from the inclusion maps $\comp \Sigma^\sigma_\infty \hookrightarrow \comp\Sigma^\tau_\infty \hookrightarrow \comp \Sigma$, we get natural restriction maps of sheaves $\i^*_{\tau\subface\sigma}\colon\Omega^k_\tau \to\Omega^k_\sigma$ on $\comp \Sigma$. Here as before, the map $\i_{\tau\subface\sigma} = \i_{\sigma \supface\tau}$ denotes the inclusion $\comp\Sigma^\sigma_\infty\hookrightarrow\comp\Sigma^\tau_\infty$.

\smallskip
We consider now the following complex of sheaves on $\comp \Sigma$
\begin{equation}
\Omega^k_\bul\colon \qquad \Omega_\conezero^k \to \bigoplus_{\varrho \in \Sigma_1}\Omega^k_{\varrho} \to \bigoplus_{\sigma \in \Sigma_2} \Omega^k_{\sigma} \to \dots \to \bigoplus_{\sigma \in \Sigma_{d-k}} \Omega^k_{\sigma}
\end{equation}
concentrated in degrees $0, 1, \dots, d-k$, given by the dimension of the cones $\sigma$ in $\Sigma$, whose boundary maps are given by
\[ \alpha\in\Omega_\sigma^k \mapsto \d\alpha := \sum_{\zeta\ssupface\sigma}\sign(\sigma,\zeta)\i^*_{\sigma\ssubface\zeta}(\alpha). \]
We will derive Theorem~\ref{thm:deligne} by looking at the hypercohomology groups $\hyp^\bul(\comp\Sigma, \Omega^k_\bul)$ of this complex and by using the following proposition.

\begin{prop} \label{prop:exactness_Omega}
The following sequence of sheaves is exact
\[ 0 \to \Omega_{\Sigma, c}^k \to \Omega_\conezero^k \to \bigoplus_{\varrho \in \Sigma_1}\Omega^k_\varrho \to \bigoplus_{\sigma \in \Sigma_2}\Omega^k_\sigma \to \dots \to \bigoplus_{\sigma \in \Sigma_{d-k}}\Omega^k_\sigma \to 0. \]
\end{prop}

\begin{proof} It will be enough to prove that taking sections over nice open sets $U$ give exact sequences of abelian groups.

If $U$ is included in $\Sigma$, clearly by definition we have $\Omega^k_{\Sigma,c}(U)\simeq\Omega^k_{\conezero}(U)$, and the other sheaves of the sequence have no nontrivial section over $U$. Thus, the sequence is exact over $U$.

It remains to prove that, for every nice open set $U$ having non-empty intersection with $\comp\Sigma\setminus\Sigma$, the sequence
\[ 0 \to \Omega_\conezero^k(U) \to \bigoplus_{\varrho \in \Sigma_1}\Omega^k_\varrho(U) \to \bigoplus_{\sigma \in \Sigma_2}\Omega^k_\sigma(U) \to \dots \to \bigoplus_{\sigma \in \Sigma_{d-k}}\Omega^k_\sigma(U) \to 0 \]
is exact. Let $\gamma\in\comp\Sigma$ be the minimum face of $U$. Let $\sigma\in\Sigma$ be the sedentarity of $\gamma$. The closed strata (for the fan stratification) of $\comp \Sigma$ which intersect $U$ are exactly those of the form $\comp\Sigma_\infty^\tau \simeq \comp\Sigma^{\tau}$ with $\tau\subface\sigma$. Moreover, if $\tau$ is a face of $\sigma$, we have
\[ \Omega^k_\tau(U)=\SF^k(\gamma). \]
Thus, the previous sequence can be rewritten in the form
\[ 0 \to \SF^k(\gamma) \to \bigoplus_{\tau\subface\sigma \\ \dims{\tau}=1}\SF^k(\gamma) \to \bigoplus_{\tau\subface\sigma \\ \dims{\tau}=2}\SF^k(\gamma) \to \dots \to \bigoplus_{\tau\subface\sigma \\ \dims{\tau}=\dims{\sigma}}\SF^k(\gamma) \to 0. \]
This is just the cochain complex of the simplicial cohomology (for the natural simplicial structure induced by the faces) of the cone $\sigma$ with coefficients in the group $\SF^k(\gamma)$. This itself corresponds to the reduced simplicial cohomology of a simplex shifted by $1$. This last cohomology is trivial, thus the sequence is exact. That concludes the proof of the proposition.
\end{proof}

\begin{proof}[Proof of Theorem~\ref{thm:deligne}]
We have
\[ H^m(\comp\Sigma,\Omega^k_{\Sigma,c}) = H^m_c(\Sigma,\Omega^k_\Sigma) = H^{k,m}_{c}(\Sigma)=
\begin{cases}
  \SF^{d-k}(\conezero)^\dual = \SF_{d-k}(\conezero) & \text{if $m=d$,}\\
  0 & \text{otherwise.}
\end{cases} \]
By proposition \ref{prop:exactness_Omega}, the cohomology of $\Omega^k_{\Sigma,c}$ becomes isomorphic to the hypercohomology $\hyp(\Sigma, \Omega^k_\bul)$. Thus, we get
\[ \hyp(\Sigma, \Omega^k_\bul)\simeq \SF_{d-k}(\conezero)[-d], \]
meaning
\[\hyp^m(\comp\Sigma, \Omega^k_\bul) = \begin{cases}
  \SF_{d-k}(\conezero) & \textrm{for $m =d$,} \\
  0 & \textrm{otherwise.}
\end{cases}\]

On the other hand, using the hypercohomology spectral sequence, combined with Proposition \ref{prop:cohomology_Omega_fan}, we infer that the hypercohomology of $\Omega^k_\bul$ is given by the cohomology of the following complex:
\[ 0 \rightarrow H^{2k}(\comp\Sigma)[-k] \to \bigoplus_{\sigma \in \Sigma_1} H^{2k}(\comp\Sigma^\sigma)[-k-1] \to \dots \to \bigoplus_{\sigma \in \Sigma_{d-k}} H^{2k}(\comp\Sigma^\sigma)[-d] \to 0. \]
We thus conclude the exactness of the sequence
\begin{equation}
0 \rightarrow H^{2k}(\comp \Sigma) \to \bigoplus_{\sigma \in \Sigma_1} H^{2k}(\comp\Sigma^\sigma) \to \bigoplus_{\sigma \in \Sigma_2} H^{2k}(\comp\Sigma^\sigma) \to \dots \to \bigoplus_{\sigma \in \Sigma_{d-k}} H^{2k}(\comp\Sigma^\sigma) \to \SF_{d-k}(\conezero) \to 0,
\end{equation}
and the theorem follows.
\end{proof}

\subsection{The non-smooth case} \label{subsec:Deligne_non_smooth}

In the case where $\Sigma$ is not necessarily smooth, we still have the following proposition that we will use in Section \ref{sec:homology_tropical_modification}.

\begin{prop} \label{prop:non_smooth_Deligne}
Let $\Sigma$ be any unimodular tropical fan and let $k$ be an integer. Then the following sequence is exact
\[ \bigoplus_{\sigma \in \Sigma_{d-k-1}}H^{k,k}(\comp\Sigma^\sigma) \longrightarrow \bigoplus_{\sigma \in \Sigma_{d-k}} H^{k,k}(\comp\Sigma^\sigma) \longrightarrow H_c^{k,d}(\Sigma) \longrightarrow 0. \]
\end{prop}

In the smooth case, one can recognize the end of the long exact sequence \eqref{eqn:Deligne_dual}: indeed, when $\Sigma$ is smooth, we have $H_c^{k,d}(\Sigma) \simeq \SF_{d-k}(\conezero)$.

\begin{proof}
Proposition \ref{prop:exactness_Omega} still holds. However, we only get that
\[ \hyp^\bul(\Sigma,\Omega_\bul^k) \simeq H_c^{k,\bul}(\Sigma), \]
with no further simplification.

If we denote by $\E_\bul^{\bul,\bul}$ the hypercohomology spectral sequence, we get
\[ \E_1^{a,b} = \bigoplus_{\sigma\in\Sigma_a} H^{k,b}(\comp\Sigma^\sigma). \]
By Theorem \ref{thm:ring_morphism} and Remark \ref{rem:ring_morphism_non_saturated}, we know that for any $\sigma\in\Sigma$, $H^{k,b}(\comp\Sigma^\sigma)$ is trivial if $b > k$. Moreover, $H^{k,b}(\comp\Sigma^\sigma)$ is clearly trivial if $\dims\sigma > d-k$. Hence $\E_1^{a,b}$ is nontrivial only for $0 \leq a \leq d-k$ and $0 \leq b \leq k$. Computing the further pages, we get that
\[ \E_\infty^{d-k,k} = \E_2^{d-k,k} = \coker\bigl(\E_1^{d-k-1,k} \to \E_1^{d-k,k}\bigr). \]
Since $\E_\infty^{d-k,k}$ is the only nontrivial term of degree $d$, we get
\[ H_c^{k,d}(\Sigma) \simeq \E_\infty^{d-k,k} = \coker\Bigl(\bigoplus_{\sigma \in \Sigma_{d-k-1}}H^{k,k}(\comp\Sigma^\sigma) \longrightarrow \bigoplus_{\sigma \in \Sigma_{d-k}} H^{k,k}(\comp\Sigma^\sigma)\Bigr), \]
which concludes the proof.
\end{proof}

\section{Compactifications of complements of hyperplane arrangements}\label{sec:comp}

The aim of this section is to discuss an application of the results we proved previously in the study of cohomology rings of wonderful compactifications.

\subsection{Tropical compactifications of subvarieties of tori} Let $N$ be a lattice of finite rank $n$ and denote by $M$ the dual lattice. Let $\torus = \spec(\C[M])$ be the corresponding complex torus of dimension $n$.

Let $X \hookrightarrow \torus$ be a connected closed subvariety of $\torus$ and let $\Trop(X) \subseteq N_\R$ be the tropicalization of $X$, which comes with a natural weight function. Any fan with support in the tropicalization $\Trop(X)$ equipped with this weight function verifies the balancing condition~\cites{BIMS, MR09, MS15}. In particular, if the weight function is equal to one, then $\Trop(X)$ is a tropical fan. That is the assumption we make in this section.

We are interested in compactifications of $X$ obtained by taking the closure in a toric variety which contains $\torus$ as open torus. The following theorem summarizes several known properties of these compactifications.

\begin{thm} Notations as above, let $X \hookrightarrow \torus$. Consider a unimodular fan $\Sigma$ in $N_\R$ and denote by $\CP_\Sigma$ the corresponding toric variety. Let $\comp X$ be the closure of $X$ in $\CP_\Sigma$.
\begin{itemize}
\item \emph{(\cites{Tev07, KKMS})} $\comp X$ is proper if and only if\/ $\supp\Sigma \supseteq \Trop(X)$.

\item \emph{(\cites{ST08, Hacking08})} The following conditions are equivalent
\begin{enumerate}
\item $\supp\Sigma = \Trop(X)$.
\item for any torus orbit $\torus^\sigma$, $\sigma\in \Sigma$, the intersection $\comp X \cap \torus^\sigma$ has pure dimension equal to $m-\dims{\sigma}$.
\end{enumerate}

\item \emph{(\cites{Tev07, Hacking08})} Assume $\supp\Sigma =\Trop(X)$. The following conditions are then equivalent
\begin{enumerate}
\item for each $\sigma \in \Sigma$, the intersection $\comp X \cap \torus^\sigma$ is smooth.
\item the multiplication map $\mult\colon \torus \times \comp X \to \CP_\Sigma$ is smooth.
\end{enumerate}

\item \emph{(\cite{Tev07})} Assume that $\mult\colon \torus \times \comp X \to \CP_\Sigma$ is smooth. Then for any unimodular fan $\Sigma'$ with support $\Trop(X)$, denoting by $\comp X'$ the closure of $X$ in $\CP_{\Sigma'}$, the multiplication map $\mult'\colon \torus \times \comp X' \to \CP_{\Sigma'}$ is smooth.
\end{itemize}
\end{thm}

We have the following definitions introduced by Tevelev~\cite{Tev07}.

\begin{defi}[Tropical compactification]\label{defn-tropcomp}
Let $X$ be a subvariety of the algebraic torus $\torus$.
\begin{itemize}
\item Let $\Sigma$ be a unimodular fan with support $\Trop(X)$. The compactification $\overline X$ obtained by taking the closure of $X$ in $\CP_\Sigma$ is called \emph{tropical} if the multiplication map $\mult\colon \torus \times \comp X \to \CP_\Sigma$ is faithfully flat.

\item The subvariety $X$ of the torus $\torus$ is called \emph{schön} if a compactification $\overline X \subseteq \CP_{\Sigma}$ with $\Sigma$ unimodular with support $\supp\Sigma =\Trop(X)$ (and so any such compactification) has the property that the multiplication map $\mult\colon \torus \times \comp X \to \CP_\Sigma$ is smooth. \qedhere
\end{itemize}
\end{defi}

\begin{remark}\label{rem:schon}
In the case where $X$ is schön, it follows that the weights are all equal to one.
\end{remark}

Tevelev proved in~\cite{Tev07} that any subvariety $X$ of a torus $\torus$ has a tropical compactification $\comp X$. Moreover, a tropical compactification $\comp X$ has the following properties:
\begin{enumerate}
\item The \emph{boundary} $\comp X \setminus X$ is divisorial.
\end{enumerate}
Denote this divisor by $D$ and let $D_1, \dots, D_r$ be its irreducible components.
\begin{enumerate}[resume]
\item The boundary divisor has combinatorial normal crossings: for any collection of $r$ irreducible divisors $D_{i_1},\dots, D_{i_r}$ in $D$, the intersection $D_{i_1} \cap \dots \cap D_{i_r}$ is either empty or it has codimension $r$ in $\comp X$.
\end{enumerate}

It follows from the properties listed above that in the case $X$ is schön, the tropical compactification $\comp X$ is smooth and the boundary divisor $D$ is strict normal crossing (we assume the underlying fan $\Sigma$ is always unimodular, so the toric varieties we consider are all smooth).

\subsection{Compactifications of the complement of a hyperplane arrangement}
Consider now a collection of hyperplanes $H_0, H_1, \dots, H_m$ in $\CP^r$ given by linear forms $\ell_0, \ell_1, \dots, \ell_m$. We assume that the vector space generated by $\ell_0,\ell_1, \dots, \ell_m$ has maximum rank, \ie, the intersection $H_0\cap H_1\cap \dots \cap H_m$ is empty. Let $X$ be the complement $\CP^r \setminus \bigcup_{j=0}^m H_j$ and consider the embedding $X \hookrightarrow \torus$ for the torus $\torus = \spec(\C[M])$ with
\[M = \ker(\deg\colon \Z^{m+1}\to \Z), \quad \deg(x_0, \dots, x_m) = x_0+ \dots + x_m.\]
The embedding is given by coordinates $[\ell_0: \dots :\ell_m]$. Denote by $N$ the dual lattice to $M$.
\begin{thm}\label{thm:tevelev} Notations as above, we have the following.
\begin{itemize}
\item \emph{(Tevelev~\cite{Tev07}*{Theorem 1.5})} The subvariety $X \hookrightarrow \torus$ is schön.

\item \emph{(Ardila-Klivans~\cite{AK06})} The tropicalization of $X$ coincides with the support of the Bergman fan $\Sigma_\Ma$ of the matroid $\Ma$ associated to the hyperplane arrangement.
\end{itemize}
\end{thm}

Here is the main theorem of this section.

\begin{thm}\label{thm:generalized_FY} Let $\Sigma$ be a unimodular fan with support $\Trop(X)= \supp{\Sigma_\Ma}$ and denote by $\comp X$ the corresponding compactification. The cohomology ring of $\comp X$ is concentrated in even degrees and it is of Hodge-Tate type, \ie, the Hodge-decomposition in degree $2p$ is concentrated in bidegree $(p,p)$. Moreover, the restriction map
\[A^\bul(\Sigma) \simeq A^\bul(\CP_\Sigma) \longrightarrow H^{2\bul}(\comp X)\]
is an isomorphism of rings.
\end{thm}
Combined with Theorem~\ref{thm:Hodge_isomorphism}, this leads to an isomorphism
\[H^{2p}(\comp \Sigma) \simeq H^{2p}(\comp X).\]

The theorem provides a \emph{tropical} generalization of the results of Feichtner-Yuzvinsky \cite{FY04} and De Concini-Procesi \cite{DP95}, by going beyond the combinatorial data of building and nested sets, assumed usually in the theory of wonderful compactifications.

\begin{proof}[Proof of Theorem~\ref{thm:generalized_FY}]
We give the proof with rational coefficients, and sketch an alternate proof in Remark~\ref{rem:generalized_FY_integral} which works with integral coefficients.

Consider a cone $\sigma \in \Sigma$ and let $X^\sigma$ be the intersection of $\comp X$ with the torus orbit $\torus^\sigma$ in $\CP_\Sigma$ associated to $\sigma$. We have the following properties:
\begin{itemize}
\item $X^\sigma$ is itself isomorphic to the complement of a hyperplane arrangement. The embedding $X^\sigma \hookrightarrow \torus^\sigma$ has tropicalization given by $\Trop(X^\sigma) = \supp{\Sigma^\sigma_\infty} \hookrightarrow N^\sigma_\infty$.
\item The closure of $X^\sigma$ in $\comp X$ coincides with the compactification $\comp X^\sigma$ of $X^\sigma$ in the toric variety $\CP_{\Sigma^\sigma_\infty}$.
\end{itemize}
Proceeding by induction on the dimension of $X$, we can thus assume that the theorem holds for all the strata $X^\sigma$ for any cone $\sigma \neq \conezero$ in $\Sigma$.

\smallskip
Consider now the smooth compactification $X \hookrightarrow \comp X$ with strict normal crossing divisor, as stated by Tevelev's theorem. The Deligne spectral sequence, which describes the mixed Hodge structure on the cohomology of $X$ in terms of the cohomology of the closed strata in the stratification of $\comp X$ given by the intersection of components of the boundary divisor, gives the following exact sequence:
\[0 \rightarrow H^p(X) \rightarrow \bigoplus_{\sigma \in \Sigma \\
\dims{\sigma} =p} H^0(\comp X^\sigma) \rightarrow \bigoplus_{\sigma \in \Sigma \\
\dims{\sigma} =p-1} H^2(\comp X^\sigma) \rightarrow \dots \rightarrow \bigoplus_{\sigma \in \Sigma \\
\dims{\sigma} =1} H^{2p-2}(\comp X^\sigma) \rightarrow H^{2p} (\comp X) \to 0, \]
where the maps between cohomology groups are again given by the Gysin maps for the embeddings $\comp X^\sigma \hookrightarrow \comp X^\tau$ for pair of faces $\tau \ssubface\sigma$. This was proved by~\cite{IKMZ} and is a consequence of a theorem of Shapiro~\cite{Shap93} which states that the $p$-th cohomology of the complement of a hyperplane arrangement is pure of weight $2p$ concentrated in bidegree $(p,p)$.

In Section~\ref{sec:deligne} we proved a tropical analogue of the above exact sequence for any smooth tropical fan, the tropical Deligne resolution:

\[0 \rightarrow \SF^p(\conezero) \rightarrow \bigoplus_{\sigma \in \Sigma \\
\dims{\sigma} =p} H^0(\comp\Sigma_\infty^\sigma) \rightarrow \bigoplus_{\sigma \in \Sigma \\
\dims{\sigma} =p-1} H^2(\comp\Sigma_\infty^\sigma) \rightarrow \dots \rightarrow \bigoplus_{\sigma \in \Sigma \\
\dims{\sigma} =1} H^{2p-2}(\comp\Sigma_\infty^\sigma) \rightarrow H^{2p} (\comp\Sigma) \to 0, \]
where the maps between cohomology groups are given by the Gysin maps for the embeddings $\comp\Sigma_\infty^\sigma \hookrightarrow \comp \Sigma_\infty^\tau$ for pair of faces $\tau \ssubface\sigma$.

It was also proved by Zharkov~\cite{Zha13} that the coefficient group $\SF^p(\conezero)$ is naturally isomorphic to the cohomology group $H^p(X)$. Under this isomorphism and the isomorphism between $H^{2p} (\comp\Sigma^\sigma)$ and $H^{2p}(\comp X^\sigma)$, obtained by combining the isomorphism of Theorem~\ref{thm:Hodge_isomorphism} and the restriction maps $A^\bul(\Sigma^\sigma) \simeq H^{2\bul}(\CP_{\Sigma^\sigma}) \to H^{2\bul}(\comp X^\sigma)$, for all $\sigma \neq \conezero$, we get that the above two Deligne exact sequences are actually isomorphic up to the last piece on the right. This shows that the natural map
\[H^{2p}(\comp\Sigma) \to H^{2p}(\comp X),\]
given by the composition of maps $H^{2p} (\comp\Sigma) \simeq A^p(\Sigma) \simeq A^p(\CP_\Sigma) \to A^{p}(\comp X) \to H^{2p}(\comp X)$, must be an isomorphism as well, which proves the theorem in even cohomological degrees.

Again, proceeding by induction, the Deligne spectral sequence implies that the cohomology groups $H^{2p+1}(\comp X)$ are all vanishing. The same statement was proved for $H^{2p+1}(\comp \Sigma)$ in Theorem~\ref{thm:Hodge_isomorphism}, so the isomorphism $H^\bul(\comp \Sigma) \simeq H^{\bul}(\comp X)$ holds as well in odd cohomological degrees.
\end{proof}

\begin{remark}\label{rem:generalized_FY_integral}
We sketch an alternative argument which allows to get the result with integral coefficients.

First, the case where the hyperplanes form the uniform matroid, and the fan $\Sigma$ has the same support as the Bergman fan $\Sigma_\Ma$, is direct. In this case, the toric variety and the compactification coincide, and the result follows from the isomorphism between the cohomology ring of the toric variety and the Chow ring of the fan, combined with the Hodge isomorphism theorem.

Proceeding now by induction, we prove the result for the case of an arbitrary hyperplane arrangement with $\Sigma$ having the same support as $\Sigma_\Ma$, for $\Ma$ the matroid of the hyperplane arrangement. Combining Keel's lemma, the weak factorization theorem, and Theorem~\ref{thm:tevelev}, it will be enough to produce one specific fan with the same support as $\Sigma_\Ma$ for which the result holds. We take the linear form $\ell$ which defines one of the hyperplanes. We take the tropicalization $\Trop(X')$ of the complement $X'$ of all the other hyperplanes with respect to the embedding given by their linear forms, and take the tropicalization of $\ell$, viewed as a rational function defined on $X'$, to obtain a piecewise integral linear function $f$ on $\Trop(X')$. The tropical modification of $\Trop(X')$ along the divisor $\div(f)$ produces a fan $\Sigma$ with the same support as the Bergman fan $\Sigma_\Ma$. Let $\Sigma'$ be the fan structure on $\Trop(X')$ and $\comp X'$ the corresponding compactification. Let $\comp X$ be the compactification of $X$ with respect to $\Sigma$. Then we get an isomorphism $\comp X \simeq \comp X'$, given by the graph of $\ell$, and an isomorphism between the Chow rings of $\Sigma$ and $\Sigma'$, given by Theorem~\ref{thm:invariant_chow_tropical_modification}. Combining these isomorphisms and applying the hypothesis of the induction, we get the theorem for $\Sigma$, and thus for any unimodular fan with the same support as $\Sigma$.
\end{remark}

\section{Homology of tropical modifications and shellability of smoothness} \label{sec:homology_tropical_modification}
The aim of this section is to study the behavior of homology and cohomology groups under tropical modifications, and to prove shellability of smoothness for tropical fans.

\subsection{Computation of the homology and cohomology of a tropical modification}

This section is devoted to the proof of the following theorem.

\begin{thm}[Tropical modification formula] \label{thm:homology_tropical_modification}
Let $\Sigma$ be a smooth tropical fan. Let $f$ be a conewise integral linear function on $\Sigma$ such that $\div(f)$ is reduced. Set $\Delta = \div(f)$. Then the following equalities between different homology and cohomology groups hold.
\begin{itemize}
\item \textup{(Open tropical modifications)} We have
\[ H^{\bul,\bul}_c(\tropmod{f}{\Sigma}) \simeq H^{\bul,\bul}_c(\Sigma,\Delta) \qquad \textrm{and} \qquad
H_{\bul,\bul}^\BM(\tropmod{f}{\Sigma}) \simeq H_{\bul,\bul}^\BM(\Sigma,\Delta), \]
where $H^{\bul,\bul}_c(\Sigma,\Delta)$ and $H^\BM_{\bul,\bul}(\Sigma,\Delta)$ denote the relative cohomology with compact support and the relative Borel-Moore homology respectively.

\item \textup{(Closed tropical modifications)} We have
\[ H^{\bul,\bul}_c(\ctropmod{f}{\Sigma}) \simeq H^{\bul,\bul}_c(\Sigma)\qquad \textrm{and} \qquad H_{\bul,\bul}^\BM(\ctropmod{f}{\Sigma}) \simeq H_{\bul,\bul}^\BM(\Sigma), \]
and these isomorphisms are compatible with the Poincaré duality.

\item \textup{(Open extended tropical modifications)} We have
\[ H^{\bul,\bul}_c(\tropmod{f}{\comp\Sigma}) \simeq H^{\bul,\bul}_c(\comp\Sigma,\comp\Delta) \qquad \textrm{and} \qquad H_{\bul,\bul}^\BM(\tropmod{f}{\comp\Sigma}) \simeq H_{\bul,\bul}^\BM(\comp\Sigma,\comp\Delta). \]

\item \textup{(Closed extended tropical modifications)}
We have
\[ H^{\bul,\bul}_c(\ctropmod{f}{\comp\Sigma}) \simeq H^{\bul,\bul}_c(\comp\Sigma) \qquad \textrm{and} \qquad H_{\bul,\bul}^\BM(\ctropmod{f}{\comp\Sigma}) \simeq H_{\bul,\bul}^\BM(\comp\Sigma), \]
and these isomorphisms are compatible with the Poincaré duality.
\end{itemize}
\end{thm}

\begin{remark}
Let us describe briefly the different isomorphisms. For non-extended tropical modifications, it is simply induced by the maps $\prtm_*\colon \SF_\bul(\basetm\sigma) \to \SF^\Sigma_\bul(\sigma)$ for $\sigma\in\Sigma$. For the extended tropical modification however, one has to add the maps induced by the inclusions $\i\colon \SF_\bul(\inftm\delta) \simeq \SF^\Delta_\bul(\delta) \hookrightarrow \SF^\Sigma_\bul(\delta)$ for $\delta\in\Delta$.
\end{remark}

Notice that the notation $\SF_\bul(\delta)$ is ambiguous for $\delta\in\Delta$. That is why, here and in the rest of this section, we precise as a superscript in which fan we consider $\SF_\bul(\delta)$: either $\SF^\Delta_\bul(\delta)$ or $\SF^\Sigma_\bul(\delta)$.

\smallskip
First note that the theorem is trivial in the case of a degenerate tropical modification, \ie, if $\Delta=0$. Indeed, since $\Sigma$ is smooth, it is div-faithful by Theorem \ref{thm:smooth_principal_div-faithful}, and then $\tropmod{f}{\Sigma}$ is isomorphic to $\Sigma$ (see Proposition \ref{prop:star_fans_tropical_modifications}) and the theorem follows easily.

In what follows, we assume the tropical modification is non-degenerate. As usual, we use the notations of Section \ref{subsec:tropical_modification}. In particular, recall that $\prtm$ is the projection associated to the tropical modification, and that $\etm$ is the unit vector of the special new ray above $\conezero$ in the tropical modification.

The key point in the proof of the theorem above is the following lemma which provides a description of the coefficient sheaves in the tropical modification.

\begin{lemma}[Local tropical modification formula] \label{lem:F_and_tropical_modification}  Let $\Sigma$ be a smooth tropical fan, and let $\Delta =\div(f)$ be a reduced divisor associated to a conewise integral linear function on $\Sigma$. Let $\~\Sigma = \tropmod{\Delta}{\Sigma}$.

Let $\sigma$ be a face in $\Sigma \setminus \Delta$, $\delta$ be a face in $ \Delta$, and let $p$ be a non-negative integer. We have the following short exact sequences.
\[ \begin{tikzcd}[row sep = tiny]
& 0 \rar& \SF^{\~\Sigma}_p(\basetm\sigma) \rar{\prtm_*}& \SF^\Sigma_p(\sigma) \rar& 0, \\
0 \rar& \SF^\Delta_{p-1}(\delta) \rar& \SF^{\~\Sigma}_p(\uptm\delta) \rar{\prtm_*}& \SF^\Delta_p(\delta) \rar& 0, \\
0 \rar& \SF^\Delta_{p-1}(\delta) \rar& \SF^{\~\Sigma}_p(\basetm\delta) \rar{\prtm_*}& \SF^\Sigma_p(\delta) \rar& 0
\end{tikzcd} \]
where the first map in the two last sequences is given by $\v \mapsto \etm \wedge \prtm^*(\v)$ with $\prtm^*(\v)$ denoting any preimage of $\v$ by $\prtm_*$.
\end{lemma}
Note that the map $\v \mapsto \etm \wedge \prtm^*(\v)$ does not depend on the chosen preimage of $\v$. The proof of this lemma is based on the tropical Deligne sequence and will be given in Section \ref{sec:proof_F_and_tropical_modification}. We now proceed with the proof of the tropical modification formula.

\begin{proof}[Proof of Theorem~\ref{thm:homology_tropical_modification}]
It is possible to derive the theorem by using the above local lemma and by following the proof of Proposition 5.5 in \cite{JRS18}. We present another proof here using cellular homology.

\smallskip
We first prove the isomorphism stated in the theorem concerning open tropical modifications. Let $\~\Sigma = \tropmod{\Delta}{\Sigma}$. Summing the short exact sequences of Lemma \ref{lem:F_and_tropical_modification} over all faces of $\~\Sigma$, we get the following short exact sequence for any pair of integers $p, q$:
\[ \begin{tikzcd}[column sep=small]
  0  \rar&  C^\BM_{p-1,q-1}(\Delta) \oplus C^\BM_{p-1,q}(\Delta)  \rar&  C^\BM_{p,q}(\~\Sigma)  \rar&  C^\BM_{p,q-1}(\Delta) \oplus C^\BM_{p,q}(\Sigma)  \rar&  0.
\end{tikzcd} \]

An inspection of the boundary operators leads to the following short exact sequence of cochain complexes
\begin{equation} \label{eqn:short_exact_sequence_tropical_modification}
\begin{tikzcd}[column sep=small]
  0  \rar&  \Cone_\bul\Bigl(C^\BM_{p-1,\bul}(\Delta) \xrightarrow{\id} C^\BM_{p-1,\bul}(\Delta)\Bigr)  \rar&  C^\BM_{p,\bul}(\~\Sigma)  \rar&  \Cone_\bul\Bigl(C^\BM_{p,\bul}(\Delta) \hookrightarrow C^\BM_{p,\bul}(\Sigma)\Bigr)  \rar&  0,
\end{tikzcd}
\end{equation}
where for two complexes $C^\bul$ and $D^\bul$ and a morphism $\phi\colon C^\bul \to D^\bul$, the complex $\Cone_\bul(\phi)$ is the mapping cone defined as follows
\[ \Cone_k(\phi) = C_{k-1} \oplus D_k,\quad\text{and}\quad \begin{array}{rccc}
\partial\colon & C_{k-1} \oplus D_k & \longrightarrow & C_{k-2} \oplus D_{k-1},          \\
               &    a    \oplus  b  & \longmapsto     & -\d a   \oplus (\phi(a) + \d b).
\end{array} \]
We refer to \cite{KS90} for more information about the mapping cones and their basic properties.

As usual, we define the relative chain complex $C^\BM_{p,\bul}(\Sigma, \Delta)$ as the cokernel of the inclusion $C^\BM_{p,\bul}(\Delta) \hookrightarrow C^\BM_{p,\bul}(\Sigma)$.

The homology of the cone of the identity map is always zero. Moreover, the homology of $\Cone_\bul\Bigl(C_{p,\bul}^\BM(\Delta) \hookrightarrow C_{p,\bul}^\BM(\Sigma)\Bigr)$ is isomorphic to the relative homology.

From the long exact sequence associated to the short exact sequence \eqref{eqn:short_exact_sequence_tropical_modification}, we get the following isomorphism
\[ H^\BM_{\bul,\bul}(\~\Sigma)  \simeq  H^\BM_{\bul,\bul}(\Sigma, \Delta), \]
which is exactly the first isomorphism of the theorem.

\medskip

We now prove the isomorphism between $H_{\bul,\bul}^\BM(\ctropmod{f}{\Sigma})$ and $H_{\bul,\bul}^\BM(\Sigma)$. The chain complex of $\ctropmod{f}{\Sigma}$ can be decomposed into the following short exact sequence
\[ \begin{split}
0 \longrightarrow \Cone_\bul\Bigl(C^\BM_{p-1,\bul}(\Delta) \overset{\id_\Delta}{\longrightarrow} C^\BM_{p-1,\bul}(\Delta)\Bigr) \oplus C^\BM_{p,\bul}(\Delta) \longrightarrow C^\BM_{p,\bul}(\ctropmod{f}{\Sigma}) &\\
& \hspace{-4cm} \longrightarrow \Cone_\bul\Bigl(C^\BM_{p,\bul}(\Delta) \overset{\i_\Delta}\longhookrightarrow C^\BM_{p,\bul}(\Sigma)\Bigr) \longrightarrow 0.
\end{split} \]
The only new term compared to \eqref{eqn:short_exact_sequence_tropical_modification} is the one coming from the stratum $\Delta_\infty := \~\Sigma^\rho_\infty$ (with $\rho$ the special ray) which is isomorphic to $\Delta$. Moreover, when we compare this short exact sequence with the short exact sequence of relative homology, we get the following commutative diagram.
\[ \begin{tikzcd}[column sep = small, row sep = scriptsize]
0 \rar& \Cone_\bul(\id_\Delta) \oplus C_{p,\bul}^\BM(\Delta) \dar{\pi_2}\rar& C_{p,\bul}^\BM(\ctropmod{f}{\Sigma}) \dar{\prtm_*+\i_{\Delta_\infty}}\rar& \Cone_\bul(\i_\Delta) \dar{\pi_1}\rar& 0\\
0 \rar& C_{p,\bul}^\BM(\Delta) \rar& C_{p,\bul}^\BM(\Sigma) \rar& C_{p,\bul}^\BM(\Sigma,\Delta) \rar& 0
\end{tikzcd} \]
Here, $\pi_2$ denotes the projection on the second part, $\prtm_*$ is the usual map coming from the projections $\prtm_*\colon \SF_p(\basetm\sigma) \to \SF^\Sigma_p(\sigma)$ for $\sigma\in\Sigma$, $\i_{\Delta_\infty}$ is given by the maps $\SF_p(\inftm\delta)\simeq\SF^\Delta_p(\delta)\hookrightarrow\SF^\Sigma_p(\delta)$ for $\delta\in\Delta$, and $\pi_1$ is the natural projection
\[ C_{p,\bul}^\BM(\Sigma) \oplus C_{p,\bul-1}^\BM(\Delta) \longrightarrow C_{p,\bul}^\BM(\Sigma) \longrightarrow \rquot{C_{p,\bul}^\BM(\Sigma)}{C_{p,\bul}^\BM(\Delta)}. \]

Hence we get a morphism between the associated long exact sequences.
\[ \begin{tikzcd}[column sep = small, row sep = scriptsize]
\cdots \rar& H^\BM_{p,q}(\Delta) \rar\dar[equal]& H^\BM_{p,q}(\ctropmod{f}{\Sigma}) \rar\dar{\prtm_*+\i_{\Delta_\infty}}& H_q(\Cone_\bul(\i_\Delta)) \rar\dar{\vsim}& H^\BM_{p,q-1}(\Delta) \rar\dar[equal]& \cdots \\
\cdots \rar& H^\BM_{p,q}(\Delta) \rar& H^\BM_{p,q}(\Sigma) \rar& H^\BM_{p,q}(\Sigma,\Delta) \rar& H^\BM_{p,q-1}(\Delta) \rar& \cdots
\end{tikzcd} \]
The maps $\pi_2$ and $\pi_1$ induce isomorphisms in homology, and we conclude by using the five lemma.

\smallskip
The other isomorphisms in the statement of the theorem can be obtained in similar ways (using if necessary, for the cohomology, the dual of the local tropical modification formula, Lemma \ref{lem:F_and_tropical_modification}). The theorem follows.
\end{proof}

\subsection{Shellability of smoothness}

As a consequence of Theorem \ref{thm:homology_tropical_modification}, we get the following.

\begin{thm} \label{thm:smooth_shellable}
Smoothness is shellable in $\tropf$.
\end{thm}

\begin{proof}
We use Lemma \ref{lem:shellability_meta_lemma}. It is a trivial fact that elements of $\Bsh_0$ are smooth. Moreover, since smoothness only depends on the support (Theorem \ref{thm:smoothness-support}), we only have to prove the closeness by products and the closeness by tropical modifications. Closeness by products is given by Proposition \ref{prop:smoothness_product}. It remains to prove the closeness by tropical modifications.

\smallskip
Let $\Sigma$ be a tropical fan. Let $f$ be a conewise integral linear function on $\Sigma$ such that $\div(f)$ is reduced. Set $\Delta = \div(f)$ and let $\~\Sigma = \tropmod{f}{\Sigma}$. Assume moreover that $\Delta$ is smooth (and non-empty). By Lemma \ref{lem:shellability_meta_lemma}, we just have to prove the Poincaré duality for $\~\Sigma$.

By Theorem \ref{thm:homology_tropical_modification}, we know that $H_{\bul,\bul}^\BM(\~\Sigma) \simeq H^\BM_{\bul,\bul}(\Sigma,\Delta)$. Using the long exact sequence of relative homology, we get the following exact sequence
\[ \dots \longrightarrow H^\BM_{p,q}(\Delta) \longrightarrow H^\BM_{p,q}(\Sigma) \longrightarrow H^\BM_{p,q}(\~\Sigma) \longrightarrow H^\BM_{p,q-1}(\Delta) \longrightarrow \cdots \]
Since $\Sigma$ (\resp $\Delta$) verifies the Poincaré duality, its Borel-Moore homology is trivial except for $q=d$ (\resp $q=d-1$). We deduce that $H^\BM_{p,\bul}(\~\Sigma)$ is trivial except in degree $d$ and that we have a short exact sequence
\[ 0 \longrightarrow H^\BM_{p,d}(\Sigma) \longrightarrow H^\BM_{p,d}(\~\Sigma) \longrightarrow H^\BM_{p,d-1}(\Delta) \longrightarrow 0. \]

The cap products $\cdot \frown \nu_\varUpsilon$ for $\varUpsilon$ any of the tropical fans $\Sigma$, $\Delta$ or $\~\Sigma$ induces the following commutative diagram.
\[ \begin{tikzcd}
0 \rar& \SF^{d-p}(\Sigma) \dar{\vsim}\rar& \SF^{d-p}(\~\Sigma) \dar  \rar& \SF^{d-p-1}(\Delta) \dar{\vsim}\rar& 0 \\
0 \rar& H^\BM_{p,d}(\Sigma)     \rar& H^\BM_{p,d}(\~\Sigma)     \rar& H^\BM_{p,d-1}(\Delta)     \rar& 0
\end{tikzcd} \]
We have already seen that the second row is exact. The first row is also a short exact sequence by the dual of Lemma \ref{lem:F_and_tropical_modification}. The first and last vertical maps are isomorphisms because $\Sigma$ and $\Delta$ verify the Poincaré duality. Hence, by the five lemma, the second vertical map is an isomorphism, and therefore, $\~\Sigma$ verifies the Poincaré duality.
\end{proof}

\subsection{Proof of Lemma \ref{lem:F_and_tropical_modification}} \label{sec:proof_F_and_tropical_modification}

In this final section, we prove Lemma \ref{lem:F_and_tropical_modification}. So let $\Sigma$ be a smooth tropical fan. Let $f$ be a conewise integral linear function on $\Sigma$ such that $\div(f)$ is reduced. Let $\Delta = \div(f)$ and $\~\Sigma = \tropmod{f}{\Sigma}$. We have to prove the following three exact sequences for any face $\sigma\in\Sigma\setminus\Delta$ and any $\delta\in\Delta$.
\[ \begin{tikzcd}[row sep = tiny]
& 0 \rar& \SF^{\~\Sigma}_p(\basetm\sigma) \rar{\prtm_*}& \SF^\Sigma_p(\sigma) \rar& 0, \\
0 \rar& \SF^\Delta_{p-1}(\delta) \rar& \SF^{\~\Sigma}_p(\uptm\delta) \rar{\prtm_*}& \SF^\Delta_p(\delta) \rar& 0, \\
0 \rar& \SF^\Delta_{p-1}(\delta) \rar& \SF^{\~\Sigma}_p(\basetm\delta) \rar{\prtm_*}& \SF^\Sigma_p(\delta) \rar& 0.
\end{tikzcd} \]

For the first exact sequence, let $\sigma \in \Sigma \setminus \Delta$. Let $\ell \in M$ be any integral linear map which coincides with $f$ on $\sigma$. Denote by $\pi^\sigma\colon N \to N^\sigma$ the natural projection. Then, $f-\ell$ induces a conewise integral linear map $f^\sigma = \pi^\sigma_*(f - \ell)$ on $\Sigma^\sigma$. Since no codimension one face of $\Sigma$ containing $\sigma$ is in $\div(f)$, we deduce that $\div(f^\sigma)$ is trivial. Since $\Sigma$ is smooth, it is div-faithful at $\sigma$. As a consequence, $f^\sigma$ must be a linear map on $\Sigma^\sigma$. Set $\~\ell = \pi^{\sigma,*}(f^\sigma) + \ell$, and note that $\~\ell$ is a linear map which coincides with $f$ on the faces containing $\sigma$. The isomorphism between $\SF^\Sigma_p(\sigma)$ and $\SF^{\~\Sigma}_p(\basetm\sigma)$ is now induced by the isomorphism between $N_\R$ and $\Im(\~\ell)$.

\smallskip
The second exact sequence is clear since $\uptm{\eta\/} \simeq \eta \times \R_+\etm$ for any cone $\eta \supface \delta$ in $\Delta$.

\smallskip
It remains to prove the last one. We prove the exactness of the following sequence
\[ 0 \longrightarrow \SF^\Delta_{p-1}(\conezero_\Delta) \longrightarrow \SF^{\~\Sigma}_p(\conezero_{\~\Sigma}) \longrightarrow \SF^\Sigma_p(\conezero_\Sigma) \longrightarrow 0.\]
The statement then follows by applying the same argument to $\Delta^\delta$, $\~\Sigma^{\basetm\delta}$, and $\Sigma^\delta$ for all $\delta\in\Delta$, and by using the compatible isomorphisms
\[ \SF^\varUpsilon_\bul(\mu) \simeq \bigwedge^\bul N_\mu \,\otimes\, \SF^{\varUpsilon^\mu}_\bul(\conezero_{\varUpsilon^\mu}) \]
for pairs $(\mu,\varUpsilon)$ equal to $(\delta,\Delta), (\basetm\delta,\~\Sigma)$, and $(\delta,\Sigma)$. Notice that the terms of the short exact sequence only depend on the support of the respective fan. Hence, we can assume without loss of generality that $\Delta$, $\Sigma$ and $\~\Sigma$ are unimodular.

\smallskip
Combining the sequence given by Proposition~\ref{prop:non_smooth_Deligne} for $\Delta, \Sigma$ and $\~\Sigma$, we get the diagram depicted in Figure \ref{fig:dual_Deligne_tropical_modification} in which the three rows are exact.
\begin{figure}[ht]
\[ \begin{tikzcd}[cells={font=\everymath\expandafter{\the\everymath\displaystyle}}, row sep=scriptsize, column sep=scriptsize]
0 \dar                                                                                 & 0 \dar \\
\displaystyle\bigoplus_{\delta \in \Delta \\ \dims\delta=d-p-2}\!\!H^{p,p}(\comp\Delta^\delta)   \rar\dar& \bigoplus_{\delta \in \Delta \\ \dims\delta=d-p-1} H^{p,p}(\comp\Delta^\delta)     \rar\dar& H_c^{p,d-1}(\Delta) \rar\dar& 0 \\
\bigoplus_{\sigma \in \~\Sigma \\ \dims\sigma=d-p-1}\!\!H^{p,p}(\comp{{\~\Sigma}}^\sigma) \rar\dar& \bigoplus_{\sigma \in \~\Sigma \\ \dims\sigma=d-p} H^{p,p}(\comp{{\~\Sigma}}^\sigma) \rar\dar& H_c^{p,d}(\~\Sigma) \rar\dar& 0 \\
\bigoplus_{\sigma \in \Sigma \\ \dims\sigma=d-p-1}\!\!H^{p,p}(\comp\Sigma^\sigma)   \rar\dar& \bigoplus_{\sigma \in \Sigma \\ \dims\sigma=d-p} H^{p,p}(\comp\Sigma^\sigma)       \rar\dar& H_c^{p,d}(\Sigma)   \rar\dar& 0 \\
0                                                                                      & 0                                                                                               & 0
\end{tikzcd} \]
\caption{End of the duals of tropical Deligne sequences for a tropical modification \label{fig:dual_Deligne_tropical_modification}}
\end{figure}

We describe the vertical maps. Let $\delta \in \Delta$. Then $\~\Sigma^{\uptm\delta} \simeq \Delta^\delta$. In the same way, if $\sigma \in \Sigma\setminus\Delta$, then $\~\Sigma^{\basetm\sigma} \simeq \Sigma^\sigma$ because $\Sigma$ is div-faithful at $\sigma$. In both cases, we get an isomorphism in cohomology. The last case concerns $\~\Sigma^\sigma$ with $\sigma \in \Delta$. In this case, $\~\Sigma^{\basetm\sigma} \simeq \tropmod{\Delta^\sigma}{\Sigma^\sigma}$. Since $\Sigma^\sigma$ is div-faithful, we can combine Theorem~\ref{thm:invariant_chow_tropical_modification} and Theorem \ref{thm:Hodge_conjecture} to get
\[ H^{p,p}(\comp{{\~\Sigma}}^{{\basetm\sigma}}) \simeq \MW_{d-p}(\~\Sigma^{\basetm\sigma}) \simeq A^{d-p}(\~\Sigma^{\basetm\sigma})^\dual \simeq A^{d-p}(\Sigma^\sigma)^\dual \simeq H^{p,p}(\comp\Sigma^\sigma). \]

Summing all these isomorphisms, we get
\[ \bigoplus_{\sigma\in\~\Sigma_k} H^{p,p}(\comp{{\~\Sigma}}^\sigma) = \bigoplus_{\delta \in \Delta_{k-1}} H^{p,p}(\comp{{\~\Sigma}}^{\uptm\delta}) \oplus \bigoplus_{\sigma \in \Sigma_k} H^{p,p}(\comp{{\~\Sigma}}^{\basetm\sigma}) \simeq \bigoplus_{\delta \in \Delta_{k-1}} H^{p,p}(\comp{\Delta}^{\delta}) \oplus \bigoplus_{\sigma \in \Sigma_k} H^{p,p}(\comp{\Sigma}^\sigma). \]
The two first columns of Figure \ref{fig:dual_Deligne_tropical_modification} are the split exact sequences obtained via the isomorphism described above. The last column can be defined in a unique way to make the diagram commutative. A diagram chasing proves that we obtain an exact sequence.

\smallskip
We now observe that we have a second commutative diagram.
\[ \begin{tikzcd}
      & H_c^{p,d-1}(\Delta)          \rar\dar[two heads]& H_c^{p,d}(\~\Sigma)                       \rar\dar[two heads]& H_c^{p,d}(\Sigma)        \rar\dar{\vsim}& 0 \\
0 \rar& \SF^\Delta_{d-p-1}(\conezero_\Delta) \rar           & \SF^{\~\Sigma}_{d-p}(\conezero_{\~\Sigma}) \rar{\prtm_*}           & \SF^\Sigma_{d-p}(\conezero_\Sigma) \rar& 0
\end{tikzcd} \]
The vertical maps are the dual of the map described in \ref{subsubsec:reformulation_PD}. In the same way that the map in \ref{subsubsec:reformulation_PD} is injective, the vertical maps are surjective. Moreover, the last vertical map is an isomorphism. Indeed, by assumption $\Sigma$ is smooth and thus verifies the Poincaré duality. Hence, we have an isomorphism $\SF^{d-p}_\Sigma(\conezero_\Sigma) \simto H^\BM_{p,d}(\Sigma)$. Moreover, since $H^{p,\bul}(\Sigma)$ is torsion-free, the Poincaré duality implies that $H^\BM_{p,\bul}(\Sigma)$ is torsion-free, and we have $H_c^{p,d}(\Sigma) \simeq H^\BM_{p,d}(\Sigma)^\dual$. Since $\SF_{d-p}^\Sigma(\conezero_\Sigma) \simeq \SF^{d-p}_\Sigma(\conezero_\Sigma)^\dual$, by dualizing, we get the isomorphism $H_c^{p,d}(\Sigma) \simto \SF_{d-p}^\Sigma(\conezero_\Sigma)$. The injectivity of the map $\SF^\Delta_{d-p-1}(\conezero_\Delta) \to \SF_{d-p}^{\~\Sigma}(\conezero_{\~\Sigma})$ is clear. The rest of the exactness of the second row then follows from a diagram chasing. This concludes the proof of Lemma~\ref{lem:F_and_tropical_modification}. \qed

\section{Examples and further questions}\label{sec:examples}

In this final section, we provide a collection of examples to which we referred in the text, and complement this with remarks and questions which we hope could clarify the concepts introduced in the paper.

\subsection{The fan over the one-skeleton of the cube} \label{subsec:cube} A rich source of examples is given by the fan defined over the $1$-skeleton of a cube. More precisely, consider the standard cube $\mbox{\mancube}$ with vertices $(\pm1, \pm1, \pm1)$, and let $\Sigma$ be the two-dimensional fan with rays generated by vertices and with facets generated by edges of the cube. The fan $\Sigma$ is locally irreducible and tropical but not unimodular. We obtain a unimodular fan changing the lattice: in what follows, we work with the lattice $N := \sum_{\varrho\in\Sigma_1}\Z\e_\varrho$.

Table \ref{tab:cohomology_cube} summarizes the main cohomological data about the cube. One can compute the homology via the universal coefficient theorem: the torsion part is shifted by one column on the left, the rest remains unchanged. Moreover the image of $A^1(\Sigma)$ inside $A^1(\Sigma)^\dual$ is a sublattice of full rank and of index two.

\begin{table}
\renewcommand{\strut}{\rule{0pt}{1.1em}}
\[ \begin{array}{c|c|c|c|cc|c|c|c|}
\cline{2-4} \cline{7-9}
\strut H^{p,q}_c(\Sigma)       & 0 & 1    & 2                         & \qquad\qquad          & H^{p,q}(\comp\Sigma) & 0  & 1    & 2    \\ \cline{1-4} \cline{6-9}
\multicolumn{1}{|c|}{\strut 0} & 0 & 0    & \Z^5                      & \multicolumn{1}{c|}{} & 0                    & \Z & 0    & 0    \\ \cline{1-4} \cline{6-9}
\multicolumn{1}{|c|}{\strut 1} & 0 & 0    & \Z^3\times\rquot{\Z}{2\Z} & \multicolumn{1}{c|}{} & 1                    & 0  & \Z^5 & 0    \\ \cline{1-4} \cline{6-9}
\multicolumn{1}{|c|}{\strut 2} & 0 & \Z^2 & \Z                        & \multicolumn{1}{c|}{} & 2                    & 0  & \Z^2 & \Z   \\ \cline{1-4} \cline{6-9}
\end{array} \]
\caption{Cohomology of the fan $\Sigma$ over the one-skeleton of the cube. \label{tab:cohomology_cube}}
\end{table}

\subsubsection*{Some general facts on the cohomology of fans} Let us start by making comments on some relations in these tables that generalize to all fans. First notice that the last rows coincide. This is a more general phenomenon: for any fan $\Sigma'$ of dimension $d$, $H^{d,q}_c(\Sigma') \simeq H^{d,q}(\comp\Sigma')$. This is due to the vanishing of the sheaf $\SF^d$ on faces of positive sedentarity. Moreover, $H^{p,0}_c(\Sigma)$ is trivial and $H_c^{p,1}(\Sigma)$ is torsion-free for any $p$. Indeed, by definition of $\SF_p(\conezero)$, the map $\bigoplus_{\varrho\in\Sigma_1}\SF_p(\varrho) \to \SF_p(\conezero)$ is surjective. Hence $H^\BM_{p,0}(\Sigma)$ is trivial and the universal coefficient theorem concludes.

\subsubsection*{Poincaré duality for fans and classical cohomology}
Pursuing our study of the cohomology of $\Sigma$, one can notice that the first row is just the reduced classical homology of the one-skeleton of the cube shifted by one. This is clear because for $p=0$ we retrieve the classical homology. This gives us a precious information: a fan $\Sigma'$ of dimension $d$ which verifies the Poincaré duality must verify that $\Sigma'\setminus\{0\}$ has its classical cohomology concentrated in degree $d-1$, and moreover, $H_{\textrm{sing}}^{d-1}(\Sigma') \simeq \Z^{\rk(\SF^d(\conezero))}$. Combined with the fact that $\rk(\SF^d(\conezero)) \leq \binom{n}{d}$, with $n$ the dimension of the space spanned by $\Sigma$, we deduce that locally irreducible tropical fans verifying the Poincaré duality are very specific.

\subsubsection*{On cohomology groups of $\comp\Sigma$ for $p>q>0$}
Notice that by Theorem \ref{thm:ring_morphism}, the cohomology groups $H^{p,q}(\comp\Sigma)$ are trivial when either $p<q$, or $q=0$ and $p>0$, and that $H^{p,p}(\comp\Sigma)\simeq A^p(\Sigma)$. The theorem however provides no information about $H^{2,1}(\comp\Sigma)$ for a general fan. For the fan $\Sigma$ consider here, this cohomology group is nontrivial, hence $\comp\Sigma$ does not verify the Poincaré duality.

\subsubsection*{Other properties of\/ $\Sigma$} The fan $\Sigma$ can be served in the following examples.

\begin{example}[A Chow-duality tropical fan which does not verify the Poincaré duality] \label{ex:cube}
The cohomology of $\Sigma$ does not verify the Poincaré duality. Its Chow ring with integral coefficients neither. On the other hand, its Chow ring with rational coefficients does verify the Poincaré duality. It also verifies the hard-Lefschetz property and the Hodge-Riemann bilinear relations.
\end{example}

\begin{example}[A non-principal locally irreducible unimodular fan]
The map $A^1(\Sigma) \to A^1(\Sigma)^\dual$ is not surjective, hence $\Sigma$ is not principal at $\conezero$ for integral coefficients (though it is principal for rational coefficients). Indeed, the divisor $D=\R\cdot(1,1,1)$ is not principal, but $2D$ is. In fact, this divisor is not a boundary for integral coefficients, and $H_{1,1}^\BM(\Sigma) \simeq \rquot{\Z}{2\Z}$ has torsion.
\end{example}

\subsection{The cross}

Consider the cross $\{xy = 0\}$ in $\R^2$ with four rays. We denote it by $\Delta$. It is the simplest singular tropical fan. It can be used in several constructions to find various counter-examples. We propose some of them below.

\begin{example}[Poincaré duality for a non normal compactified fan] \label{ex:PD_chow_not_smooth}
Recall that $\Lambda^k$ is the complete fan on $\R^k$ whose facets are the $2^k$ orthants. Consider $\Sigma:=\tropmod{\Delta}{\Lambda^2}$. Since $\Lambda^2$ is div-faithful, by Theorem \ref{thm:invariant_chow_tropical_modification}, $A^\bul(\Sigma) \simeq A^\bul(\Lambda^2)$. Moreover, since $\Lambda^2$ is smooth, $H^{\bul,\bul}(\comp\Sigma) \simeq H^{\bul,\bul}(\comp{\Lambda^2})$. Hence, despite $\Sigma$ being not normal, both the Chow ring of $\Sigma$ and the cohomology of $\comp\Sigma$ verify the Poincaré duality. In particular, $\Sigma$ is irreducible, principal and div-faithful at $\conezero$. On the other hand, it is neither locally irreducible nor div-faithful globally. However, $\Sigma$ does not verify the homological Poincaré duality: $H^\BM_{0,2}(\Sigma) \simeq \Z^4$ but $\SF^2(\conezero)$ is of rank $3$.
\end{example}

\begin{example}[A nontrivial degenerate tropical modification] \label{ex:nontrivial_degenerate_tropical_modification}
Let $f$ be the conewise integral linear function on $\Delta$ given by $f=x-y$ on $\{x + y \geq 0\}$ and $0$ on the complement. Then $\div(f)$ is trivial. The tropical modification $\tropmod{f}{\Delta}$ is a smooth tropical line. In particular it verifies the Poincaré duality, which is not the case of $\Delta$. This example is quite interesting: it seems that tropical modifications tend to \emph{desingularize} tropical varieties as blow-ups do in algebraic geometry.
\end{example}

\begin{example}[Two non-principal unimodular tropical fans] \label{ex:non-principal_fan}
Consider the two-dimensional fan $\Sigma=\Delta\times\Lambda$ in $\R^3$, \ie, it has the six rays of the axes and support $\{x = 0\} \cup \{y = 0\} \subset \R^3$. The divisor $\{y = z = 0 \}$ in $\Sigma$ is not principal. This means $\Sigma$ is not principal at $\conezero$.

Following the idea of Example \ref{ex:PD_chow_not_smooth}, one can go further and create a fan which is principal at $\conezero$ but not globally principal: it suffices to take $\tropmod{\Sigma}{\Lambda^3}$.
\end{example}

\begin{example}[Non-irreducible unimodular fans]
We can find higher dimensional analogues of the cross in $\R^2$. Consider any unimodular fan $\Delta'$ of support $\{x_1 = x_2 = 0\} \cup \{x_3 = x_4 = 0\}$ in $\R^4$. The fan $\Sigma$ is normal, div-faithful and principal, but it is not irreducible.
\end{example}

\begin{example}[A normal div-faithful unimodular fan which is neither irreducible nor principal] \label{ex:div-faithful_not_irreducible}
Let $\Delta'$ be the fan defined in the previous example. Set $\Sigma = \Delta'\times\Lambda$. We get a tropical fan $\Sigma$ of dimension three which is still normal and div-faithful but it is neither irreducible, nor principal: the divisor $\{x_1 = x_2 = 0\}\times\{0\}$ is not a principal divisor in $\Sigma'$.
\end{example}

\subsection{Non saturated and non unimodular fans} \label{subsec:examples_non_saturated_non_unimodular}

Let $(\e_1,\e_2)$ be the standard basis of $\Z^2$ and let $\e_0=-\e_1-\e_2$. Denote by $\rho_i = \R_+\e_i$, $i\in\{0,1,2\}$ the corresponding rays. Let $N$ be the lattice $\Z(1,0)+\Z(-1/3,-1/3)$. Then the lattice $\Z^2$ generated by $\e_0,\e_1,\e_2$ is of index three in $N$.

\begin{example} \label{ex:F1_vs_N}
Let $\Sigma$ in $\R^2$ be the one-dimensional fan with rays $\rho_i$, for $i=0,1,2$. Then $\Sigma$ is tropical, and tropically smooth, but it is not saturated. Any element of $(\Z^2)^\dual \subset M$ induces an integral linear form on $\Sigma$. As a consequence, $A^1(\Sigma)$ has torsion. For instance, $\x_{\rho_1}-\x_{\rho_2}$ is non-zero, but $3\cdot(\x_{\rho_1}-\x_{\rho_2})$ is zero in $A^1(\Sigma)$.
\end{example}

\begin{example} \label{ex:non_unimodular_projective_fan}
Let $\Sigma$ be the complete fan in the plane $\R^2$, with $N =\Z^2$, whose rays are the $\rho_i$. This time, $\Sigma$ is not unimodular. As in example \ref{ex:F1_vs_N}, $A^1(\Sigma)$ has torsion. Moreover, the cohomology group $H^{1,2}(\comp\Sigma) \simeq \rquot\Z{3\Z}$ is nontrivial even though we are in bidegree $(p,q) = (1,2)$ with $p < q$.
\end{example}

\begin{example} \label{ex:non_saturated_smooth_fan}
Consider again the complete fan $\Sigma$ of the previous Example. Let $f$ be the conewise integral linear function on $\Sigma$ which maps $\e_0$ on $-3$ and $\e_1$ and $\e_2$ on $0$. Then $\div(f)$ is the reduced divisor $\Delta$ with rays $\rho_i$ for $i=0,1,2$.

Let $\Sigma'$ be a unimodular subdivision of $\Sigma$; for instance one can add the rays along the vectors $\tfrac13\e_i+\tfrac23\e_j$ for distinct $i,j \in \{0,1,2\}$. This fan is smooth. Set $\~\Sigma = \tropmod{f}{\Sigma'}$ and denote by $\rho$ the new ray. Then $\~\Sigma$ is a smooth unimodular fan which is saturated at $\conezero$ but not at $\rho$. Even worse, there is no way to change the lattice $N$ to make it saturated without changing $\~\Sigma \cap N$.

Since $\~\Sigma$ is saturated at $\conezero$, $A^1(\~\Sigma)$ has no torsion. However, the element $\x_\rho(\x_{\rho_1}-\x_{\rho_2}) \in A^2(\~\Sigma)$ is an element of order three.
\end{example}

The three previous examples show that the saturation and unimodularity assumptions in Theorem~\ref{thm:ring_morphism} are needed in general.

\subsection{Further examples}

\begin{example}[A shellable tropical fan which is not generalized Bergman] \label{ex:shellable_not_Bergman}
Let $\Sigma$ be the two-dimensional skeleton of the fan of the projective space of dimension three: it has four rays generated by vectors $\e_1, \e_2, \e_3$ and $\e_0= -\e_1-\e_2-\e_3$, where $(\e_1,\e_2,\e_3)$ is a basis of $N$, and the six facets $\R_+\e_i+\R_+e_j$ for $0\leq i<j\leq 3$. We have $\supp\Sigma = \supp{\Sigma_\Ma}$ for the uniform matroid $\Ma = U^3_{4}$ of rank three on four elements. Let $\Delta$ be the tropical curve in $\supp\Sigma$ with four rays $\rho_1, \rho_2, \rho_3$ and $\rho_4$ with
\[\rho_1 =\R_+ (\e_1+\e_2), \qquad \rho_2 =\R_+ \e_2, \qquad \rho_3 = \R_+(\e_0+ 2 \e_3), \qquad \rho_4 = \R_+ (\e_0 + \e_1).\]
Then we have $\Delta \simeq \Sigma_{\Ma'}$ with $\Ma' = U^2_4$ the uniform matroid of rank two on four elements. The tropical modification $\tropmod\Delta\Sigma$ is shellable but does not have the same support as the Bergman fan of any matroid. It means it is not generalized Bergman.
\end{example}

\begin{example}[A non-shellable smooth fan]
We study here an interesting example discovered by Babaee and Huh \cite{BH17}. We would like to thank Edvard Aksnes and Kristin Shaw for pointing out the relevance of this example to us.

A full description of the fan is given in \cites{BH17, Aks19}. We just mention its main properties. It is a normal unimodular tropical fan $\Sigma$ of dimension two living in $\R^4$. Aksnes computes its rational homology in \cite{Aks19} and concludes that it is smooth (this remains true with integral coefficients). However, $\Sigma$ does not verify the Hodge-Riemann bilinear relations: it is quasi-projective, but the pairing $A^1(\Sigma) \times A^1(\Sigma) \to \Z$ has more than one positive eigenvalue. In particular it is not shellable. Indeed, by Theorem \ref{thm:examples_shellable} \ref{thm:examples_shellable:HR_Chow_ring}, any quasi-projective shellable fan verifies the Hodge-Riemann bilinear relations.
\end{example}

\begin{example}[Irreducible components of non normal fans] \label{ex:irreducible_components_not_well_defined}
Let $\Sigma = \tikz[scale=.18]\draw (0:-1)--(0:1) (45:-1.2)--(45:1.2) (90:-1)--(90:1);$ be the $1$-skeleton of $U_{3,3}$ in $\R^2$. There are five irreducible components in the sense of minimal support of an element of $H^\BM_{1,1}(\Sigma)$, but they do not induce a partition of $\Sigma_1$.
\end{example}

\subsection{Questions} We end this section with a list of questions naturally arising in our study.

\smallskip
\noindent\ --\, Is the smoothness of a tropical fan $\Sigma$ implied by the surjectivity of the map $\SF^p(\conezero_{\Sigma^\sigma}) \to H^{\BM}_{d-p, d}(\Sigma^\sigma)$, $\sigma \in \Sigma$. In other words, does the vanishing of $H_{p,q}^{\BM}(\Sigma)$ follow from the surjectivity?

\smallskip
\noindent\ --\, Provide examples of unimodular fans $\Sigma$ with torsion in $H^{p,q}(\comp \Sigma)$, for $p>q$.

\smallskip
\noindent\ --\, Provide a classification of smooth tropical fans.

\let\v\cech
\bibliographystyle{alpha}
\bibliography{$HOME/bibliography/bibliography}

\end{document}